\documentclass{amsart}
\usepackage{graphicx}
\usepackage{color}
\usepackage[emtex,matrix,arrow]{xy}
\usepackage{latexsym}
\usepackage{epsfig}
\newtheorem{thm}{Theorem}[section]
\newtheorem{cor}[thm]{Corollary}
\newtheorem{lem}[thm]{Lemma}
\newtheorem{prop}[thm]{Proposition}
\theoremstyle{definition}
\newtheorem{defn}[thm]{Definition}
\theoremstyle{remark}
\newtheorem{rem}[thm]{Remark}
\theoremstyle{remark}
\newtheorem{ex}[thm]{Example}
\numberwithin{equation}{section}
\newcommand{\To}{\longrightarrow}

\newcommand{\C}{\mathcal{C}}
\newcommand{\D}{\mathcal{D}}

\newcommand{\F}{\mathcal{F}}
\newcommand{\G}{\mathcal{G}}
\newcommand{\Hl}{\mathcal{H}}

\newcommand{\B}{\mathcal{B}}

\font\gotica= eufm10
\font\Gotica= eufm10

\def\Bg{{\Gotica\char 66}}
\def\Cg{{\Gotica\char 67}}
\def\Dg{{\Gotica\char 68}}

\def\gg{{\gotica\char 103}}

\def\mm{{\gotica\char 109}}
\def\nn{{\gotica\char 110}}
\begin{document}

\title[]{Cohomology and deformation theory of monoidal 2-categories I}
\author{Josep Elgueta}
\address{Dept. Matem\`{a}tica Aplicada II, Universitat Polit\`{e}cnica de Catalunya}
\email{jelgueta@ma2.upc.es}

\thanks{}
\subjclass{}
\keywords{}
\date{March 2002}
\dedicatory{}
\commby{}

%%% ----------------------------------------------------------------------
\begin{abstract}
In this paper we define a cohomology for an arbitrary $K$-linear
semistrict semigroupal 2-category $(\text\Cg,\otimes)$ (called for
short a Gray semigroup) and show that its first order (unitary)
deformations, up to the suitable notion of equivalence, are in
bijection with the elements of the second cohomology group.
Fundamental to the construction is a double complex, similar to the
Gerstenhaber-Schack double complex for bialgebras, the role of the
multiplication and the comultiplication being now played by the
composition and the tensor product of 1-morphisms. We also identify
the cohomologies describing separately the deformations of the
tensor product, the associator and the pentagonator. To obtain the above results, a
cohomology theory for an arbitrary $K$-linear (unitary)
pseudofunctor is introduced describing its purely
pseudofunctorial deformations, and generalizing Yetter's cohomology
for semigroupal functors \cite{dY98}. The corresponding higher
order obstructions will be considered in detail in a future paper.
\end{abstract}

%%% ----------------------------------------------------------------------
\maketitle
%%% ----------------------------------------------------------------------

\section{Introduction}

%{\bf To rewrite}

This is the first of two papers where we intend to give a
cohomological description of the infinitesimal deformations of a
monoidal 2-category.

Roughly speaking, a monoidal 2-category is
a 2-category equipped with a binary operation, usually called the
{\it tensor product}, defined at the three levels existing in any 2-category, i.e., objects,
1-morphisms and 2-morphisms, and which is associative and with a unit
up to suitable 2-isomorphisms. Actually, in this paper we will
consider the more general structure of a semigroupal 2-category
\Cg, namely, a 2-category with a tensor product as above but
which is only associative
(up to a suitable 2-isomorphism) with no unit. More explicitly,
we show that the first order (unitary) deformations of such an
object
\Cg\ can be identified with the elements of some cohomology group
associated to \Cg. The generalization to the case of monoidal
categories and the question of the obstructions will be treated in
a future paper.

This work is an extension to the context of 2-categories of the
theory developped by Crane and Yetter for semigroupal categories
\cite{CY981} and by Yetter for braided monoidal categories
\cite{dY98},\cite{dY01}, which are in turn an extension to the
context of (1-)categories of Gerstenhaber's work on deformations of
algebras \cite{mG63},\cite{mG64}, later generalized to the case of
Hopf algebras by Gerstenhaber and Schack \cite{GS90} (see also
\cite{GS88},\cite{GS92}). These classical works should be viewed as
the corresponding theories in the so-called {\it 0-dimensional
algebra} setting \cite{BD95}, which is the algebra in the context
of sets. The situation can be schematically represented as in the
table below. This table is the $K$-linear version of the first two
rows in the table of {\it k-tuply monoidal n-categories} of Baez
and Dolan; see \cite{BD95}, Table I. The $n$ here denotes the
``dimensionality'' of the algebraic framework we work with. So,
dimension $n$ corresponds to work in the context of {\it
n-categories}, a natural generalization of the notion of 2-category
where we also have 3-morphisms between the 2-morphisms, and so on,
until $n$-morphisms between $(n-1)$-morphisms.

\vspace{0.5 truecm}
\begin{tabular}{|c|c|c|c|} \hline
   & n=0 & n=1 & n=2 \\ \hline
k=0 & $\begin{array}{c} K\text{-vector} \\ \text{space}
\end{array}$ & $\begin{array}{c} K\text{-linear} \\ \text{category} \end{array}$ &
$\begin{array}{c} K\text{-linear} \\ \text{2-category} \end{array}$
\\ \hline k=1 & $\begin{array}{c} K\text{-algebra} \\
{\mathbf (Gerstenhaber's\ work)} \end{array}$ & $\begin{array}{c}
K\text{-linear} \\ \text{monoidal category} \\ {\mathbf
(Crane\text{-}Yetter's\ work)}
\end{array}$ & $\begin{array}{c} K\text{-linear} \\ \text{monoidal 2-category} \end{array}$ \\
\hline
\end{tabular}

\vspace{0.5 truecm}
 Notice that going from the top to the
bottom row along the diagonal corresponds to taking the one object
case (for example, a $K$-linear category of only one object is
exactly the same thing as a $K$-algebra). For a more expanded
explanation of this table, see the reference above.

This idea of generalizing Gerstenhaber's work to 1- and
2-categories comes from the important notion of {\it
categorification}, which in the table corresponds to moving to the
right. It first appears in the work by Crane and Frenkel
\cite{CF94} on Hopf categories, although it seems it was already
present in Grothendieck's thought. A Hopf category is an analog of
a Hopf algebra where the underlying $K$-vector space has been
substituted by a special kind of $K$-linear category, usually
called a 2-vector space over $K$. For a more precise definition,
the reader is refered to Neuchl's thesis \cite{mN97}. The basic
idea behind the notion of categorification is that constructions
from usual algebra can be translated to the level of categories and
to higher levels, the {\it n-categories}, for $n\geq 1$, making one
essential step: to substitute equations for isomorphisms. The price
to pay is that it is necessary to simultaneously impose equations
on these isomorphisms, which are the so-called {\it coherence
relations}. These ideas have been developed by different authors,
in particular, Crane and Yetter \cite{CY982} and in the more
general language of n-categories, Baez and Dolan (see, for ex.,
\cite{BD98}).

Apart from its own interest, our motivation for studying
deformations of monoidal 2-categories has to be found in its
potential applications to the construction of interesting four
dimensional Topological Quantum Field Theories (briefly, TQFT's).
Indeed, in \cite{CF94}, Crane and Frenkel suggest that Hopf
categories may be used to construct four dimensional TQFT's, in a
similar way as three dimensional TQFT's can be constructed from
Hopf algebras (see, for ex., \cite{gK91},\cite{CFS94}). Now, it is
well-known that three dimensional TQFT's can also be obtained using
the monoidal categories of representations of suitable Hopf
algebras (see, for ex., \cite{TV92}, \cite{dY94},\cite{BW96}). This
clearly suggests the possibility that four dimensional TQFT's could
be obtained from the category of representations of suitable Hopf
categories, which will be, going up in the categorification
process, some kind of monoidal 2-categories (actually, Neuchl
\cite{mN97} has proved that the 2-category of representations of a
Hopf category is indeed monoidal). That idea has been made explicit
by Mackaay \cite{mM99}, who develops a method to construct
invariants of piecewise linear four manifolds from a special kind
of monoidal 2-categories he calls {\it spherical 2-categories}. His
construction parallels that of Barrett and Westbury \cite{BW96} for
three manifolds. This explains the interest of monoidal
2-categories in the construction of four dimensional TQFT's. But,
why are we interested in their deformations? The answer is again an
expected analogy between the cases of dimension three and four. In
dimension three, we can get a state sum invariant of a piecewise
linear three-manifold using irreducible representations of an
arbitrary semisimple Lie algebra. The method comes from the
classical work of Ponzano-Regge \cite{PR68}. The problem is that
the sum turns out to be infinite. Progress was made possible only
when the corresponding quantum group was discovered, which is a
deformation (as a {\it braided bialgebra}; see \cite{cK95}) of the
classical universal envelopping algebra of the Lie algebra. Using
the representations of the quantum group at a root of unity instead
of those of the classical version, the state sum invariant becomes
a convergent sum. That's what Turaev and Viro do in their paper
\cite{TV92}. The reader can also find more details, for example, in
the book by Carter, Flath and Saito \cite{CFS95}. The hope is that
a similar situation reproduces in dimension four. So, instead of
having a Lie algebra or, equivalently, its universal envelopping
algebra, which has a natural structure of a (trivially braided)
Hopf algebra, we should now have a Hopf category, and instead of
having the monoidal category of representations of the Lie algebra,
we should have the monoidal 2-category of representations of the
Hopf category. Via the reconstruction theorems of the Tannaka-Krein
type, the deformations of the universal envelopping algebra of the
Lie algebra correspond to deformations of its category of
representations. Similarly, deformations of the Hopf category
should correspond to deformations of its 2-category of
representations. Therefore, we are indeed led in this way to
consider the theory of deformations for monoidal 2-categories. In
the above mentioned paper \cite{CF94}, Crane and Frenkel already
outline a method for constructing interesting Hopf categories out
of the quantum groups and their canonical bases. A difficult point
is to find the analog of the quantum groups in this new framework,
which could be called {\it 2-quantum groups} and which would
correspond to non trivial deformations of these Hopf categories.

Let's say a few words about what it means to give a cohomological
description, in the sense of Gerstenhaber, of the theory of
deformations of a semigroupal 2-category. Although we will be thinking of this case, the situation is similar
in all of the above mentioned settings. Given an arbitrary 2-category \Cg, it will be possible to
define more than one semigroupal structure on it. So, we can think
of a space X(\Cg) whose points are in 1-1 correspondence with all
such possible semigroupal structures on
\Cg, up to a suitable notion of equivalence.
The ultimate goal should be to have a description of such a {\it
moduli space} X(\Cg) in terms, for example, of a suitable
parametrization of its points. However, this is difficult. The idea
is then to focus the attention on one particular point in that
space and to study the corresponding ``tangent space''. That's why
we speak of {\it infinitesimal} deformations of the (reference)
semigroupal 2-category. Clearly, the first point is how to
formalize that idea of a tangent space, because a priori we have no
differentiable manifold structure on X(\Cg). In the sequel, we will
see how to do that. We will need to assume some $K$-linear
structure on the 2-category, for some commutative unitary ring $K$,
and to have some local $K$-algebra extending $K$, and over which
the deformations will take place. In the classical algebra setting,
this is accomplished by considering, instead of the original
$K$-algebra $A$, its $K[[h]]$-linear extension $A[[h]]$ (see
\cite{mG63}). According to Gerstenhaber's foundational work, to
give a cohomological description of such infinitesimal deformations
amounts then to find a suitable cohomology $H^{\bullet}(\text\Cg)$
such that the so-called first order deformations (with respect to a
formal deformation parameter) are classified, up to equivalence, by
the elements of one of the cohomology groups $H^n(\text\Cg)$, for
some $n$. But this is only the first point. According to
Gerstenhaber, a nice cohomological description is required to
further satisfy the property that the obstructions to extending
such a first order deformation to higher order deformations or even
to formal series deformations also live in some of the groups
$H^m(\text\Cg)$. In the 0-dimensional setting of algebras, it turns
out that the corresponding obstructions are described by a graded
Lie algebra structure on the cochain complex governing the
deformations \cite{mG64}, and, after Gerstenhaber, this should be a
basic principle of any obstruction theory. As mentioned before,
however, in this paper we will not consider the question of higher order
obstructions, whose treatment is defered to a future paper.
Therefore, the goal of the present work is to just develop the
first of the above points, i.e., to identify the first order
deformations of a semigroupal 2-category with the cocycles of a
suitable cohomology theory.

An important point is how the infinitesimal deformations of a
semigroupal 2-category are defined. In the classical algebra
setting \cite{mG64},\cite{GS88}, recall that the deformation
consists of taking a new (deformed) product $\mu_h$ of the form
$$
\mu_h(a,a')=\mu(a,a')+\mu_1(a,a')h+\mu_2(a,a')h^2+\cdots
$$
where $\mu:A\times A\To A$ denotes the original (undeformed)
product and the $\mu_i:A\times A\To A$, $i\geq 1$, are suitable
$K$-bilinear maps such that $\mu_h$ is indeed associative and with
unit. In the category setting, this should correspond to
considering a new (deformed) tensor product $\otimes_h$ between
morphisms of the form
$$
f\otimes_h g=f\otimes g+(f\otimes_1 g)h+(f\otimes_2 g)h^2+\cdots
$$
where
$\otimes=\otimes_{(X,Y),(X',Y')}:\C\times\C((X,Y),(X',Y'))\To\C(X\otimes
Y,X',\otimes Y')$ corresponds to the original tensor product and
the
$\otimes_i=(\otimes_i)_{(X,Y),(X',Y')}:\C\times\C((X,Y),(X',Y'))\To\C(X\otimes
Y,X',\otimes Y')$, $i\geq 1$, are suitable $K$-bilinear functors.
In the category setting, however, we
should further consider possible deformations of the structural
isomorphisms taking account of the associativity and unit character
of the deformed tensor product, i.e., we should consider, for
example, a new (deformed) associator $a_h$ of the form
$$
(a_h)_{X,Y,Z}=a_{X,Y,Z}+a^{(1)}_{X,Y,Z}h+a^{(2)}_{X,Y,Z}h^2+\cdots
$$
for suitable morphisms $a^{(i)}_{X,Y,Z}:X\otimes(Y\otimes Z)\To
(X\otimes Y)\otimes Z$. Now, in the definition of an infinitesimal
deformation of a monoidal category $\C$, as given by Crane and
Yetter \cite{CY981}, the only structure susceptible of being
deformed is that defined by these structural isomorphisms
\footnote{ As shown by Yetter \cite{dY01}, it is enough to only
consider deformations of the associator, since they already induce
a deformation of the unital structure.}. In other words, the tensor
product $\otimes:\C\times\C\To\C$ is assumed to remain the same
(except for a trivial linear extension). Apart from the fact that
this clearly simplifies the theory, there is another reason that
may induce to adopt this point of view. Indeed, in his milestone
paper \cite{vD91}, Drinfeld proved that the category of
representations of the quantum group $U_h(\text\gg)$ associated to
a simple Lie algebra \gg, which corresponds to a certain
deformation in the above generic sense of the category of
representations of $U(\text\gg)$, is in fact equivalent to the
category of representations of $U(\text\gg)[[h]]$ but with a non
trivially deformed associator (see also \cite{cK95}). Hence, at least in this case, it is enough to consider those
deformations where only the isomorphisms included in the monoidal
structure are deformed, keeping the tensor product undeformed.

In defining the infinitesimal deformations of a semigroupal
2-category we will adopt the same point of view as Crane and
Yetter. So, an infinitesimal deformation of a semigroupal 2-category will be defined in such a way that the only things
susceptible of deformation are the {\it 2-isomorphisms}
defining the semigroupal structure on the
2-category. Contrary to the case of monoidal categories, however,
this involves many things. So, among the
2-isomorphisms susceptible of deformation, we can distinguish three
groups: (1) the 2-isomorphisms included in the tensor product,
coming from the weakening of the definition of the tensor product
as a bifunctor, (2) the 2-isomorphisms included in the associator,
and coming from the weakening of the naturality of the maps
$a_{X,Y,Z}$, and (3) the 2-isomorphisms included in the so called
pentagonator, coming from the weakening of the pentagon axiom on
the associator. In a generic infinitesimal deformation, all of them
will be deformed.

The outline of the paper is as follows. In Section 2 we recall the
basic definitions from bicategory theory, together with the
corresponding strictification theorem (MacLane-Pare's theorem). In
Section 3 we generalize to (unitary) pseudofunctors Epstein's
coherence theorem for semigroupal functors \cite{dE66} and
introduce the analog of Crane-Yetter's ``padding'' composition
operators \cite{CY981} in this setting. They are essential in the
development of the theory. Section 4 is devoted to reviewing in
detail the definition of a semigroupal 2-category, giving a
formulation adapted to our purposes, and we also give an explicit
definition of the corresponding notion of morphism, deduced from
the notion of morphisms between tricategories as it appears in the
paper by Gordon, Power and Street \cite{GPS95}. In Section 5, we
give the precise definition of deformation of a semigroupal
2-category we will work with, together with the notion of
equivalence of deformations. For later use, we also define in this
section the notion of purely pseudofunctorial infinitesimal
deformation of a pseudofunctor. In Section 6, and using results
from Section 3, we develop a cohomology theory for the purely
pseudofunctorial infinitesimal deformations of a pseudofunctor,
which partially generalizes Yetter's theory for monoidal functors
\cite{dY98}. The cohomological description of the deformations of a
semigroupal 2-category is then initiated in Section 7, where we
consider the particular case of the {\it
pentagonator-deformations}, i.e., those deformations where only the
pentagonator is deformed, all the other structural 2-isomorphisms
remaining undeformed. The next section is devoted to determine a
cohomological description of the infinitesimal deformations
involving both the tensor product and the associator. We do that in
the special case where the deformations are {\it unitary}, i.e.,
such that the structural 2-isomorphisms $\otimes_0(X,Y)$ remain
undeformed. They will be called {\it unitary
(tensorator,associator)-deformations}. We also identify
cohomologies describing the deformations separately of both
structures. For the sake of simplicity, in this section we restrict
ourselves to the case of a Gray semigroup. Finally, in Section 9 we
show how the cohomologies in Sections 7 and 8 fit together to give
a cohomology which describes the generic (unitary) deformations.

\section{Basic concepts from bicategory theory}

\subsection{}
Recall that a bicategory, also called lax or weak 2-category, and
first defined by B\'enabou \cite{jB67}, can be obtained from a
category after doing the following two steps: (1) enrich the sets
of morphisms with the category of small categories in the sense of
Kelly \cite{gK82}, and (2) weaken the associativity and unit axioms
on the composition by substituting 2-isomorphisms for the
equations, with the consequent introduction of coherence relations
on these 2-isomorphisms, as explained in the introduction. When we
do that, we obtain the following definition.

\begin{defn} \label{bicategory}
A bicategory \Cg\ consists of:

\begin{itemize}

\item
A class $|\text\Cg|$ of objects or 0-cells.

\item
For any ordered pair of objects $X,Y\in|\text\Cg|$, a small {\it category} \Cg$(X,Y)$.

\noindent{The} objects of \Cg$(X,Y)$, denoted by $f:X\longrightarrow Y$, are called
1-morphisms or 1-cells, and its morphisms, denoted by
$\tau:f\Longrightarrow f'$, are called 2-morphisms or 2-cells.
Remark that, included in this data \Cg$(X,Y)$, there is a
distinguished identity 2-morphism $1_f:f\Longrightarrow f$ for any
1-morphism $f:X\longrightarrow Y$, and an {\it associative}
composition between 2-morphisms, called vertical composition. Given
1-morphisms $f,f',f'':X\longrightarrow Y$, the vertical composite
of two 2-morphisms $\tau:f\Longrightarrow f'$ and
$\tau':f'\Longrightarrow f''$ will be denoted by $\tau'\cdot\tau$.

\item
For any ordered triple of objects $X,Y,Z\in|\text\Cg|$, a {\it functor}
$$
c_{X,Y,Z}:\text\Cg(X,Y)\times \text\Cg(Y,Z)\longrightarrow \text\Cg(X,Z)
$$

\noindent{These} functors provide us not only the composition
$c_{X,Y,Z}(f,g):=g\circ f$ of two 1-morphisms $f:X\longrightarrow
Y$ and  $g:Y\longrightarrow Z$, but also a second composition
between 2-morphisms, called the horizontal composition, which
involves three objects. If $f,f':X\longrightarrow Y$ and
$g,g':Y\longrightarrow Z$, the horizontal composition
$c_{X,Y,Z}(\tau,\eta)$ between the two 2-morphisms
$\tau:f\Longrightarrow f'$ and $\eta:g\Longrightarrow g'$ will be
denoted by $\eta\circ\tau$. In a general bicategory, this
composition may be non associative.

\item
For any object $X\in|\text\Cg|$, a distinguished 1-morphism
$id_X\in|\text\Cg(X,X)|$.

\item
For any objects $X,Y,Z,T\in|\text\Cg|$ and any composable
1-morphisms $f:X\longrightarrow Y$, $g:Y\longrightarrow Z$,
$h:Z\longrightarrow T$, a 2-isomorphism $\alpha_{h,g,f}:h\circ
(g\circ f)\Longrightarrow (h\circ g)\circ f$,  called the
associator or associativity constraint on $f,g,h$.

\item
For any 1-morphism $f:X\longrightarrow Y$, two 2-isomorphisms
$\lambda_f:id_Y\circ f\Longrightarrow f$ and $\rho_f:f\circ
id_X\Longrightarrow f$, called the left and right unit constraints
on $f$, respectively.
\end{itemize}
Moreover, these data must satisfy the following axioms:

\begin{enumerate}
\item
The $\alpha_{h,g,f}$ are natural in $f,g,h$ and the $\lambda_f$ and
$\rho_f$ natural in $f$.

\item
The associator $\alpha=\{\alpha_{h,g,f}\}$ is such that the
following diagram commutes:

$$
\xymatrix{
k\circ(h\circ(g\circ f))\ar[d]_{\alpha_{k,h,g\circ f}}
\ar[rr]^{1_k\circ \alpha_{h,g,f}} & & k\circ((h\circ g)\circ
f)\ar[d]^{\alpha_{k,h\circ g,f}}
\\ (k\circ h)\circ (g\circ f)\ar[dr]_{\alpha_{k\circ h,g,f}} &
& (k\circ(h\circ g))\circ f \ar[dl]^{\alpha_{k,h,g}\circ 1_f}
\\ & ((k\circ h)\circ g)\circ f }
$$

\item
The left and right unit constraints $\lambda=\{\lambda_f\}$ and
$\rho=\{\rho_f\}$ make commutative the following diagram:

$$
\xymatrix{
(g\circ id_Y)\circ f\ar[rr]^{\alpha_{g,id_Y,f}}
\ar[dr]_{\rho_g\circ 1_f} & & g\circ(id_Y\circ f)\ar[dl]^{1_g\circ \lambda_f} \\ & g\circ f }
$$

\end{enumerate}

When all the associators $\alpha_{h,g,f}$ and left and right unit
constraints $\lambda_f$, $\rho_f$ are identities, which in
particular means that the composition of 1-morphisms is strictly
associative and the identity 1-morphisms are strict units, we will
speak of a 2-category.
\end{defn}

The reader should check that a bicategory \Cg\ with only one object
corresponds exactly to the notion of a monoidal category. If $X$ is
the only object of \Cg, the monoidal category is $\text\Cg(X,X)$
with the composition functor as tensor product. This fact will be
used repeatedly in what follows. We also leave to the reader to
check that in a 2-category horizontal composition is strictly
associative and that the identity 2-morphisms of the identity
1-morphisms act as strict units with respect to horizontal
composition. Both facts will also be frequently used.

Given two bicategories \Bg\ and \Cg, their cartesian product is
defined as the bicategory $\text\Bg\times\text\Cg$ such that
\begin{align*} &
|\text\Bg\times\text\Cg|:=|\text\Bg|\times|\text\Cg| \\
&(\text\Bg\times\text\Cg)((X,Y),(X',Y')):=\text\Bg(X,X')\times\text\Cg(Y,Y')
\\ &c_{(X,Y),(X',Y'),X'',Y'')}^{\text\Bg\times\text\Cg}:=(c_{X,X',X''}^{\text\Bg}\times
c_{Y,Y',Y''}^{\text\Cg})\circ P_{23}
\end{align*}
with identity 1-morphisms $id_{(X,Y)}=(id_X,id_Y)$ and whose
structural 2-isomorpisms $\alpha_{(f'',g''),(f',g'),(f,g)}$,
$\lambda_{(f,g)}$ and $\rho_{(f,g)}$ are componentwise given by
those of
\Bg\ and \Cg\ ($P_{23}$ denotes the functor which permutes factors
2 and 3). The same construction obviously extends to a finite
number of bicategories.

We will mainly work with 2-categories. This means no loss of
generality because of the following strictification theorem for
bicategories, due to MacLane and Pare \cite{MP85} (see also
\cite{GPS95}, $\S 1.3$):

\begin{thm} \label{MacLane-Pare}
Any bicategory is biequivalent (in the sense defined below) to a
2-category.
\end{thm}

Diagramatically, a 2-category differs from a category in that it
has vertices (the objects) and edges (the 1-morphisms) but also
faces between pairs of edges. In other words, while a category can
be represented as a 1-dimensional cellular complex, a 2-category is
a 2-dimensional cellular complex. As a consequence, when working
with 2-categories, a generic diagram will be a three dimensional
one, with a new ``pasting'' game where both vertical and horizontal
compositions are combined. We will find some examples in the
sequel. Another significant difference is that in 2-categories (and
in bicategories in general), we have 1-isomorphisms (i.e.,
invertible 1-morphisms), but also {\it equivalences}, i.e.,
1-morphisms which are invertible only up to a 2-isomorphism. This
leads to the notion of equivalent objects in a 2-category, which is
weaker than the notion of isomorphic objects.

\subsection{}
We will also need the corresponding notion of morphism between
bicategories. There are in the literature various versions and
names for this notion.  Following Gray \cite{jG74} I will call them
pseudofunctors, although our definition differs slightly from that
of Gray.

\begin{defn} \label{pseudofunctor}
If \Cg\ and \Dg\ are two bicategories, a {\sl pseudofunctor} from
\Cg\ to
\Dg\ is any quadruple $\F=(|\F|,\F_*,\widehat{\F}_*,\F_0)$, where
\begin{itemize}
\item
$|\F|:|\text\Cg|\longrightarrow|\text\Dg|$ is an object map (the
image $|\F|(X)$ of $X\in|\text\Cg|$ will be denoted by $\F(X)$);

\item
$\F_*=\{\F_{X,Y}:\text\Cg(X,Y)\longrightarrow\text\Dg(\F(X),\F(Y))\}$
is a collection of functors, indexed by ordered pairs of objects
$X,Y\in|\text\Cg|$;

\item
$\widehat{\F}_*=\{\widehat{\F}_{X,Y,Z}:
c^{\text\Dg}_{\F(X),\F(Y),\F(Z)}\circ(\F_{X,Y}\times\F_{Y,Z})
\Longrightarrow\F_{X,Z}\circ c^{\text\Cg}_{X,Y,Z}\}$
is a family of natural isomorphisms,
indexed by triples of objects $X,Y,Z\in|\text\Cg|$. Explicitly,
this corresponds to having, for all composable 1-morphisms
$X\stackrel{f}{\longrightarrow} Y\stackrel{g}{\longrightarrow} Z$,
a 2-isomorphism
$$
\widehat{\F}_{X,Y,Z}(g,f):\F_{Y,Z}(g)\circ\F_{X,Y}(f)\Longrightarrow\F_{X,Z}(g\circ f)
$$
natural in $(f,g)$, and

\item
$\F_0=\{\F_0(X):\F_{X,X}(id_X)\Longrightarrow id_{\F(X)}\}$ is a
collection of 2-isomorphisms, indexed by objects $X\in|\text\Cg|$.
\end{itemize}
\noindent{Moreover}, this data must satisfy the following conditions:

\begin{enumerate}
\item
({\sl hexagonal axiom}) for all composable 1-morphisms
$X\stackrel{f}{\longrightarrow} Y\stackrel{g}{\longrightarrow}
Z\stackrel{h}{\longrightarrow} T$, it commutes
$$
\xymatrix{
\F(h)\circ(\F(g)\circ\F(f))\ar[r]^{1_{\F(h)}\circ\widehat{\F}(g,f)}
\ar[d]_{\alpha_{\F(h),\F(g),\F(f)}} & \F(h)\circ\F(g\circ f)\ar[r]^{\widehat{\F}(h,g\circ f)} &
\F(h\circ (g\circ f))\ar[d]^{\F(\alpha_{h,g,f})}
\\ (\F(h)\circ\F(g))\circ\F(f)\ar[r]_{\widehat{\F}(h,g)\circ 1_{\F(f)}} &
\F(h\circ g)\circ\F(f)\ar[r]_{\widehat{\F}(h\circ g,f)} & \F((h\circ g)\circ f) }
$$

\item
({\sl triangular axioms}) for any 1-morphism $f:X\longrightarrow
Y$, the following diagrams commute:
$$
\xymatrix{
\F(f)\circ id_{\F(X)}\ar[dr]_{\rho_{\F(f)}} & \F(f)\circ\F(id_X)
\ar[l]_{1_{\F(f)}\circ \F_0(X)} \ar[r]^{\widehat{\F}(f,id_X)} &
\F(f\circ id_X)\ar[ld]^{\F(\rho_f)} \\ & \F(f) & \\
id_{\F(Y)}\circ\F(f)\ar[dr]_{\lambda_{\F(f)}} &
\F(id_Y)\circ\F(f)\ar[l]_{\F_0(Y)\circ 1_{\F(f)}}
\ar[r]^{\widehat{\F}(id_Y,f)} & \F(id_Y\circ f)\ar[ld]^{\F(\lambda_f)} \\ & \F(f) &
}
$$

\end{enumerate}
(here, and from now on, we just write $\widehat{\F}(g,f)$ and
$\F(f)$, the indexing objects being omitted for short).

The $\widehat{\F}(g,f)$ and $\F_0(X)$, for all objects $X$ and
composable 1-morphisms $f,g$, will be called the {\sl structural
2-isomorphisms of} $\F$, and the whole set will be called the {\sl
pseudofunctorial structure on} $\F$. When they are all identities,
which in particular means that the functors $\F_{X,Y}$ preserve the
composition of 1-morphisms and the identity 1-morphisms, the
pseudofunctor will be called a {\sl 2-functor}. When only the
$\F_0(X)$ are identities, we will call it a {\sl unitary
pseudofunctor}.

If no confusion arises, a pseudofunctor
$\F=(|\F|,\F_*,\widehat{\F}_*,\F_0)$ from \Cg\ to \Dg\ will be
denoted by $\F:\text\Cg\longrightarrow\text\Dg$ or simply by $\F$.
\end{defn}

\begin{rem}
The only difference between this definition and the one by Gray in
\cite{jG74} is our assumption that all structural 2-morphisms
$\widehat{\F}(g,f)$ and $\F_0(X)$ are actually 2-isomorphisms.
\end{rem}

From the above definition, it follows immediately that a
pseudofunctor between one object bicategories amounts to a monoidal
functor between the corresponding monoidal categories, as the
reader should check.

Given two pseudofunctors $\F:\text\Bg\To\text\Cg$ and
$\G:\text\Cg\To\text\Dg$, the composite pseudofunctor $\G\circ\F$
is defined by
\begin{eqnarray*}
&|\G\circ\F|=|\G|\circ|\F| \\
&(\G\circ\F)_{X,Y}=\G_{\F(X),\F(Y)}\circ\F_{X,Y} \\
&(\widehat{\G\circ\F})(g,f)=\G(\widehat{\F}(g,f))\cdot\widehat{\G}(\F(g),\F(f))
\\ &(\G\circ\F)_0(X)=\G_0(\F(X))\cdot\G(\F_0(X))
\end{eqnarray*}
The direct product of pseudofunctors can also be defined, whose
source and target bicategories are the corresponding product
bicategories. We leave to the reader to write out the explicit
definition. Finally, let's recall that a pseudofunctor
$\F:\text\Cg\To\text\Dg$ is called a {\sl biequivalence} if for any
object $Y\in|\text\Dg|$, there exists an object $X\in|\text\Cg|$
whose image $\F(X)$ is equivalent to $Y$, and for any pair of
objects $X,X'\in|\text\Cg|$ the functor $\F_{X,X'}$ is an
equivalence of categories.

\subsection{}
As in the case of categories, there is a notion of morphism between
pseudofunctors, which I will call pseudonatural transformations.

\begin{defn} \label{trans_quasi}
Let \Cg\ and \Dg\ be two bicategories, and
$\F,\G:\text\Cg\longrightarrow\text\Dg$ two pseudofunctors. Then a
{\sl pseudonatural transformation from} $\F$ {\sl to} $\G$ is any
pair $\xi=(\xi_*,\widehat{\xi}_*)$, where

\begin{itemize}
\item
$\xi_*=\{\xi_X:\F(X)\longrightarrow\G(X)\}$ is a collection of
1-morphisms, indexed by objects $X\in|\text\Cg|$;

\item
$\widehat{\xi}_*=\{\widehat{\xi}_{X,Y}:c^{\text\Dg}_{\F(X),\G(X),\G(Y)}(\xi_X,-)\circ
\G_{X,Y}\Longrightarrow
c^{\text\Dg}_{\F(X),\F(Y),\G(Y)}(-,\xi_Y)\circ\F_{X,Y}\}$ is a
family of natural isomorphisms, indexed by pairs of objects
$X,Y\in|\text\Cg|$. Explicitly, this means having for any
1-morphism $f:X\longrightarrow Y$ a 2-isomorphism
$\widehat{\xi}_{X,Y}(f):\G_{X,Y}(f)\circ\xi_X\Longrightarrow\xi_Y\circ\F_{X,Y}(f)$,
natural in $f$.
\end{itemize}
\noindent{Moreover}, this data must satisfy the conditions

\begin{enumerate}
\item
for all composable 1-morphisms $X\stackrel{f}{\longrightarrow}
Y\stackrel{g}{\longrightarrow} Z$, the following diagram commutes
$$
\xymatrix{
\G(g)\circ(\G(f)\circ\xi_X)\ar[d]_{\alpha_{\G(g),\G(f),\xi_X}}
\ar[r]^{1_{\G(g)}\circ\widehat{\xi}(f)} & \G(g)\circ(\xi_Y\circ\F(f))
\ar[r]^{\alpha_{\G(g),\xi_Y,\F(f)}} & (\G(g)\circ\xi_Y)\circ\F(f)
\ar[r]^{\widehat{\xi}(g)\circ 1_{\F(f)}}  & (\xi_Z\circ\F(g))\circ\F(f)
\ar[d]^{\alpha^{-1}_{\xi_Z,\F(g),\F(f)}}
\\ (\G(g)\circ\G(f))\circ\xi_X\ar[d]_{\widehat{\G}(g,f)\circ 1_{\xi_X}} &  & &
\xi_Z\circ(\F(g)\circ\F(f))\ar[d]^{1_{\xi_Z}\circ\widehat{\F}(g,f)} \\
\G(g\circ f)\circ\xi_X\ar[rrr]^{\widehat{\xi}(g\circ f)} & & & \xi_Z\circ\F(g\circ f) }
$$

\item
for all objects $X$, the following diagram commutes
$$
\xymatrix{
\xi_X\circ id_{\F(X)}\ar[r]^{\rho_{\xi_X}} & \xi_X
\ar[r]^{\lambda^{-1}_{\xi_X}} & id_{\G(X)}\circ\xi_X
\\ \xi_X\circ\F(id_X)\ar[u]^{1_{\xi_X}\circ\F_0(X)} & &
\G(id_X)\circ\xi_X\ar[u]_{\G_0(X)\circ 1_{\xi_X}}\ar[ll]^{\widehat{\xi}(id_X)} }
$$
\end{enumerate}
(for short again, here and from now on, we will omit the indexing
objects in $\widehat{\xi}_{X,Y}(f)$).

When all the 2-isomorphisms $\widehat{\xi}(f)$ are identities,
$\xi$ will be called a {\sl 2-natural transformation}. On the other
hand, if all the $\xi_X$ are 1-isomorphisms (resp. equivalences),
we will speak of a {\sl pseudonatural isomorphism} (resp. {\sl
pseudonatural equivalence}).

If no confusion arises, a pseudonatural transformation
$\xi=(\xi_*,\widehat{\xi}_*)$ between two pseudofunctors $\F$ and
$\G$ will be denoted by $\xi:\F\Longrightarrow\G$ or simply by
$\xi$.
\end{defn}

\begin{rem}
As in the case of pseudofunctors, there is also a more general
notion of pseudonatural transformation, where the 2-cells
$\widehat{\xi}(f)$ are not assumed to be invertible.
\end{rem}

The reader should be familiar with the formulas for the vertical
and horizontal compositions of pseudonatural transformations.
Explicitly, given pseudofunctors
$\F,\G,\Hl:\text\Cg\longrightarrow\text\Dg$, recall that the
vertical composite of $\xi:\F\Longrightarrow\G$ and
$\zeta:\G\Longrightarrow\Hl$, denoted $\zeta\cdot\xi$, is defined
by
\begin{eqnarray*}
&(\zeta\cdot\xi)_X=\zeta_X\circ\xi_X \\
&\widehat{\zeta\cdot\xi}(f)=\alpha^{\text\Dg}_{\zeta_Y,\xi_Y,\F(f)}
\cdot(1_{\zeta_Y}\circ\widehat{\xi}(f))\cdot
(\alpha^{\text\Dg})^{-1}_{\zeta_Y,\G(f),\xi_X}\cdot
(\widehat{\zeta}(f)\circ
1_{\xi_X})\cdot\alpha^{\text\Dg}_{\Hl(f),\zeta_X,\xi_X}
\end{eqnarray*}
On the other hand, given pseudofunctors
$\F,\F':\text\Bg\To\text\Cg$ and $\G:\text\Cg\To\text\Dg$, and a
pseudonatural transformation $\xi:\F\Longrightarrow\F'$, the
horizontal composition $1_{id_{\G}}\circ\xi$ is defined by
\begin{eqnarray*}
&(1_{id_{\G}}\circ\xi)_X=\G_{\F(X),\F'(X)}(\xi_X) \\
&\widehat{1_{id_{\G}}\circ\xi}(f)=\widehat{\G}^{-1}(\xi_Y,\F(f))\cdot
\G(\widehat{\xi}(f))\cdot\widehat{\G}(\F'(f),\xi_X)
\end{eqnarray*}
Finally, the horizontal composition $\zeta\circ 1_{id_{\F}}$, where
$\zeta:\G\Longrightarrow\G':\text\Cg\To\text\Dg$ is any
pseudonatural transformation and $\F:\text\Bg\To\text\Cg$ any
pseudofunctor, is given by
\begin{eqnarray*}
&(\zeta\circ 1_{id_{\F}})_X=\zeta_{\F(X)} \\ &\widehat{\zeta\circ
1_{id_{\F}}}(f)=\widehat{\zeta}(\F(f))
\end{eqnarray*}

\subsection{}
Let's finish this section by recalling that, in the context of
bicategories, there is still a notion of morphism between two
pseudonatural transformations, usually called a modification, and
which has no analog in the category setting.

\begin{defn} \label{modification}
Let \Cg\ and \Dg\ be two bicategories,
$\F,\G:\text\Cg\longrightarrow\text\Dg$ two pseudofunctors and
$\xi,\zeta:\F\Longrightarrow\G$ two pseudonatural transformations.
Then, a {\sl modification from} $\xi$ {\sl to} $\zeta$ is any
family of 2-morphisms
$\text\nn=\{\text\nn_X:\xi_X\Longrightarrow\zeta_X\}$, indexed by
the objects of \Cg, such that for any 1-morphism
$f:X\longrightarrow Y$ in \Cg, it holds
$$
\widehat{\zeta}(f)\cdot(1_{\G(f)}\circ\text\nn_X)=(\text\nn_Y\circ 1_{\F(f)})
\cdot\widehat{\xi}(f).
$$
This condition expresses the fact that the 2-morphisms $\text\nn_X$
are natural in $X$. A modification from $\xi$ to $\zeta$ will be
denoted by $\text\nn:\xi\Longrightarrow\zeta$ or simply by \nn\ if
no confusion arises.

A family of 2-morphisms as above which not necessarily satisfy the
previous naturality condition will be called a {\sl
pseudomodification from} $\xi$ {\sl to} $\zeta$. This more general
notion will be needed later.
\end{defn}

\section{Coherence and padding operators for unitary pseudofunctors}

\subsection{}
Before giving the definition of a semigroupal 2-category and the
corresponding notion of morphism, we consider in this section a
coherence theorem for unitary pseudofunctors which generalizes to
the many objects setting Epstein's coherence theorem for
semigroupal functors \cite{dE66}. Such a coherence result allows us
to introduce the analog of Crane-Yetter's ``padding'' composition
operators \cite{CY981} in this setting. As the reader will realize
later, these results are essential in what follows. So, they are
first used in Section 6 to associate a cochain complex to a unitary
pseudofunctor describing its purely pseudofunctorial infinitesimal
deformations, and which is a key ingredient in the definition of the
double complex of a $K$-linear Gray semigroup introduced in Section
8. The coherence result is also needed to prove that this double
complex is indeed a double complex.

\subsection{}
Recall that in \cite{dE66}, for any pair of
semigroupal\footnote{The term semigroupal applied to categories
means a category with a tensor product and a coherent associator,
but without unit constraints; when applied to functors, it means a
functor $F$ between semigroupal (or monoidal) categories together
with a coherent natural isomorphism
$\widehat{F}:\otimes\circ(F\times F)\Longrightarrow F\circ\otimes$,
but without the isomorphism $F_0:F(I)\longrightarrow I$. Note also
that Epstein actually considers semigroupal categories equipped
with a symmetry.} categories $(\B,\otimes,a)$ and
$(\hat{\B},\hat{\otimes},\hat{a})$ and a semigroupal functor
$(G,\widehat{G})$, the author defines the $G$-{\it iterates} (of
multiplicity $n$, $n\geq 1$) as the set of all functors
$\B^n\longrightarrow\hat{\B}$ that can be obtained as compositions
of product functors $G^i:\B^i\longrightarrow\hat{\B}^i$ ($i\leq n$)
and $j$-iterates ($j\leq n$) of the tensor products $\otimes$ and
$\hat{\otimes}$ in $\B$ and $\hat{\B}$, respectively, which are
functors $\B^j\longrightarrow\B$ or
$\hat{\B}^j\longrightarrow\hat{\B}$. So, a generic $G$-iterate (of
multiplicity $n$) will apply the object $(A_1,\ldots,A_n)$ of
$\B^n$ to an object of $\hat{\B}$ of the form
$$
G(A_1\otimes\cdots\otimes A_{i_1})\hat{\otimes}
G(A_{i_1+1}\otimes\cdots\otimes A_{i_2})\hat{\otimes}\cdots\hat{\otimes}
G(A_{i_r}\otimes\cdots\otimes A_n),
$$
with a suitable parenthesization of the $A_i$'s inside each group
and of the $\hat{\otimes}$-factors that we omit because it will
depend on the $\otimes$- and $\hat{\otimes}$-iterates used. Let $\C
at(G,\otimes,\hat{\otimes})$ be the category whose objects are all
these $G$-iterates and whose morphisms are all the natural
transformations between them. Then, the structural natural
isomorphisms $\widehat{G},a,\hat{a}$ define a subcategory $\C
at(G,\otimes,\hat{\otimes},\widehat{G},a,\hat{a})$ with the same
objects as $\C at(G,\otimes,\hat{\otimes})$ but whose morphisms are
only those induced by these natural isomorphisms
$\widehat{G},a,\hat{a}$, which are called the canonical ones. More
precisely, a canonical morphism is any morphism obtained as a
compositions of expansions of instances of $\widehat{G},a,\hat{a}$
or its inverses, where by an expansion of a morphism $f$ one means
any morphism obtained from $f$ by tensorially multiplying it by
identity morphisms. For example, a canonical morphism from
$G((A\otimes B)\otimes(C\otimes D))$ to $G((A\otimes B)\otimes
C)\hat{\otimes}G(D)$ is
$$
(\widehat{G}(A\otimes B,C)\hat{\otimes} id_{G(D)})\circ
\hat{a}_{G(A\otimes B),G(C),G(D)}\circ (id_{G(A\otimes B}
\hat{\otimes}\widehat{G}(C,D)^{-1})\circ\widehat{G}(A\otimes B,C\otimes D)^{-1}.
$$
A priori, there are other canonical isomorphisms between the same
two objects, as the reader should check. Epstein's coherence
theorem states that, when $\widehat{G},a,\hat{a}$ satisfy the
appropiate coherence relations, for any two objects of $\C
at(G,\otimes,\hat{\otimes},\widehat{G},a,\hat{a})$ there is at most
one morphism (actually an isomorphism).

Let's consider now the case of a unitary pseudofunctor
$(|\G|,\G_*,\widehat{\G}_*)$ between two bicategories \Bg,
$\hat{\text\Bg}$. To emphasize the similarity between both
situations, we present here the ``conversion table'':
\begin{eqnarray*}
- & \longleftrightarrow & |\text\Bg|=\{X,Y,Z,\ldots\} \\
- & \longleftrightarrow & |\hat{\text\Bg}|=\{U,V,W,\ldots\} \\
\B & \longleftrightarrow & \text\Bg_*=\{\text\Bg(X,Y)\}_{X,Y} \\
\hat{\B} & \longleftrightarrow & \hat{\text\Bg}_*=\{\hat{\text\Bg}(U,V)\}_{U,V} \\
\otimes & \longleftrightarrow & c_*=\{c_{X,Y,Z}\}_{X,Y,Z} \\
\hat{\otimes} & \longleftrightarrow & \hat{c}_*=\{\hat{c}_{U,V,W}\}_{U,V,W} \\
a=\{a_{A,B,C}\}_{A,B,C} & \longleftrightarrow &
\alpha_*=\{\alpha(X,Y,Z,T)=\{\alpha_{h,g,f}\}_{h,g,f}\}_{X,Y,Z,T}
\\
\hat{a}=\{\hat{a}_{D,E,F}\}_{D,E,F} & \longleftrightarrow & \hat{\alpha}_*=
\{\hat{\alpha}(U,V,W,S)=\{\hat{\alpha}_{t,s,r}\}_{t,s,r}\}_{U,V,W,S} \\
- & \longleftrightarrow & |\G|:|\text\Bg|\longrightarrow|\hat{\text\Bg}| \\
G:\B\longrightarrow\hat{\B} & \longleftrightarrow &
\G_*=\{\G_{X,Y}:\text\Bg(X,Y)\longrightarrow\hat{\text\Bg}(\G(X),\G(Y))\}_{X,Y}
\\
\widehat{G}=\{\widehat{G}(A,B)\}_{A,B} & \longleftrightarrow & \widehat{\G}_*=
\{\widehat{\G}_{X,Y,Z}=\{\widehat{\G}(g,f)\}_{g,f}\}_{X,Y,Z}
\end{eqnarray*}
Looking at this table, we see that to go to the
bicategory-pseudofunctor setting simply means substituting any
thing in the left-hand side by a family of things of exactly the
same type and indexed by objects of the appropriate bicategory. One
can now proceed in the same way as Epstein does. The $\G$-iterates
of multiplicity $n$, $n\geq 2$, will now be functors
$\text\Bg(X_1,X_2)\times\text\Bg(X_2,X_3)\times\cdots\times
\text\Bg(X_n,X_{n+1})\longrightarrow\hat{\text\Bg}(\G(X_1),\G(X_{n+1}))$,
indexed by a collection $X_1,\ldots,X_{n+1}$ of $n+1$ objects of
\Bg, and obtained as compositions of product functors
$\G_{X_i,X_{i+1}}\times\G_{X_{i+1},X_{i+2}}\times\cdots\times\G_{X_j,X_{j+1}}$
with suitable iterates of the composition functors $c_{X,Y,Z}$ and
$\hat{c}_{U,V,W}$. So, a generic $\G$-iterate will apply the
1-morphisms
$(f_n,\ldots,f_1)\in\text\Bg(X_1,X_2)\times\text\Bg(X_2,X_3)
\times\cdots\times\text\Bg(X_n,X_{n+1})$
to a 1-morphism in $\hat{\text\Bg}(\G(X_1),\G(X_{n+1}))$ of the
form
$$
\G(f_1\circ\cdots\circ f_{i_1})\circ\G(f_{i_1+1}\circ\cdots
\circ f_{i_2})\circ\cdots\circ\G(f_{i_{r+1}}\circ\cdots\circ f_{n-1})\ \ \ \ \ \ \ \ (*)
$$
with the appropriate parenthesization according to the used
composition functors. Similarly, a canonical 2-morphism will be any
2-morphism (actually, a 2-isomorphism) obtained as a vertical
composition of expansions of instances of the given
$\widehat{\G}_{X,Y,Z},\alpha(X,Y,Z),\hat{\alpha}(U,V,W)$ or its
inverses, for all $X,Y,Z\in|\text\Bg|$ and all
$U,V,W\in|\hat{\text\Bg}|$, where expansion now means the
horizontal composition with identity 2-morphisms. We can then
consider the analogs of the above categories, namely, $\C
at(|\G|,\G_*,c_*,\hat{c}_*)$ and
$
\C at(|\G|,\G_*,c_*,\hat{c}_*,\widehat{\G}_*,\alpha_*,\hat{\alpha}_*).
$
We have then the following generalization of Epstein's theorem to
unitary pseudofunctors:

\begin{thm} \label{coherence_pseudofunctor}
Let \Bg, $\hat{\text\Bg}$ be two bicategories and let
$(|\G|,\G_*,,\widehat{\G}_*)$ be a unitary pseudofunctor between
them. Then for any pair of objects of $\C
at(|\G|,\G_*,c_*,\hat{c}_*,\widehat{\G}_*,\alpha_*,\hat{\alpha}_*)$
there is at most one morphism.
\end{thm}
\begin{proof}
Formally, the proof is the same as that of Epstein, but ignoring
the permutations which appear in his paper because we do not
consider commutativity constraints. The main difference is that we
work simultaneously with various functors and natural isomorphisms.
\end{proof}

\begin{rem}
This coherence theorem already appears in a different formulation
in \cite{GPS95}, $\S 1.6$.
\end{rem}

\subsection{}
\begin{figure}
\centering
\input{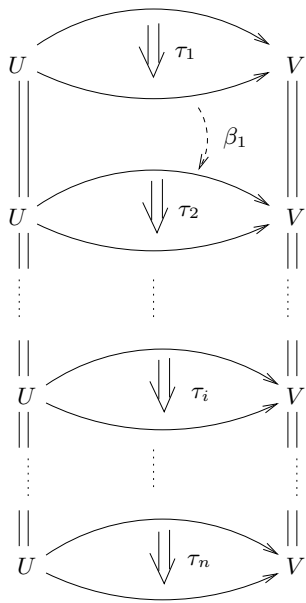}
\caption{Padding for the vertical composition of 2-morphisms.}
\label{figurapadding}
\end{figure}
The previous result allows us to introduce the analog of
Crane-Yetter's ``padding'' composition operators \cite{CY981} in
the context of a unitary pseudofunctor $\G$ between two
bicategories
\Bg\ and $\hat{\text\Bg}$. The main difference is that now we have
a whole collection of such padding operators, indexed by pairs of
objects of the target bicategory $\hat{\text\Bg}$. So, given two
such objects $U,V$, the situation is that depicted in
Fig.~\ref{figurapadding}. We have a sequence $\tau_1,\ldots,\tau_n$
of 2-morphisms in $\hat{\text\Bg}$ such that the source 1-morphism
of $\tau_{i+1}$ is canonically 2-isomorphic to the target
1-morphism of $\tau_i$ (i.e., they are 2-isomorphic through a
composition of expansions of the structural 2-isomorphisms coming
from $\G,\text\Bg,\hat{\text\Bg}$). Then, define
$$
\lceil\tau_n\cdot\tau_{n-1}\cdot\cdots\cdot\tau_1\rceil_{U,V}:=
\beta_n\cdot\tau_n\cdot\beta_{n-1}\cdot\tau_{n-1}\cdot\cdots\cdot
\beta_1\cdot\tau_1\cdot\beta_0,
$$
where the $\beta_i$'s are the canonical 2-isomorphisms between the
target of $\tau_i$ and the source of $\tau_{i+1}$, $\beta_0$ is the
canonical 2-isomorphism whose source 1-morphism has no identity
composition factors and it is completely right-parenthesized and
free from images of composite morphisms under $\G$, and $\beta_n$
is the canonical 2-isomorphism whose target 1-morphism has no
identity composition factors and it is completely
left-parenthesized and free from compositions both of whose factors
are images under $\G$. Note that these are the padding operators
when one chooses as ``references'' the $\G$-iterates
$c^{(n)}\circ\G^{(n)}$ and $\G\circ ^{(n)}c$, where $c^{(n)}$
denotes the appropriate iterate of the composition functors of
$\hat{\text\Bg}$ for the resulting composition to be completely
right-parenthesized, $^{(n)}c$ the same thing but using the
composition functors of \Bg\ and so that the resulting composition
is completely left-parenthesized, and $\G^{(n)}$ denotes the
appropriate $\G$-iterate (probably with some factor equal to an
identity functor). Other choices of references are also possible.
That the above 2-morphism is well defined is a consequence of the
previous coherence theorem.

\begin{ex}
Let $\G=(|\G|,\G_*,\widehat{\G}_*)$ be a unitary pseudofunctor
between two bicategories \Bg\ and $\hat{\text\Bg}$. Let $X,Y,Z,T$
be objects of \Bg\ and let us consider 1-morphisms
$X\stackrel{f}{\longrightarrow}
Y\stackrel{g}{\longrightarrow}Z\stackrel{h}{\longrightarrow} T$.
Taking $U=\G(X)$ and $V=\G(T)$, we have
$
\lceil1_{\G(h)}\circ\widehat{\G}(g,f)\rceil_{\G(X),\G(T)}=
\G(\alpha^{\text\Bg}_{h,g,f})\cdot
\widehat{\G}(h,g\circ f)\cdot(1_{\G(h)}\circ\widehat{\G}(g,f)).
$
\end{ex}

\section{Semigroupal 2-categories and their morphisms}

\subsection{}
From now on, and unless otherwise indicated, all bicategories will
be assumed to be 2-categories. This assumption does not imply loss
of generality due to MacLane-Pare's strictification theorem for
bicategories (see Theorem~\ref{MacLane-Pare}).

\subsection{}
The objects of our interest are the semigroupal 2-categories.
Recall that semigroupal 2-category is a monoidal 2-category without
the unit object for the tensor product, and hence without the
structural 1- and 2-isomorphisms related to the unital structure.

A standard reference on monoidal 2-categories is the paper by
Kapranov-Voevodsky \cite{KV94}. In that paper, however, they give
an unraveled definition which involves many data and an even
greater number of axioms. To make things more intelligible, it is
worth to point out that a semigroupal 2-category is just the
categorification of the definition of a semigroupal category. This
naturally leads to the following definition (except for the $K_5$
coherence condition on the pentagonator).

%\parindent=1 truecm
\begin{defn}
A {\sl semigroupal 2-category} consists of the following data {\bf
SBDi} and axiom {\bf SBA}:
\begin{description}
\item[SBD1]
A 2-category \Cg.

\item[SBD2]
A pseudofunctor $\otimes:\text\Cg\times\text\Cg\To\text\Cg$, called
the {\sl tensor product}.

\item[SBD3]
A pseudonatural isomorphism $a:\otimes^{(3)}\Longrightarrow
^{(3)}\otimes:
\text\Cg\times\text\Cg\times\text\Cg\To\text\Cg$, called the {\sl associator},
where $\otimes^{(3)}$ and $^{(3)}\otimes$ denote, respectively, the
composite pseudofunctors $\otimes\circ(id_{\text\Cg}\times\otimes)$
and $\otimes\circ(\otimes\times id_{\text\Cg})$.

\item[SBD4]
An invertible modification $\pi:a^{(4)}\Rightarrow ^{(4)}a:
\otimes^{(4)}\Longrightarrow ^{(4)}\otimes:
\text\Cg\times\text\Cg\times\text\Cg\times\text\Cg\To\text\Cg$, called the {\sl pentagonator},
where $\otimes^{(4)}$ and $^{(4)}\otimes$ denote, respectively, the
composite pseudofunctors
$\otimes\circ(id_{\text\Cg}\times\otimes)\circ(id_{\text\Cg}\times
id_{\text\Cg}\times\otimes)$ and $\otimes\circ(\otimes\times
id_{\text\Cg})\circ(\otimes\times id_{\text\Cg}\times
id_{\text\Cg})$, and $a^{(4)}$, $^{(4)}a$ are the pseudonatural
isomorphisms
\begin{eqnarray*}
a^{(4)}&=&(1_{\otimes}\circ(a\times 1_{id}))\cdot(a\circ
1_{id\times\otimes\times id})\cdot(1_{\otimes}\circ(1_{id}\times
a)) \\ ^{(4)}a&=&(a\circ 1_{\otimes\times id\times id})\cdot(a\circ
1_{id\times id\times\otimes})
\end{eqnarray*}
(here, $id$ denotes the identity 2-functor of \Cg). See
Fig.~\ref{figura_pentagonador}

\begin{figure}
\centering
\input{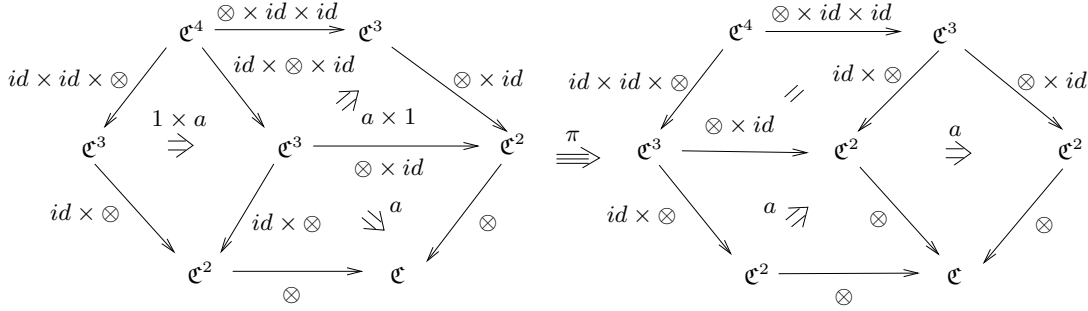}
\caption{Definition of the pentagonator}
\label{figura_pentagonador}
\end{figure}

\item[SBA]
The data $(\otimes,a,\pi)$ is such that the equality in
Fig.~\ref{relacioK5} holds (to simplify notation, the tensor
product of objects or 1-morphisms is denoted by simple
juxtaposition, and the identity 1-morphisms are represented by the
corresponding objects; for more details about the notations in this
Figure, see the next Proposition). This condition will be called
the $K_5$ {\sl coherence relation} (the name comes from the fact
that the two pastings in Fig.~\ref{relacioK5} respresent together a
realization of the $K_5$ Stasheff polytope; see \cite{jS63}).
\end{description}

\begin{figure}
\centering
\input{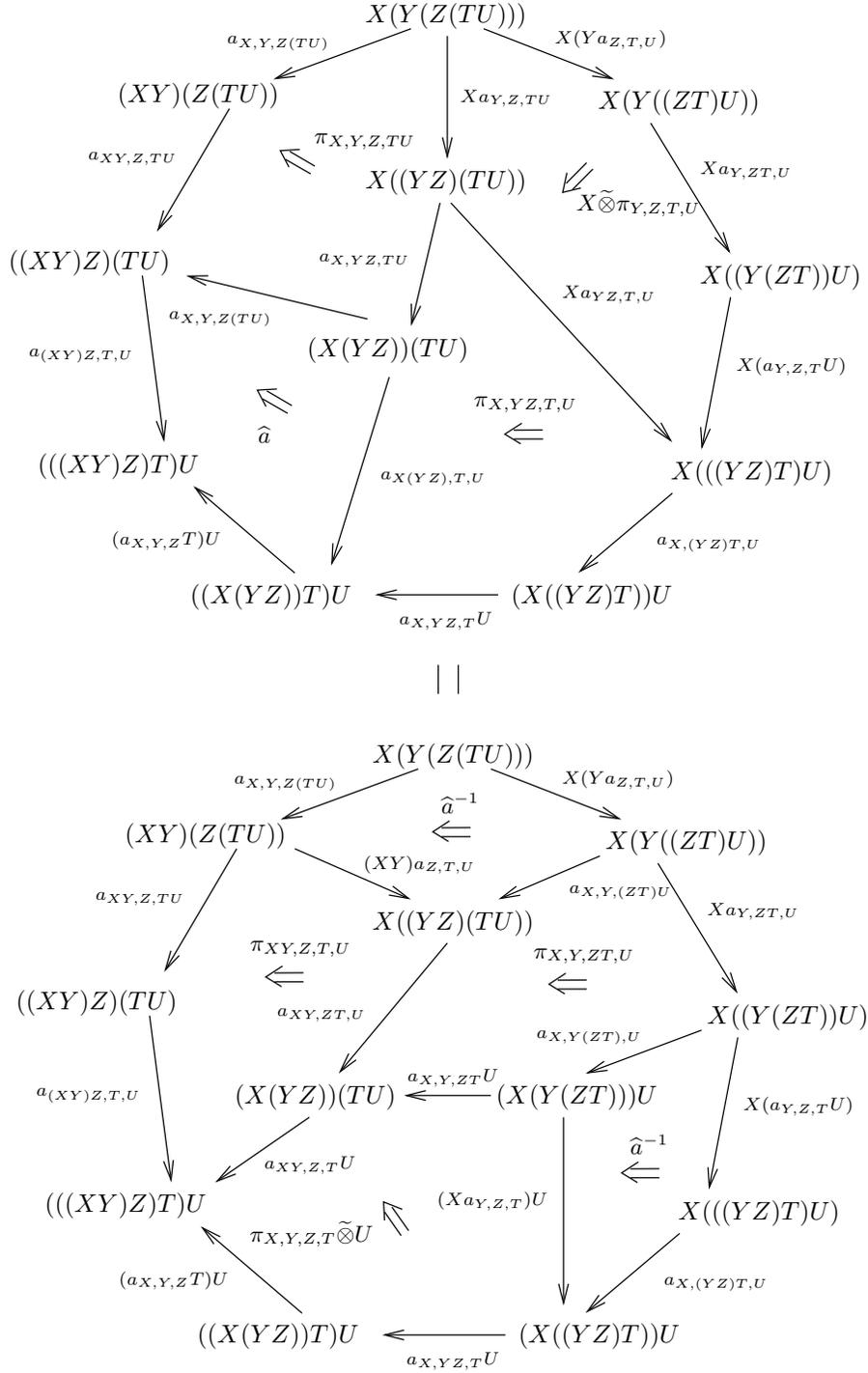}
\caption{$K_5$ coherence relation on the pentagonator}
\label{relacioK5}
\end{figure}

\noindent{A} semigroupal 2-category will be denoted by
$(\text\Cg,\otimes,a,\pi)$ and the triple $(\otimes,a,\pi)$ will be
called a {\sl semigroupal structure} on the 2-category \Cg.
\end{defn}

For convenience in what follows, we give an explicit description of
the structural 1- and 2-isomorphisms involved in the previous
definition, as well as the whole list of equations they must
satisfy.

%\newpage
\begin{prop} \label{semigroupal_2_category}
Let \Cg\ be a 2-category. Then, a semigroupal structure
$(\otimes,a,\pi)$ on \Cg\ consists of the following data:

\begin{description}

\item[$|\otimes|$]
An object $X\otimes Y$ for any object $(X,Y)$ of
$\text\Cg\times\text\Cg$.

\item[$\otimes$]
A collection of functors
$\otimes_{(X,Y),(X',Y')}:\text\Cg(X,X')\times\text\Cg(Y,Y')\To\text\Cg(X\otimes
Y,X'\otimes Y')$ for all $(X,Y),(X',Y')$ objects of
$\text\Cg\times\text\Cg$. As usual, the image of the 1-morphism
$(f,g):(X,Y)\To(X',Y')$ and the 2-morphism
$(\tau,\sigma):(f,g)\Longrightarrow(\tilde{f},\tilde{g}):(X,Y)\To(X',Y')$
by this functor $\otimes_{(X,Y),(X',Y')}$ will be denoted by
$f\otimes g$ and $\tau\otimes\sigma$, respectively.

\item[$\mathbf{\widehat{\otimes}}$]
A collection of 2-isomorphisms
$\widehat{\otimes}((f',g'),(f,g)):(f'\otimes g')\circ(f\otimes
g)\Longrightarrow (f'\circ f)\otimes (g'\circ g)$, for all
composable 1-morphisms $(f,g):(X,Y)\To (X',Y')$ and
$(f',g'):(X',Y')\To(X'',Y'')$ of $\text\Cg\times\text\Cg$.

\item[$\mathbf{\otimes_0}$]
A collection of 2-isomorphisms $\otimes_0(X,Y):id_X\otimes
id_Y\Longrightarrow id_{X\otimes Y}$, for all objects $(X,Y)$ of
$\text\Cg\times\text\Cg$.

\item[a]
A collection of 1-isomorphisms $a_{X,Y,Z}:X\otimes(Y\otimes Z)\To
(X\otimes Y)\otimes Z$, for all objects $(X,Y,Z)$ of
$\text\Cg\times\text\Cg\times\text\Cg$.

\item[$\mathbf{\widehat{a}}$]
A collection of 2-isomorphisms $\widehat{a}(f,g,h):((f\otimes
g)\otimes h)\circ a_{X,Y,Z}\Longrightarrow a_{X',Y',Z'}\circ
(f\otimes(g\otimes h))$, for all 1-morphisms $(f,g,h):(X,Y,Z)\To
(X',Y',Z')$ of $\text\Cg\times\text\Cg\times\text\Cg$.

\item[$\mathbf{\pi}$]
A collection of 2-isomorphisms $\pi_{X,Y,Z,T}:(a_{X,Y,Z}\otimes
id_T)\circ a_{X,Y\otimes Z,T}\circ(id_X\otimes
a_{Y,Z,T})\Longrightarrow a_{X\otimes Y,Z,T}\circ a_{X,Y,Z\otimes
T}$, for all objects $(X,Y,Z,T)$ of
$\text\Cg\times\text\Cg\times\text\Cg\times\text\Cg$.

\end{description}
\noindent{Moreover}, all the above 1- and 2-isomorphisms must satisfy the
following equations:
\begin{description}

\item[A$\mathbf{\widehat{\otimes}}$1]
For all 2-morphisms
$(\tau,\sigma):(f,g)\Longrightarrow(\tilde{f},\tilde{g}):(X,Y)\To(X',Y')$
and
$(\tau',\sigma'):(f',g')\Longrightarrow(\tilde{f}',\tilde{g}'):(X',Y')\To(X'',Y'')$
of $\text\Cg\times\text\Cg$
$$
((\tau'\circ\tau)\otimes(\sigma'\circ\sigma))\cdot
\widehat{\otimes}((f',g'),(f,g))=
\widehat{\otimes}((\tilde{f}',\tilde{g}'),(\tilde{f},\tilde{g}))\cdot
((\tau'\otimes\sigma')\circ(\tau\otimes\sigma))
$$

\item[A$\mathbf{\widehat{\otimes}}$2]
For all composable 1-morphisms
$(X,Y)\stackrel{(f,g)}{\To}(X',Y')\stackrel{(f',g')}{\To}(X'',Y'')
\stackrel{(f'',g'')}{\To}(X''',Y''')$ of $\text\Cg\times\text\Cg$
\begin{eqnarray*}
\lefteqn{\widehat{\otimes}((f'',g''),(f'\circ f,g'\circ g))\cdot
(1_{f''\otimes g''}\circ\widehat{\otimes}((f',g')(f,g)))=} \\ & &
=\widehat{\otimes}((f''\circ f',g''\circ g'),(f,g))\cdot
(\widehat{\otimes}((f'',g''),(f',g'))\circ 1_{f\otimes g})
\end{eqnarray*}

\item[A$\mathbf{\widehat{\otimes}}$3]
For any 1-morphism $(f,g):(X,Y)\To(X',Y')$ of
$\text\Cg\times\text\Cg$
\begin{eqnarray*}
&\widehat{\otimes}((id_{X'},id_{Y'}),(f,g))=\otimes_0(X',Y')\circ
1_{f\otimes g} \\ &\widehat{\otimes}((f,g),(id_X,id_Y))=1_{f\otimes g}\circ\otimes_0(X,Y)
\end{eqnarray*}

\item[A$\mathbf{\widehat{a}}$1]
For all 2-morphisms
$(\tau,\sigma,\eta):(f,g,h)\Longrightarrow(\tilde{f},\tilde{g},\tilde{h}):(X,Y,Z)\To(X',Y',Z')$
of $\text\Cg^3$
$$
(1_{a_{X',Y',Z'}}\circ(\tau\otimes(\sigma\otimes\eta)))\cdot\widehat{a}(f,g,h)=
\widehat{a}(\tilde{f},\tilde{g},\tilde{h})\cdot(((\tau\otimes\sigma)\otimes\eta)
\circ 1_{a_{X,Y,Z}})
$$

\item[A$\mathbf{\widehat{a}}$2]
For all composable 1-morphisms
$(X,Y,Z)\stackrel{(f,g,h)}\To(X',Y',Z')\stackrel{(f',g',h')}{\To}(X'',Y'',Z'')$
of $\text\Cg\times\text\Cg\times\text\Cg$
\begin{align*}
\widehat{a}(f'\circ f,g'\circ g,h'\circ
h)&\cdot((\widehat{\otimes}((f',g'),(f,g))\otimes 1_{h'\circ
h})\circ 1_{a_{X,Y,Z}})
\\ &\cdot(\widehat{\otimes}((f'\otimes
g',h'),(f\otimes g,h))\circ 1_{a_{X,Y,Z}})=
\\ &\ \ \ \ \ \ \ \ \ \ \ \ =(1_{a_{X'',Y'',Z''}}\circ(1_{f'\circ
f}\otimes\widehat{\otimes}((g',h')(g,h))))
\\ &\ \ \ \ \ \ \ \ \ \ \ \ \ \ \ \ \cdot(1_{a_{X',Y',Z'}}\circ
\widehat{\otimes}((f',g'\otimes h'),(f,g\otimes h)))
\\ &\ \ \ \ \ \ \ \ \ \ \ \ \ \ \ \ \cdot(\widehat{a}(f',g',h')\circ 1_{f\otimes(g\otimes
h)})\cdot(1_{(f'\otimes g')\otimes h'}\circ\widehat{a}(f,g,h))
\end{align*}
%
%
%\begin{eqnarray*}
%\lefteqn{\widehat{a}(f'\circ
%f,g'\circ g,h'\circ
%h)\cdot([(\widehat{\otimes}((f',g'),(f,g))\otimes 1_{h'\circ
%h})\cdot\widehat{\otimes}((f'\otimes g',h'),(f\otimes g,h))]\circ
%1_{a_{X,Y,Z}})=} \\ & &=(1_{a_{X'',Y'',Z''}}\circ[(1_{f'\circ
%f}\otimes\widehat{\otimes}((g',h')(g,h)))\cdot
%\widehat{\otimes}((f',g'\otimes h'),(f,g\otimes h))])\cdot \\ & &
%\ \ \ \ \ \cdot(\widehat{a}(f',g',h')\circ 1_{f\otimes(g\otimes
%h)})\cdot(1_{(f'\otimes g')\otimes h'}\circ\widehat{a}(f,g,h))
%\end{eqnarray*}

\item[A$\mathbf{\widehat{a}}$3]
For all objects $(X,Y,Z)$ of $\text\Cg\times\text\Cg\times\text\Cg$
\begin{align*}
\widehat{a}(id_X,id_Y,id_Z)=&(1_{a_{X,Y,Z}}\circ(\otimes_0(X,Y\otimes Z)
\cdot(1_{id_X}\otimes\otimes_0(Y,Z)))^{-1})\cdot \\ &\
\ \ \cdot((\otimes_0(X\otimes Y,Z)\cdot(\otimes_0(X,Y)\otimes 1_{id_Z}))\circ
1_{a_{X,Y,Z}})
\end{align*}
%\newpage

\item[A$\pi$1]
For any 1-morphism $(f,g,h,k):(X,Y,Z,T)\To(X',Y',Z',T')$ of
$\text\Cg^4$
\begin{align*}
(1_{a_{X'\otimes Y',Z',T'}}\circ\widehat{a}&(f,g,h\otimes
k))\cdot(\widehat{a}(f\otimes g,h,k)\circ 1_{a_{X,Y,Z\otimes
T}})\cdot(1_{((f\otimes g)\otimes h)\otimes k}\circ\pi_{X,Y,Z,T})=
\\ &=(\pi_{X',Y',Z',T'}\circ 1_{f\otimes(g\otimes(h\otimes
k))})\cdot
\\ &\ \ \ \cdot(1_{a_{X',Y',Z'}\otimes id_{T'}\circ a_{X,Y\otimes
Z,T}}\circ\widehat{\otimes}((id_{X'},a_{Y',Z',T'}),(f,g\otimes(h\otimes
k)))^{-1})\cdot
\\ &\ \ \ \cdot(1_{(a_{X',Y',Z'}\otimes id_{T'})\circ a_{X,Y\otimes Z,T}}\circ
(1_f\otimes\widehat{a}(g,h,k)))\cdot
\\ &\ \ \ \cdot(1_{(a_{X',Y',Z'}\otimes id_{T'})\circ a_{X,Y\otimes Z,T}}\circ
\widehat{\otimes}((f,(g\otimes h)\otimes k),(id_X,a_{Y,Z,T})))\cdot
\\ &\ \ \ \cdot(1_{a_{X',Y',Z'}\otimes id_{T'}}\circ
\widehat{a}(f,g\otimes h,k)\circ 1_{id_X\otimes a_{Y,Z,T}})\cdot
\\ &\ \ \ \cdot(\widehat{\otimes}((a_{X',Y',Z'},id_{T'}),(f\otimes(g\otimes h),k))^{-1}
\circ 1_{a_{X,Y\otimes Z,T}\circ(id_X\otimes a_{Y,Z,T})})\cdot
\\ &\ \ \ \cdot((\widehat{a}(f,g,h)\otimes 1_k)\circ 1_{a_{X,Y\otimes Z,T}
\circ(id_X\otimes a_{Y,Z,T})})\cdot
\\ &\ \ \ \cdot(\widehat{\otimes}(((f\otimes g)\otimes h,k),(a_{X,Y,Z},id_T))
\circ 1_{a_{X,Y\otimes Z,T}\circ(id_X\otimes a_{Y,Z,T})})
\end{align*}
%
%
%\begin{eqnarray*}
%\lefteqn{(1_{a_{X'\otimes
%Y',Z',T'}}\circ\widehat{a}(f,g,h\otimes
%k))\cdot(\widehat{a}(f\otimes g,h,k)\circ 1_{a_{X,Y,Z\otimes
%T}})\cdot(1_{((f\otimes g)\otimes h)\otimes k}\circ\pi_{X,Y,Z,T})=}
%\\ & &=(\pi_{X',Y',Z',T'}\circ 1_{f\otimes(g\otimes(h\otimes
%k))})\cdot(1_{a_{X',Y',Z'}\otimes id_{T'}}\circ[(1_{a_{X,Y\otimes
%Z,T}}\circ\{\widehat{\otimes}((id_{X'},a_{Y',Z',T'}),(f,g\otimes(h\otimes
%k)))^{-1}\cdot
%\\ & &\ \ \cdot(1_f\otimes\widehat{a}(g,h,k))\cdot
%\widehat{\otimes}((f,(g\otimes h)\otimes k),(id_X,a_{Y,Z,T}))\})\cdot
%(\widehat{a}(f,g\otimes h,k)\circ 1_{id_X\otimes a_{Y,Z,T}})])\cdot
%\\ & &\ \ \cdot([\widehat{\otimes}((a_{X',Y',Z'},id_{T'}),(f\otimes(g\otimes h),k))^{-1}
%\cdot(\widehat{a}(f,g,h)\otimes 1_k)\cdot \\ & &
%\ \ \cdot\widehat{\otimes}(((f\otimes g)\otimes h,k),(a_{X,Y,Z},id_T))]
%\circ 1_{a_{X,Y\otimes Z,T}\circ(id_X\otimes a_{Y,Z,T})})
%\end{eqnarray*}

\item[A$\mathbf{\pi}$2]
For any object $(X,Y,Z,T,U)$ of
$\text\Cg\times\text\Cg\times\text\Cg\times\text\Cg\times\text\Cg$
\begin{align*}
(\pi_{X\otimes Y,Z,T,U}&\circ 1_{a_{X,Y,Z\otimes(T\otimes
Z)}})\cdot
\\ &\cdot(1_{(a_{X\otimes Y,Z,T}\otimes id_U)\circ a_{X\otimes Y,Z\otimes
T,U}}\circ\widehat{a}(id_X,id_Y,a_{Z,T,U})^{-1})\cdot
\\ &\cdot(1_{a_{X\otimes Y,Z,T}\otimes id_U}\circ\pi_{X,Y,Z\otimes
T,U}\circ 1_{id_X\otimes(id_Y\otimes a_{Z,T,U})})\cdot
\\ &\cdot((\pi_{X,Y,Z,T}\tilde{\otimes} 1_{id_U})\circ
1_{a_{X,Y\otimes(Z\otimes T),U}\circ(id_X\otimes a_{Y,Z\otimes
T,U})\circ(id_X\otimes(id_Y\otimes a_{Z,T,U}))})=
\\ &\ \ \ \ \ \ \ =(1_{a_{(X\otimes Y)\otimes Z,T,U}}\circ\pi_{X,Y,Z,T\otimes U})\cdot
\\ &\ \ \ \ \ \ \ \ \ \ \ \cdot
(\widehat{a}(a_{X,Y,Z},id_T,id_U)\circ 1_{a_{X,Y\otimes Z,T\otimes
U}\circ(id_X\otimes a_{Y,Z,T\otimes U})})\cdot
\\ &\ \ \ \ \ \ \ \ \ \ \ \cdot
(1_{(a_{X,Y,Z}\otimes id_T)\otimes id_U}\circ
\pi_{X,Y\otimes Z,T,U}\circ 1_{id_X\otimes a_{Y,Z,T\otimes U}})\cdot
\\ &\ \ \ \ \ \ \ \ \ \ \ \cdot
(1_{((a_{X,Y,Z}\otimes id_T)\otimes id_U)\circ(a_{X,Y\otimes Z,T}\otimes id_U)
\circ a_{X,(Y\otimes Z)\otimes T,U}}\circ(1_{id_X}\tilde{\otimes}\pi_{Y,Z,T,U}))\cdot
\\ &\ \ \ \ \ \ \ \ \ \ \ \cdot
(1_{((a_{X,Y,Z}\otimes id_T)\otimes id_U)\circ
(a_{X,Y\otimes Z,T}\otimes id_U)}\circ\widehat{a}(id_X,a_{Y,Z,T},id_U)\circ
\\ &\ \ \ \ \ \ \ \ \ \ \ \ \ \ \ \ \ \ \ \ \ \ \circ
1_{(id_X\otimes a_{Y,Z\otimes T,U})\circ(id_X\otimes(id_Y\otimes
a_{Z,T,U}))})
\end{align*}
%
%\begin{eqnarray*}
%\lefteqn{(\pi_{X\otimes Y,Z,T,U}\circ 1_{a_{X,Y,Z\otimes(T\otimes Z)}})\cdot
%(1_{(a_{X\otimes Y,Z,T}\otimes id_U)\circ a_{X\otimes Y,Z\otimes
%T,U}}\circ\widehat{a}(id_X,id_Y,a_{Z,T,U})^{-1})\cdot} \\ &
%&\cdot(1_{a_{X\otimes Y,Z,T}\otimes id_U}\circ\pi_{X,Y,Z\otimes
%T,U}\circ 1_{id_X\otimes(id_Y\otimes a_{Z,T,U})})\cdot \\ &
%&\cdot((\pi_{X,Y,Z,T}\tilde{\otimes} 1_{id_U})\circ
%1_{a_{X,Y\otimes(Z\otimes T),U}\circ(id_X\otimes a_{Y,Z\otimes
%T,U})\circ(id_X\otimes(id_Y\otimes a_{Z,T,U}))})=
%\\ & &
%\ \ =(1_{a_{(X\otimes Y)\otimes Z,T,U}}\circ\pi_{X,Y,Z,T\otimes U})\cdot
%(\widehat{a}(a_{X,Y,Z},id_T,id_U)\circ 1_{a_{X,Y\otimes Z,T\otimes
%U}\circ(id_X\otimes a_{Y,Z,T\otimes U})})\cdot \\ & &
%\ \ \ \ \ \cdot(1_{(a_{X,Y,Z}\otimes id_T)\otimes id_U}\circ
%\pi_{X,Y\otimes Z,T,U}\circ 1_{id_X\otimes a_{Y,Z,T\otimes U}})\cdot \\ & &
%\ \ \ \ \ \cdot(1_{((a_{X,Y,Z}\otimes id_T)\otimes id_U)\circ(a_{X,Y\otimes Z,T}\otimes id_U)
%\circ a_{X,(Y\otimes Z)\otimes T,U}}\circ(1_{id_X}\tilde{\otimes}\pi_{Y,Z,T,U}))\cdot \\ & &
%\ \ \ \ \ \cdot(1_{((a_{X,Y,Z}\otimes id_T)\otimes id_U)\circ
%(a_{X,Y\otimes Z,T}\otimes id_U)}\circ\widehat{a}(id_X,a_{Y,Z,T},id_U)\circ
%1_{(id_X\otimes a_{Y,Z\otimes T,U})\circ(id_X\otimes(id_Y\otimes
%a_{Z,T,U}))})
%\end{eqnarray*}
\end{description}
\end{prop}

We will refer to the previous equations as the {\sl structural
equations} of a semigroupal 2-category. Notice that, in Equation
(A$\pi$2), both terms $1_{id_X}\tilde{\otimes}\pi_{Y,Z,T,U}$ and
$\pi_{X,Y,Z,T}\tilde{\otimes} 1_{id_U}$ denote pastings of the
corresponding terms $1_{id_X}\otimes\pi_{Y,Z,T,U}$ and
$\pi_{X,Y,Z,T}\otimes 1_{id_U}$. For example, the reader may check
that the first term is given by
\begin{eqnarray*}
1_{id_X}\tilde{\otimes}\pi_{Y,Z,T,U}&=&\widehat{\otimes}((id_X,a_{Y\otimes
Z,T,U}),(id_X,a_{Y,Z,T\otimes
U}))^{-1}\cdot(1_{id_X}\otimes\pi_{Y,Z,T,U})\cdot \\ &
&\cdot\widehat{\otimes}((id_X,a_{Y,Z,T}\otimes
id_U),(id_X,a_{Y,Z\otimes T,U}\circ(id_Y\otimes a_{Z,T,U})))\cdot
\\ & &\cdot(1_{id_X\otimes(a_{Y,Z,T}\otimes
id_U)}\circ\widehat{\otimes}((id_X,a_{Y,Z\otimes
T,U}),(id_X,id_Y\otimes a_{Z,T,U})))
\end{eqnarray*}
A similar expression gives us the pasting
$\pi_{X,Y,Z,T}\tilde{\otimes} 1_{id_U}$.

\begin{proof}
Equations $(A\widehat{\otimes}i)$ correspond to the naturality of
the 2-isomorphisms $\widehat{\otimes}((f',g'),(f,g))$ and the two
axioms on the pseudofunctorial structure which appear in the
definition of pseudofunctor. Similarly, equations $(A\widehat{a}i)$
correspond to the naturality of the 2-isomorphisms
$\widehat{a}(f,g,h)$ and the axioms appearing in the definition of
pseudonatural transformation. On the other hand, Equation $(A\pi
1)$ corresponds to the naturality condition on the pentagonator
$\pi_{X,Y,Z,T}$ in $(X,Y,Z,T)$, namely,
$$
\widehat{^{(4)}a}(f,g,h,k)\cdot(1_{((f\otimes g)\otimes h)\otimes k}
\circ\pi_{X,Y,Z,T})=(\pi_{X',Y',Z',T'}\circ 1_{f\otimes(g\otimes(h\otimes k))})
\cdot\widehat{a^{(4)}}(f,g,h,k)
$$
after making explicit the 2-isomorphisms
$\widehat{^{(4)}a}(f,g,h,k)$ and $\widehat{a^{(4)}}(f,g,h,k)$ using
the definitions in Section 2. Finally, Equation $(A\pi 2)$ is the
algebraic expression of the $K_5$ coherence relation.
\end{proof}

\begin{defn}
A semigroupal 2-category $(\text\Cg,\otimes,a,\pi)$ is called {\sl
strict} when all the above structural isomorphisms are identities
(notice that this is not possible for an arbitrary tensor product
$\otimes$; for example, it must satisfy that $X\otimes(Y\otimes
Z)=(X\otimes Y)\otimes Z$, etc.).
\end{defn}

\begin{figure}
\centering
\input{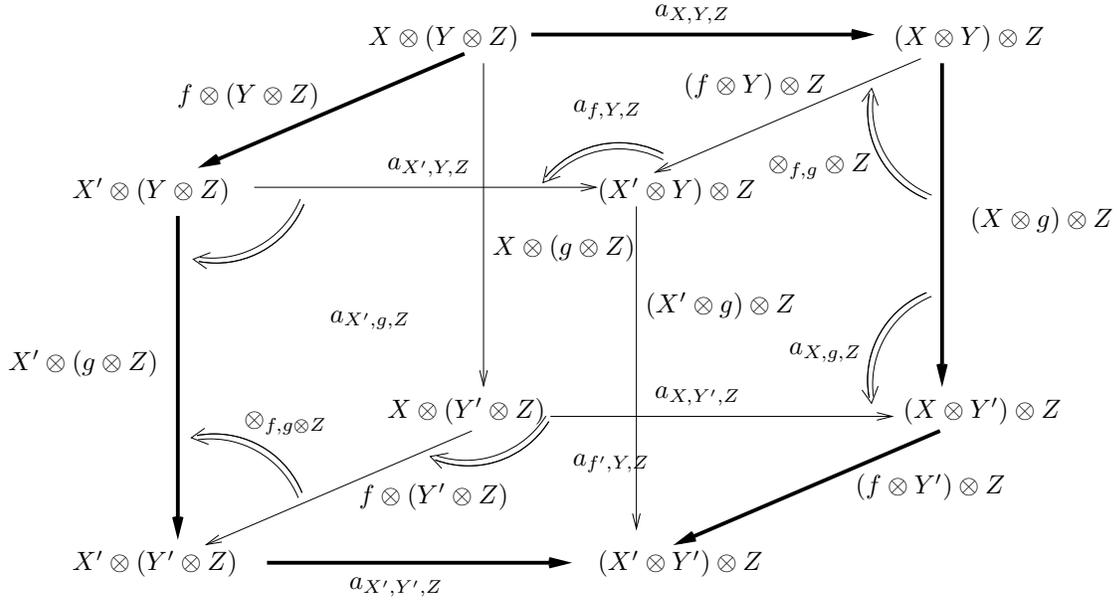}
\caption{Kapranov-Voevodsky's axiom $(\rightarrow\otimes\rightarrow\otimes
\bullet)$}
\label{cub}
\end{figure}

\begin{rem}
Let us remark that, apart from the structural 1- and 2-isomorphisms
related to the unital structure that do not appear in our
definition above, Kapranov-Voevodsky's definition of a monoidal
2-category (see \cite{KV94}) includes a different collection of
structural 2-isomorphisms. So, instead of our
$\widehat{\otimes}((f',g'),(f,g)):(f'\otimes g')\circ(f\otimes
g)\Longrightarrow (f'\circ f)\otimes (g'\circ g)$, they introduce
the two sets of 2-isomorphisms $\otimes_{f,f',Y}:(f'\otimes Y)\circ
(f\otimes Y)\Longrightarrow (f'\circ f)\otimes Y$ and
$\otimes_{X,g,g'}:(X\otimes g')\circ(X\otimes g)\Longrightarrow
X\otimes (g'\circ g)$ together with the basic 2-isomorphisms
$\otimes_{f,g}:(f\otimes Y')\circ(X\otimes
g)\Longrightarrow(X'\otimes g)\circ(f\otimes Y)$. Similarly,
instead of our $\widehat{a}(f,g,h):((f\otimes g)\otimes h)\circ
a_{X,Y,Z}\Longrightarrow a_{X',Y',Z'}\circ(f\otimes(g\otimes h))$,
they use 2-isomorphisms $a_{f,Y,Z}:((f\otimes Y)\otimes Z)\circ
a_{X,Y,Z}\Longrightarrow a_{X',Y',Z'}\circ(f\otimes(Y\otimes Z))$
and the similarly defined $a_{X,g,Z}$, $a_{X,Y,h}$. This obviously
implies a different set of axioms. However, both formulations are
equivalent, and correspond to the two possible ways of defining a
``bipseudofunctor'' directly as a pseudofunctor of two variables or
as two collections of pseudofunctors of one variable. Although it
is possible to work with Kapranov-Voevodsky's 2-isomorphisms, the
cohomological nature of the axioms is much more clear when working
with those of the previous proposition. Let us further remark that
the special case where the pseudofunctor
$\otimes:\text\Cg^2\To\text\Cg$ in our definition is {\it cubical}
(see the definition below) corresponds, in the Kapranov-Voevodsky's
formulation, to the notion of a {\it quasifunctor of 2-variables}
introduced by Gray in \cite{jG74}, p.56. The equivalence between
both notions, cubical pseudofunctor and quasifunctor of
2-variables, is in fact the content of Proposition I.4.8. in Gray's
book. We do not enter into the details of the equivalence, but let
us mention that, in terms of our
$\widehat{\otimes}((f',g'),(f,g))$, the above Kapranov-Voevodsky's
2-isomorphisms $\otimes_{f,g}$ correspond to
$$
\otimes_{f,g}=\widehat{\otimes}((id_{X'},g),(f,id_Y))^{-1}\cdot
\widehat{\otimes}((f,id_{Y'}),(id_X,g))
$$
and conversely, our $\widehat{\otimes}((f',g'),(f,g))$ are given by
$$
\widehat{\otimes}((f',g'),(f,g))=1_{f'\otimes Y''}\circ(\otimes_{f,g'})^{-1}\circ 1_{X\otimes g}
$$
The reader may also check that our structural equation
$(A\widehat{a}2)$ exactly corresponds to the axiom Kapranov and
Voevodsky denote by $(\rightarrow\otimes\rightarrow\otimes\bullet)$
(see Fig.~\ref{cub}) together with two more similar axioms.
\end{rem}

%\end{enumerate}
%

\subsection{}
A fundamental fact in the theory of semigroupal 2-categories is the
corresponding {\it strictification theorem}, due to
Gordon-Power-Street \cite{GPS95}. In fact, they proved a much more
general strictification theorem, valid for an arbitrary {\it
tricategory} (i.e., the categorification of the notion of a
bicategory). In the same way as a monoidal category just
corresponds to a bicategory of only one object, a monoidal
2-category is just a tricategory of only one object \footnote{
Strictly speaking, tricategories of one object correspond to the
more general notion of a monoidal {\it bicategory}.}. Now, contrary
to the case of bicategories, not all tricategories are equivalent
to the corresponding 3-categories (the reader may figure out the
precise definition of such objects). Indeed, some of the structural
3-isomorphisms can not be strictified in general, i.e., made equal
to identitites. In our case, this means that an arbitrary
semigroupal 2-category is in some sense equivalent to a particular
kind of semigroupal 2-categories, which, following Day and Street
\cite{DS97}, we will call {\it Gray semigroups}
\footnote{ Actually, they use the name {\it Gray monoid}, because
they consider monoidal 2-categories.}, and which are not the strict
semigroupal 2-categories. Since this theorem plays an essential
role in what follows, allowing us to greatly simplify the theory,
we review here the precise definitions.

\begin{defn}
Let \Cg\ be any 2-category. A pseudofunctor
$\F:\text\Cg\times\text\Cg\To\text\Cg$ is called {\sl cubical} if
its structural 2-isomorphisms
$\widehat{\F}((f',g'),(f,g)):\F(f',g')\circ\F(f,g)\Longrightarrow
\F(f'\circ f,g'\circ g)$ and
$\F_0(X,Y):\F(id_X,id_Y)\Longrightarrow id_{\F(X,Y)}$ are such
that:
\begin{enumerate}
\item
The $\widehat{\F}((f',id_{Y'}),(f,g))$ and
$\widehat{\F}((f,g'),(id_X,g))$ are all identity 2-morphisms.

\item
$\F$ is a unitary pseudofunctor, i.e., for all $(X,Y)$,
$\F_0(X,Y)=1_{id_{\F(X,Y)}}$.
\end{enumerate}
\end{defn}

Notice that our definition here differs from that in \cite{GPS95},
p.31, in that we explicitly require the pseudofunctor to be
unitary. Indeed, although the authors say that this condition
follows from the cubical condition of $\F$, it seems that this is
not the case, and the assumption must be included in the definition
\footnote{We would like to thank J. Power and R. Street for the emails 
interchanged about this point.}.

\begin{defn}
A {\sl cubical semigroupal 2-category} is any semigroupal
2-category $(\text\Cg,\otimes,a,\pi)$ such that the tensor product
$\otimes$ is a cubical pseudofunctor. A cubical semigroupal
2-category will be called a {\sl Gray semigroup} whenever its
structural 2-isomorphisms included in the associator $a$ and the
pentagonator $\pi$ are all identities.
\end{defn}
\begin{rem}
The analogous notions in the more general context of tricategories
are respectively called {\it cubical tricategories} and {\it Gray
categories} in \cite{GPS95}.
\end{rem}

A Gray semigroup will be simply denoted by $(\text\Cg,\otimes)$,
the $a$ and $\pi$ being trivial. We leave to the reader to make
explicit this definition. Notice that the set of structural
2-isomorphisms reduces in this case to the
$\widehat{\otimes}((f',g'),(f,g))$, most of which are moreover
trivial by the cubical condition. This is the reason a Gray
semigroup is usually described in terms of Kapranov-Voevodsky's
2-isomorphisms $\otimes_{f,g}$. It is worth to point out that not
every cubical pseudofunctor $\otimes$ defines a structure of Gray
semigroup on a 2-category.

The fundamental strictification theorem for semigroupal
2-categories can now be stated as follows:

\begin{thm} (\cite{GPS95}) \label{strictification}
Every semigroupal 2-category is equivalent (in a sense we do not
make precise) to a Gray semigroup.
\end{thm}

After reading the next section, where the notion of a morphism
between semigroupal 2-categories is defined, the reader may figure
out by himself the sense in which this equivalence should be
understood.

\subsection{}
Let us finish this section by giving the corresponding notion of
morphisms between semigroupal 2-categories, which will be needed in
the next section in order to define equivalence of deformations. In
the case of Gray semigroups, this definition appears, for example,
in \cite{DS97}, Def.2 (in fact, they define morphism between Gray
monoids). Our definition below follows from the general definition
of morphism between tricategories which appears in \cite{GPS95}
when restricted to the one object case (and forgetting the unital
structure).
%Furthermore, if $f:X\To X'$ and $g:Y\To
%Y'$ are any 1-morphisms of a semigroupal 2-category, we will denote
%the composition $(f\otimes Y')\circ (X\otimes g)$ simply by
%$f\otimes g$.

%
\begin{defn}
Let $(\text\Cg,\otimes,a,\pi)$ and $(\text\Cg',\otimes',a',\pi')$
be semigroupal 2-categories. A {\sl semigroupal pseudofunctor} from
\Cg\ to
$\text\Cg'$ is a pseudofunctor $\F:\text\Cg\To\text\Cg'$ together
with the following data {\bf SPDi} and axiom {\bf SPA}:

\begin{description}
\item[SPD1]
A pseudonatural isomorphism $\psi:\otimes'\circ
(\F\times\F)\Longrightarrow
\F\circ\otimes:\text\Cg\times\text\Cg\to\text\Cg'$.

\item[SPD2]
An invertible modification $\omega:(1_{\F}\circ a)\cdot(\psi\circ
1_{id_{\text\Cg}\times\otimes})\cdot(1_{\otimes'}\circ(1_{\F}\times\psi))\Rightarrow(\psi\circ
1_{\otimes\times id_{\text\Cg}})\cdot(1_{\otimes'}\circ(\psi\times
1_{\F}))\cdot(a'\circ 1_{\F\times\F\times\F})$ (see
Fig.~\ref{omega_figura})

\begin{figure}
\centering
\input{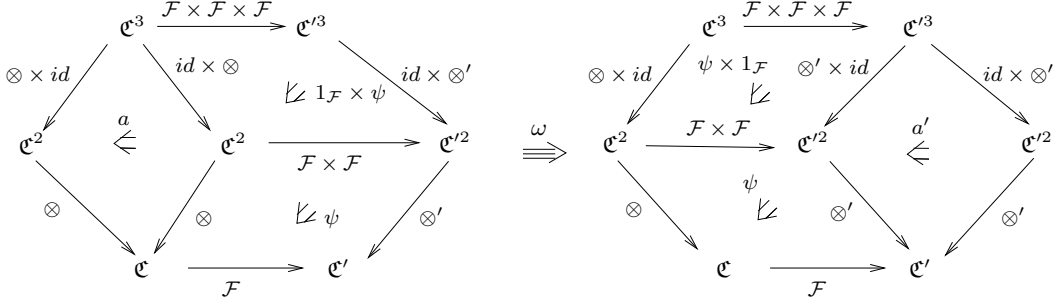}
\caption{Definition of the modification $\omega$}
\label{omega_figura}
\end{figure}

\item[SPA]
The pair $(\psi,\omega)$ is such that the equation in
Fig.~\ref{semigrupal_figura} holds (to simplify notation, the
tensor product of objects and 1-morphisms is again denoted by
simple juxtaposition and the identity 1-morphisms are represented
by the corresponding objects; furthermore, the action of the
pseudofunctor on objects, 1-morphisms or 2-morphisms is indicated
by the symbols $[-]$, so that, for example, $\psi_{XY,Z}([Z][T])$
denotes the 1-morphism $\psi_{X\otimes
Y,Z}\otimes'(\F(id_Z)\otimes'\F(id_T))$. For more details about the
notations appearing in this Figure, see the next Proposition).
\end{description}

\begin{figure}
\centering
\input{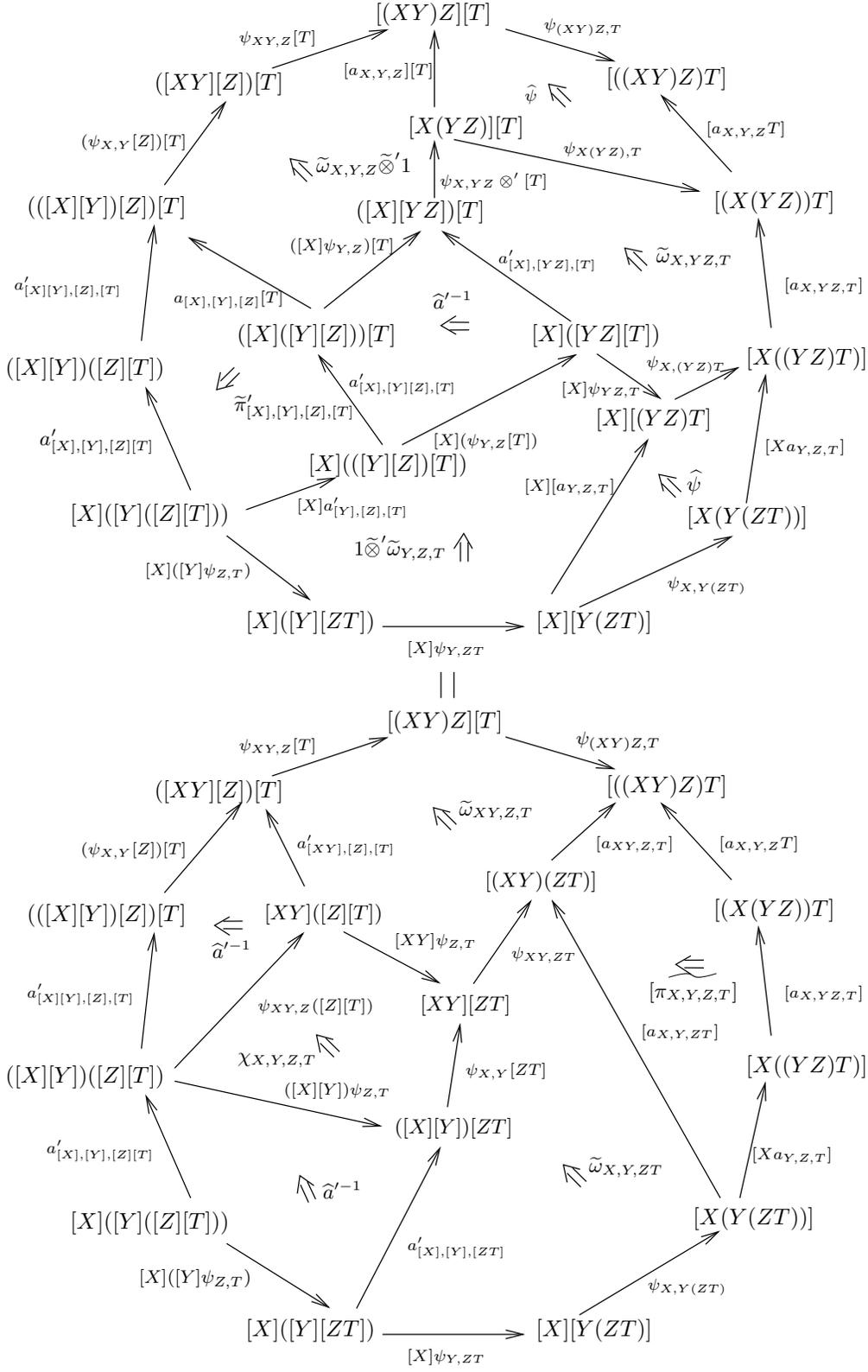}
\caption{The coherence relation $(SPA)$}
\label{semigrupal_figura}
\end{figure}

\noindent{A} semigroupal pseudofunctor will be denoted by the triple
$(\F,\psi,\omega)$ and the pair $(\psi,\omega)$ will be called a
{\sl semigroupal structure on} $\F$.
\end{defn}
Observe that the above definition indeed corresponds to
categorifying the definition of a semigroupal functor between
semigroupal categories: the axiom on the semigroupal structure is
substituted for the modification $\omega$, which in turn must
satisfy the additional coherence relation (SPA). A more explicit
description of a semigroupal structure $(\psi,\omega)$ on $\F$,
with the whole list of equations on the structural 1- and
2-isomorphisms, is as follows:

\begin{prop} \label{estructura_semigrupal}
Let $(\text\Cg,\otimes,a,\pi)$ and $(\text\Cg',\otimes',a',\pi')$
be semigroupal 2-categories and $\F:\text\Cg\To\text\Cg'$ a
pseudofunctor. Then, a semigroupal structure $(\psi,\omega)$ on
$\F:\text\Cg\To\text\Cg'$ consists of:

\begin{description}
\item[$\mathbf{\psi}$]
A collection of 1-isomorphisms $\psi_{X,Y}:\F(X)\otimes'\F(Y)\To
\F(X\otimes Y)$ for all objects $(X,Y)$ of $\text\Cg\times\text\Cg$.
\item[$\mathbf{\widehat{\psi}}$]
A collection of 2-isomorphisms $\widehat{\psi}(f,g):\F(f\otimes
g)\circ\psi_{X,Y}\Longrightarrow\psi_{X',Y'}\circ(\F(f)\otimes'\F(g))$
for all 1-morphisms $(f,g):(X,Y)\To (X',Y')$ of
$\text\Cg\times\text\Cg$.

\item[$\mathbf{\omega}$]
A collection of 2-isomorphisms
$\omega_{X,Y,Z}:\F(a_{X,Y,Z})\circ\psi_{X,Y\otimes
Z}\circ(id_{\F(X)}\otimes'\psi_{Y,Z})\Longrightarrow\psi_{X\otimes
Y,Z}\circ(\psi_{X,Y}\otimes' id_{\F(Z)})\circ
a'_{\F(X),\F(Y),\F(Z)}$ for all objects $(X,Y,Z)$ of
$\text\Cg\times\text\Cg\times\text\Cg$.

%$$
%\xymatrix{
%& \F(X)\otimes'(\F(Y)\otimes'\F(Z))
%\ar[ld]^{id_{\F(X)}\otimes'\F_2(Y,Z)} \ar[r]^{a'_{\F(X),\F(Y),\F(Z)}} &
%(\F(X)\otimes'\F(Y))\otimes'\F(Z) \ar[rd]^{\F_2(X,Y)\otimes'
%id_{\F(Z)}} & \\ \F(X)\otimes'\F(Y\otimes Z)
%\ar[rd]^{\F_2(X,Y\otimes z)} & & & \F(X\otimes Y)\otimes'\F(Z)
%\ar[ld]^{\F_2(X\otimes Y,Z)}
%\\ & \F(X\otimes(Y\otimes Z)) \ar[r]_{a_{X,Y,Z}} & \F((X\otimes Y)\otimes Z) &
%}
%$$
\end{description}
\noindent{Moreover}, these data must satisfy the following equations:

\begin{description}
\item[A$\mathbf{\widehat{\psi}}$1]
For all 2-morphisms
$(\tau,\sigma):(f,g)\Longrightarrow(\tilde{f},\tilde{g}):(X,Y)\To(X',Y')$
of $\text\Cg\times\text\Cg$
$$
(1_{\psi_{X',Y'}}\circ(\F(\tau)\otimes'\F(\sigma)))\cdot\widehat{\psi}(f,g)=
\widehat{\psi}(\tilde{f},\tilde{g})\cdot(\F(\tau\otimes\sigma)\circ 1_{\psi_{X,Y}})
$$
\item[A$\mathbf{\widehat{\psi}}$2]
For all composable 1-morphisms
$(X,Y)\stackrel{(f,g)}{\To}(X',Y')\stackrel{(f',g')}{\To}(X'',Y'')$
of $\text\Cg\times\text\Cg$
\begin{eqnarray*}
\lefteqn{\widehat{\psi}(f'\circ
f,g'\circ
g)\cdot([\F(\widehat{\otimes}((f',g'),(f,g)))\cdot\widehat{\F}(f'\otimes
g',f\otimes g)]\circ 1_{\psi_{X,Y}})=} \\ &
&=(1_{\psi_{X'',Y''}}\circ[(\widehat{\F}(f',f)\otimes'\widehat{\F}(g',g))\cdot
\widehat{\otimes}'((\F(f'),\F(g')),(\F(f),\F(g)))])\cdot \\ & &
\ \ \ \ \cdot(\widehat{\psi}(f',g')\circ
1_{\F(f)\otimes'\F(g)})\cdot(1_{\F(f'\otimes
g')}\circ\widehat{\psi}(f,g))
\end{eqnarray*}
\item[A$\mathbf{\widehat{\psi}}$3]
For all objects $(X,Y)$ of $\text\Cg\times\text\Cg$
\begin{align*}
\widehat{\psi}(id_X,id_Y)=(1_{\psi_{X,Y}}&\circ
[\otimes'_0(\F(X),\F(Y))\cdot(\F_0(X)\otimes'\F_0(Y))])^{-1}\cdot
\\ &\cdot([\F_0(X\otimes Y)\cdot\F(\otimes_0(X,Y))]\circ 1_{\psi_{X,Y}})
\end{align*}

\item[A$\mathbf{\omega}$1]
For all 1-morphisms $(f,g,h):(X,Y,Z)\To(X',Y',Z')$ of
$\text\Cg\times\text\Cg\times\text\Cg$
\begin{align*}
(\omega_{X',Y',Z'}&\circ
1_{\F(f)\otimes'(\F(g)\otimes'\F(h))})\cdot
\\ &\cdot(1_{\F(a_{X',Y',Z'})\circ\psi_{X',Y'\otimes
Z'}}\circ\widehat{\otimes}'((id_{\F(X')},\psi_{Y',Z'}),(\F(f),\F(g)\otimes'\F(h)))^{-1})\cdot
\\ &\cdot(1_{\F(a_{X',Y',Z'})\circ\psi_{X',Y'\otimes
Z'}}\circ(1_{\F(f)}\otimes'\widehat{\psi}(g,h)))\cdot
\\ &\cdot(1_{\F(a_{X',Y',Z'})\circ\psi_{X',Y'\otimes
Z'}}\circ\widehat{\otimes}'((\F(f),\F(g\otimes
h)),(id_{\F(X)},\psi_{Y,Z})))\cdot
\\ &\cdot(1_{\F(a_{X',Y',Z'})}\circ\widehat{\psi}(f,g\otimes h)
\circ 1_{id_{\F(X)}\otimes'\psi_{Y,Z}})\cdot
\\ &\cdot(\widehat{\F}(a_{X',Y',Z'},f\otimes(g\otimes
h))^{-1}\circ 1_{\psi_{X,Y\otimes Z}}\circ
1_{id_{\F(X)}\otimes'\psi_{Y,Z}})\cdot
\\ &\cdot(\F(\widehat{a}(f,g,h))\circ
1_{\psi_{X,Y\otimes Z}}\circ
1_{id_{\F(X)}\otimes'\psi_{Y,Z}})\cdot
\\ &\cdot(\widehat{\F}((f\otimes
g)\otimes h,a_{X,Y,Z})\circ 1_{\psi_{X,Y\otimes Z}}\circ
1_{id_{\F(X)}\otimes'\psi_{Y,Z}})=
\\ &\ \ \ \ \ =(1_{\psi_{X'\otimes Y',Z'}}\circ 1_{\psi_{X',Y'}\otimes'
id_{\F(Z')}}\circ\widehat{a}'(\F(f),\F(g),\F(h)))\cdot
\\ &\ \ \ \ \ \ \ \ \ \cdot(1_{\psi_{X'\otimes Y',Z'}}\circ
\widehat{\otimes}'((\psi_{X',Y'},id_{\F(Z')}),(\F(f)\otimes'\F(g),\F(h)))^{-1}
\circ 1_{a'_{\F(X),\F(Y),\F(Z)}})\cdot
\\ &\ \ \ \ \ \ \ \ \ \cdot(1_{\psi_{X'\otimes Y',Z'}}\circ(\widehat{\psi}(f,g)\otimes'
1_{\F(h)})\circ 1_{a'_{\F(X),\F(Y),\F(Z)}})\cdot
\\ &\ \ \ \ \ \ \ \ \ \cdot(1_{\psi_{X'\otimes Y',Z'}}\circ
\widehat{\otimes}'((\F(f\otimes g),\F(h)),(\psi_{X,Y},id_{\F(Z)}))\circ
1_{a'_{\F(X),\F(Y),\F(Z)}})\cdot
\\ &\ \ \ \ \ \ \ \ \ \cdot(\widehat{\psi}(f\otimes
g,h)\circ 1_{(\psi_{X,Y}\otimes' id_{\F(Z)})\circ
a'_{\F(X),\F(Y),\F(Z)}})\cdot
\\ &\ \ \ \ \ \ \ \ \ \cdot(1_{\F((f\otimes g)\otimes
h)}\circ\omega_{X,Y,Z})
\end{align*}

\item[A$\mathbf{\omega}$2]
For all objects $(X,Y,Z,T)$ of
$\text\Cg\times\text\Cg\times\text\Cg\times\text\Cg$
\begin{align*}
(1&_{\psi_{(X\otimes Y)\otimes Z,T}\circ(\psi_{X\otimes
Y,Z}\otimes'\F(id_T))}\circ\widehat{a}'(\psi_{X,Y},\F(id_Z),\F(id_T))^{-1}\circ
1_{a'_{\F(X),\F(Y),\F(Z)\otimes'\F(T)}})\cdot
\\ &\cdot(\tilde{\omega}_{X\otimes Y,Z,T}\circ
1_{(\psi_{X,Y}\otimes'(\F(id_Z)\otimes'\F(id_T)))\circ
a'_{\F(X),\F(Y),\F(Z)\otimes'\F(T)}})\cdot
\\ &\cdot(1_{\F(a_{X\otimes Y,Z,T})\circ\psi_{X\otimes Y,Z\otimes T}}\circ\chi_{X,Y,Z,T}\circ
1_{a'_{\F(X),\F(Y),\F(Z)\otimes'\F(T)}})\cdot
\\ &\cdot(1_{\F(a_{X\otimes Y,Z,T})\circ\psi_{X\otimes Y,Z\otimes T}\circ
(\psi_{X,Y}\otimes'\F(id_{Z\otimes
T}))}\circ\widehat{a}'(\F(id_X),\F(id_Y),\psi_{Z,T}))^{-1}\cdot
\\ &\cdot(1_{\F(a_{X\otimes Y,Z,T})}\circ\tilde{\omega}_{X,Y,Z\otimes
T}\circ 1_{\F(id_X)\otimes'(\F(id_Y)\otimes'\psi_{Z,T})})\cdot
\\ &\cdot(\tilde{\F}(\pi_{X,Y,Z,T})\circ 1_{\psi_{X,Y\otimes(Z\otimes T)}\circ(\F(id_X)
\otimes'\psi_{Y,Z\otimes T})\circ(\F(id_X)\otimes'(\F(id_Y)\otimes'\psi_{Z,T}))})=
\\ &\ \ \ \ \ =(1_{\psi_{(X\otimes Y)\otimes Z,T}\circ(\psi_{X\otimes
Y,Z}\otimes'\F(id_T))\circ((\psi_{X,Y}\otimes'\F(id_Z))\otimes'\F(id_T))}
\circ\tilde{\pi}'_{\F(X),\F(Y),\F(Z),\F(T)})\cdot
\\ &\ \ \ \ \ \ \ \ \ \cdot(1_{\psi_{(X\otimes Y)\otimes Z,T}}\circ
(\tilde{\omega}_{X,Y,Z}\tilde{\otimes} 1_{\F(id_T)})\circ
1_{a'_{\F(X),\F(Y)\otimes'\F(Z),\F(T)}\circ(\F(id_X)\otimes'
a'_{\F(Y),\F(Z),\F(T)})})\cdot
\\ &\ \ \ \ \ \ \ \ \ \cdot(\widehat{\psi}(a_{X,Y,Z},id_T)\circ
1_{(\psi_{X,Y\otimes
Z}\otimes'[id_T])\circ(([id_X]\otimes'\psi_{Y,Z})\otimes'[id_T])\circ
a'_{[X],[Y]\otimes'[Z],[T]}\circ([id_X]\otimes'
a'_{[Y],[Z],[T]})})\cdot
\\ &\ \ \ \ \ \ \ \ \ \cdot(1_{[a_{X,Y,Z}\otimes id_T]\circ\psi_{X\otimes(Y\otimes Z),T}
\circ(\psi_{X,Y\otimes Z}\otimes'[id_T])}\circ\widehat{a}'([id_X],\psi_{Y,Z},[id_T])^{-1}
\circ 1_{[id_X]\otimes' a'_{[Y],[Z],[T]}})\cdot
\\ &\ \ \ \ \ \ \ \ \ \cdot(1_{\F(a_{X,Y,Z}\otimes id_T)}\circ\tilde{\omega}_{X,Y\otimes
Z,T}\circ
1_{(\F(id_X)\otimes'(\psi_{Y,Z}\otimes'\F(id_T)))\circ(\F(id_X)\otimes'
a'_{\F(Y),\F(Z),\F(T)})})\cdot
\\ &\ \ \ \ \ \ \ \ \ \cdot(1_{\F(a_{X,Y,Z}\otimes id_T)
\circ\F(a_{X,Y\otimes Z,T})\circ\psi_{X,(Y\otimes Z)\otimes T}}\circ
(1_{\F(id_X)}\tilde{\otimes}\tilde{\omega}_{Y,Z,T}))\cdot
\\ &\ \ \ \ \ \ \ \ \ \cdot(1_{[a_{X,Y,Z}\otimes id_T]
\circ[a_{X,Y\otimes Z,T}]}\circ\widehat{\psi}(id_X,a_{Y,Z,T})\circ
1_{([id_X]\otimes'\psi_{Y,Z\otimes
T})\circ([id_X]\otimes'([id_Y]\otimes'\psi_{Z,T}))})
\end{align*}
(to simplify notation, $\F(-)$ is denoted in some places by $[-]$).
\end{description}
\end{prop}

\noindent{In} the last equation (A$\omega$2), the term $\chi_{X,Y,Z,T}$
denotes the 2-isomorphism
\begin{align*}
\chi_{X,Y,Z,T}=&\widehat{\otimes}((\F(id_{X\otimes Y}),\psi_{Z,T}),
(\psi_{X,Y},\F(id_Z)\otimes'\F(id_T)))^{-1}\cdot((\F(\otimes_0(X,Y))\circ
1)\otimes' 1)\cdot \\
&\cdot(\widehat{\psi}(id_X,id_Y)^{-1}\otimes'\widehat{\psi}(id_Z,id_T))\cdot
(1\otimes'(\F(\otimes_0(Z,T)^{-1})\circ 1))\cdot
\\ &\cdot\widehat{\otimes}'((\psi_{X,Y},\F(id_{Z\otimes
T})),(\F(id_X)\otimes'\F(id_Y),\psi_{Z,T}))
\end{align*}
Notice that, in the particular case of a Gray semigroup, they just
reduce to Kapranov-Voevodsky's 2-isomorphisms
$\otimes_{\psi_{X,Y},\psi_{Z,T}}$, as it appears in Day-Street's
definition mentioned above \cite{DS97}. Observe also that all the
terms $\tilde{\omega}_{X,Y,Z}$, $\tilde{\pi}_{X,Y,Z,T}$,
$\tilde{\F}(\pi_{X,Y,Z,T})$,
$\tilde{\omega}_{X,Y,Z}\tilde{\otimes}' 1_{\F(id_T)}$
$1_{\F(id_X)}\tilde{\otimes}'\tilde{\omega}_{Y,Z,T}$ are pastings
of the corresponding 2-isomorphisms with the appropriate structural
2-isomorphisms from $\widehat{\F}$, $\F_0$ and
$\widehat{\otimes}'$. We leave to the reader to find out the
explicit formulas.

\begin{proof}
The Proposition again follows from the definitions in Section 2. In
particular, Equation (A$\omega$1) corresponds to the naturality
condition on $\omega_{X,Y,Z}$ in the object $(X,Y,Z)$, and
(A$\omega$2) is the algebraic expression of the coherence relation
(SPA).
\end{proof}

Later on, we will need this Proposition in the very special case
where $\text\Cg'=\text\Cg$ (but with different semigroupal
structures $(\otimes,a,\pi)$ and $(\otimes',a',\pi')$) and
$\F=id_{\text\Cg}$, the identity 2-functor of \Cg.

\section{Deformations of pseudofunctors and semigroupal 2-categories}

\subsection{}
In this section we formalize the idea outlined in the introduction,
i.e., we ``linearize'' the problem of deforming a semigroupal
2-category $(\text\Cg,\otimes,a,\pi)$. To do that, we will need to
assume that
\Cg\ has some $K$-linear structure, for some commutative ring with
unit $K$. Before that, however, we introduce the notion of a purely
pseudofunctorial infinitesimal deformation of a ($K$-linear)
pseudofunctor, a notion which appears later in Section 8 when we
study the deformations of the tensor product in a semigroupal
2-category. The corresponding notions of equivalent deformations
are also introduced, and they are made explicit in the case of
first order deformations for later use.

\subsection{}
Recall that a category $\C$ is called $K$-linear when all its
hom-sets $\C(A,B)$, $A,B\in|\C|$, are $K$-modules and all
composition maps are $K$-bilinear. On the other hand, a $K$-linear
functor between two $K$-linear categories $\C, \D$ is any functor
$F:\C\longrightarrow\D$ such that all maps
$F_{A,B}:\C(A,B)\longrightarrow\D(F(A),F(B))$, $A,B\in|\C|$, are
$K$-linear. The analogous definitions for 2-categories are as
follows.

\begin{defn}
Let $K$ a commutative ring with unit. A $K$-{\sl linear 2-category}
is a 2-cateogry \Cg\ such that all its hom-categories \Cg$(X,Y)$
are $K$-linear, and all the composition functors
$c_{X,Y,Z}:\text\Cg(X,Y)\times
\text\Cg(X,Z)\longrightarrow \text\Cg(X,Z)$ are $K$-bilinear.
Given $K$-linear 2-categories \Cg\ and \Dg, a $K$-{\sl linear
pseudofunctor} between them is any pseudofunctor
$\F:\text\Cg\longrightarrow \text\Dg$ such that all functors
$\F_{X,Y}:\text\Cg(X,Y)\To\text\Dg(\F(X),\F(Y))$,
$X,Y\in|\text\Cg|$, are $K$-linear.
\end{defn}

Notice that, according to this definition, we only have a structure of
$K$-module on the sets of 2-morphisms. This will mean that, in our
definition of deformation below, all structural 1-morphisms will
remain undeformed, and the only thing susceptible to be deformed
will be the 2-morphisms.

The following result brings together some easy facts about
$K$-linear 2-categories and pseudofunctors whose proof is left to
the reader.

\begin{prop} \label{prop_K_linear}
Let \Bg, \Cg\ and \Dg\ be $K$-linear 2-categories, and
$\F,\F':\text\Bg\To\text\Cg$ and $\G:\text\Cg\To\text\Dg$
$K$-linear pseudofunctors. Then:

(i) The product 2-category $\text\Bg\times\text\Cg$ is $K$-linear.

(ii) The composition pseudofunctor $\G\circ\F:\text\Bg\To\text\Dg$
is $K$-linear.

(iii) For any $\xi,\zeta:\F\Longrightarrow\F'$, the set ${\rm
Mod}(\xi,\zeta)$ (resp. ${\rm PseudMod}(\xi,\zeta)$) of
modifications (resp. pseudomodifications) between $\xi$ and $\zeta$
is a $K$-vector space.
\end{prop}
Our main objects of interest are the $K$-linear semigroupal
2-categories, defined as follows:

\begin{defn}
A $K$-{\sl linear semigroupal 2-category} is a semigroupal
2-category $(\text\Cg,\otimes,a,\pi)$ such that both \Cg\ and
$\otimes$ are $K$-linear.
\end{defn}
%On the
%other hand, given two $K$-linear semigroupal 2-categories
%$(\text\Cg,\otimes,a,\pi)$ and $(\text\Cg',\otimes',a',\pi')$, a
%$K$-{\sl linear semigroupal pseudofunctor} between them is any
%semigroupal pseudofunctor $(\F,\psi,\omega)$ such that the
%pseudofunctor $\F:\text\Cg\To\text\Cg'$ is $K$-linear.

\subsection{}
Fundamental for the definitions of infinitesimal deformation given
later are the notions of $R$-linear extension for $K$-linear
(semigroupal) 2-categories and $K$-linear pseudofunctors, for any
$K$-algebra $R$. As the reader will see, in the first case it provides 
us with the necessary ``tangent space'' at the point of $X(\text\Cg)$ defined 
by the semigroupal structure in question.

\begin{defn} \label{extensio_categoria}
Let \Cg\ be a $K$-linear 2-category. Given a $K$-algebra $R$, the
$R$-{\sl linear extension of} \Cg\ is the $R$-linear 2-category
$\text\Cg^0_R$ defined as follows: (1) its objects and 1-morphisms
are the same as in
\Cg, (2) its sets of 2-morphisms are given by
$(\text\Cg^0_R(X,Y))(f,f'):=(\text\Cg(X,Y))(f,f')\otimes_K R$, (3)
the vertical composition is defined by $(\tau\otimes
r)\cdot(\tilde{\tau}\otimes\tilde{r}):=(\tau\cdot\tilde{\tau})\otimes(r\tilde{r})$
and by linear extension, (4) the composition functors
$(c_{X,Y,Z})_R\equiv
\circ_R:\text\Cg^0_R(X,Y)\times\text\Cg^0_R(Y,Z)\To\text\Cg^0_R(X,Z)$
are defined on 1-morphisms as in \Cg\ and on 2-morphisms by
$(\eta\otimes s)\circ_R(\tau\otimes r):=(\eta\circ\tau)\otimes(rs)$
and by linear extension, and (5) the identity 1-morphisms are the
same as in $\text\Cg$.
\end{defn}
The reader may easily check that these data indeed define an
$R$-linear 2-category. The reason to add the zero superscript in
$\text\Cg^0_R$ will be soon understood.

\begin{rem}
If $K$ is a topological ring and $R$ is an \mm-adically complete
local $K$-algebra (for example, $R=K[[h]]$), it can be defined the
\mm-adically complete $R$-linear extension of \Cg. This extension is the starting point
for the definition of the \mm-adically complete infinitesimal
deformations. In this work, however, we are mainly interested in
the non topological case, and we leave to the reader to figure out
the corresponding definitions in this topological setting.
\end{rem}

We are specially interested in the case
$R=K[\epsilon]/<\epsilon^{n+1}>$. The corresponding $R$-linear
extension will be denoted by $\text\Cg^0_{(n)}$. In this case, a
generic 2-morphism $\tau_{\epsilon}:f\Longrightarrow
f':X\longrightarrow Y$ in the linear extension can be written in
the form
\begin{equation} \label{tau_n}
\tau_{\epsilon}=\tau_0+\tau_1\epsilon+\cdots+\tau_n\epsilon^n
\end{equation}
where $\tau_0,\ldots,\tau_n\in\text\Cg(X,Y)(f,f')$.

The above definition is part of a functor of extension of scalars
for $K$-linear 2-categories, a fact which allows us to further
introduce the required $R$-linear extension of a $K$-linear
pseudofunctor.

\begin{prop} \label{extensio}
Let \Cg, \Dg\ be two $K$-linear 2-categories and $R$ a $K$-algebra.
Then, any $K$-linear pseudofunctor
$\F=(|\F|,\F_*,\widehat{\F}_*,\F_0)$ from
\Cg\ to \Dg\ extends to an $R$-linear pseudofunctor
$\F^0_R=(|\F|^0_R,(\F_*)^0_R,(\widehat{\F}_*)^0_R,(\F_0)^0_R)$ from
$\text\Cg^0_R$ to $\text\Dg^0_R$. Furthermore, if
$\xi=(\xi_*,\widehat{\xi}_*)$ is a pseudonatural transformation
between two $K$-linear pseudofunctors $\F$ and $\G$, it extends to
a pseudonatural transformation
$\xi^0_R=((\xi_*)^0_R,(\widehat{\xi}_*)^0_R))$ between the
$R$-linear extensions $\F^0_R$ and $\G^0_R$, and the same thing for
modifications between pseudonatural transformations.
\end{prop}
\begin{proof}
Take $|\F|^0_R=|\F|$, and for any pair of objects
$X,Y\in|\text\Cg|$, define the functor $(\F_{X,Y})^0_R$ as follows:
$(\F_{X,Y})^0_R(f)=\F_{X,Y}(f)$ for all 1-morphisms $f:X\To Y$, and
on 2-morphisms, take $(\F_{X,Y})^0_R(\tau\otimes
r)=\F_{X,Y}(\tau)\otimes r$ and extend by linearity. Finally,
define a pseudofunctorial structure on $\F^0_R$ by taking
$(\widehat{\F})^0_R(g,f)=\widehat{\F}(g,f)\otimes 1$ and
$(\F_0)^0_R(X)=\F_0(X)\otimes 1$ for all objects $X,Y,Z$ and
1-morphisms $f,g$. The rest of the proposition is proved similarly
and is left to the reader.
\end{proof}

As a by-product, we obtain the notion of $R$-linear extension for
$K$-linear semigroupal 2-categories. Indeed, we have:
\begin{cor} \label{extensio_semigrupal}
Let $(\text\Cg,\otimes,a,\pi)$ be a $K$-linear semigroupal
2-category. Then, for any $K$-algebra $R$, the extension
$\text\Cg^0_R$ inherits a structure $(\otimes^0_R,a^0_R,\pi^0_R)$
of $R$-linear semigroupal 2-category.
\end{cor}
\begin{proof}
Indeed, to give a semigroupal structure on a 2-category means to
give a pseudofunctor, a pseudonatural isomorphism and a
modification, and all of them can be extended according to the
previous Proposition. We leave to the reader to check that this
extensions satisfy the appropriate axioms.
\end{proof}

\subsection{}
We can now define the corresponding notions of infinitesimal
deformation. Let us begin with the case of a $K$-linear
pseudofunctor. According to Proposition~\ref{extensio}, given such
a pseudofunctor $\F:\text\Cg\To\text\Dg$, we have a ``copy'' of it
$\F_R^0:\text\Cg_R^0\To\text\Dg^0_R$ in the ``category of
$R$-linear pseudofunctors''. The reason to consider such a copy is
that the infinitesimal deformations of $\F$ will actually be,
strictly speaking, deformations of that copy, for some local
$K$-algebra $R$. The copy itself will be called the {\sl null
deformation of} $\F$ {\sl over} $R$. A generic deformation is then
defined as follows.

\begin{defn}
Let \Cg, \Dg\ be two $K$-linear 2-categories, and
$\F=(|\F|,\F_*,\widehat{\F}_*,\F_0)$ a $K$-linear pseudofunctor
between them. Given a local $K$-algebra $R$, a {\sl purely
pseudofunctorial infinitesimal deformation of} $\F$ {\sl over} $R$
is the pair $(|\F|^0_R,(\F_*)^0_R)$ of Proposition~\ref{extensio}
equipped with a pseudofunctorial structure
$((\widehat{\F}_*)_R,(\F_0)_R)$ which reduces mod. \mm\ to that of
the null deformation. When $R=K[\epsilon]/<\epsilon^{n+1}>$, the
corresponding deformations are called {\sl purely pseudofunctorial}
$n^{th}$-{\sl order deformations of} $\F$.
\end{defn}
The terms ``purely pseudofunctorial'' in this definition refer to
the fact that the only deformed thing is the pseudofunctorial
structure of $\F_R^0$, the source and target 2-categories remaining
undeformed, in the sense that they are simply substituted for the
corresponding $R$-linear extensions.

For example, it is easy to see that to give a purely
pseudofunctorial $n^{th}$-order deformation of $\F$ simply amounts
to give new families of 2-isomorphisms of the form
\begin{eqnarray*}
&\widehat{\F}_{\epsilon}(g,f)=\widehat{\F}(g,f)+\widehat{\F}^{(1)}(g,f)\epsilon+\cdots+
\widehat{\F}^{(n)}(g,f)\epsilon^n \\
&(\F_0)_{\epsilon}(X)=\F_0(X)+\F_0^{(1)}(X)\epsilon+\cdots+\F_0^{(n)}(X)\epsilon^n
\end{eqnarray*}
where $\widehat{\F}^{(i)}(g,f):\F(g)\circ
\F(f)\Longrightarrow\F(g\circ f)$ and $\F_0^{(i)}(X):\F(id_X)\Longrightarrow id_{\F(X)}$,
for all $i=1,\ldots,n$, are suitable 2-morphisms in \Dg\ such that
the above 2-isomorphisms indeed define a pseudofunctorial structure
on the pair $(|\F|_R^0,(\F_*)_R^0)$. To emphasize that, a purely
pseudofunctorial $n^{th}$-order deformation will be denoted by the
pair
$(\{\widehat{\F}^{(i)}\}_{i=1,\ldots,n},\{\F_0^{(i)}\}_{i=1,\ldots,n})$.
In particular, when all these 2-morphisms are zero, we recover the
null deformation $\F^0_R$.

We are only interested in the equivalence classes of such purely
pseudofunctorial infinitesimal deformations, two such deformations
being considered equivalent in the following sense:
\begin{defn}
Let $\F=(|\F|,\F_*,\widehat{\F}_*,\F_0)$ be a $K$-linear
pseudofunctor. Then, given two purely pseudofunctorial
infinitesimal deformations $\F_R=((\widehat{\F}_*)_R,(\F_0)_R)$ and
$\F'_R=((\widehat{\F}_*)'_R,(\F_0)'_R)$, they are called {\sl
equivalent} if there exists a pseudonatural isomorphism
$\xi:\F_R\Longrightarrow\F'_R$ such that
\begin{enumerate}
\item
$\xi_X=id_{\F(X)}$ for all objects $X\in|\text\Cg|$, and
\item
$\widehat{\xi}(f)=1_{\F(f)}\ \ (\text{mod. \mm})$, for all
1-morphisms $f$ of \Cg.
\end{enumerate}
\end{defn}

\noindent{For} later use, let us make explicit what this definition means in
the case of first order deformations.

\begin{prop} \label{equivalent_pseudo}
Two purely pseudofunctorial first order deformations
$\F_{\epsilon}$, $\F'_{\epsilon}$ of a $K$-linear pseudofunctor
$\F:\text\Cg\To\text\Dg$, defined by 2-morphisms
$(\widehat{\F}^{(1)}(g,f),\F_0^{(1)}(X))$ and
$((\widehat{\F}^{(1)})'(g,f),(\F_0^{(1)})'(X))$, respectively, are
equivalent if and only if there exists 2-morphisms
$\widehat{\xi}^{(1)}(f):\F(f)\Longrightarrow\F(f)$, for all
1-morphisms $f$ of \Cg, satisfying the following conditions:
\begin{enumerate}
\item
They are natural in $f$.
%, i.e., for all 2-morphisms
%$\tau:f\Longrightarrow g$, we have
%$$
%\F(\tau)\cdot\widehat{\xi}^{(1)}(f)=\widehat{\xi}^{(1)}(g)\cdot\F(\tau)
%$$
\item
For all composable 1-morphisms
$X\stackrel{f}{\To}Y\stackrel{g}{\To}Z$, it holds
\begin{align*}
(\widehat{\F}^{(1)})'(g,f)-\widehat{\F}^{(1)}(g,f)=\widehat{\F}&(g,f)
\cdot(1_{\F(g)}\circ\widehat{\xi}^{(1)}(f))-\widehat{\xi}^{(1)}(g\circ
f)\cdot\widehat{\F}(g,f)+ \\ &+\widehat{\F}(g,f)\cdot
(\widehat{\xi}^{(1)}(g)\circ 1_{\F(f)})
\end{align*}
\item
For all objects $X$ of \Cg, it holds
$$
(\F_0^{(1)})'(X)-\F_0^{(1)}(X)=\F_0(X)\cdot\widehat{\xi}^{(1)}(id_X)
$$
\end{enumerate}
\end{prop}
\begin{proof}
Indeed, let us go back to the definition of pseudonatural
transformation (see Definition~\ref{trans_quasi}) and take
$\F=\F_{\epsilon}$, $\G=\F_{\epsilon}'$, and $\xi$ defined by
$\xi_X=id_{\F(X)}$ and
$\widehat{\xi}(f)=1_{\F(f)}+\widehat{\xi}^{(1)}(f)\epsilon$. The
conditions above follow then by writing out the first order terms
in $\epsilon$ in every condition satisfied by $\xi$.
\end{proof}

\subsection{}
Let us consider now the deformations of a $K$-linear semigroupal
2-category $(\text\Cg,\otimes,a,\pi)$. As in the case of a
pseudofunctor, the first thing we need is its ``copy''
$(\text\Cg^0_R,\otimes^0_R,a^0_R,\pi^0_R)$ in the ``category of
$R$-linear semigroupal 2-categories'' (see
Corollary~\ref{extensio_semigrupal}). It will be called the {\sl
null deformation of} $(\text\Cg,\otimes,a,\pi)$ {\sl over} $R$. A
generic infinitesimal deformation is then a deformation of that
copy. More precisely:

\begin{defn} \label{deformacio_2_categoria_semigrupal}
Let $(\text\Cg,\otimes,a,\pi)$ be a $K$-linear semigroupal
2-category, and let $R$ be a local $K$-algebra, with maximal ideal
\mm. An {\sl infinitesimal deformation of} $(\text\Cg,\otimes,a,\pi)$ {\sl over} $R$ is the
$R$-linear extension 2-category $\text\Cg^0_R$ equipped with a
semigroupal structure $(\otimes_R,a_R,\pi_R)$ which reduces mod.
\mm\ to that of the null deformation. More explicitly,
$(\otimes_R,a_R,\pi_R)$ must be such that:
(1) $\otimes_R$ only differs from $\otimes^0_R$ in the
pseudofunctorial structure, (2) the 1-isomorphisms $(a_R)_{X,Y,Z}$
coincide with those of $a_R^0$, and (3) all structural
2-isomorphisms
$\widehat{\otimes}_R((f',g'),(f,g)),(\otimes_0)_R(X,Y),
\widehat{a}_R(f,g,h),(\pi_R)_{X,Y,Z,T}$
reduce mod. \mm\ to those of the null deformation.

When $R=K[\epsilon]/<\epsilon^{n+1}>$, the corresponding
infinitesimal deformations are called $n^{th}$-{\sl order
deformations of} \Cg.
\end{defn}

\begin{rem}
Since the associator as well as the left and right unit constraints
on the composition of 1-morphisms (which we are assuming trivial)
both live in the ``deformable world'' of 2-morphisms, it would be
possible to modify the above definitions of infinitesimal
deformation in such a way that also the bicategory structure of
\Cg\ is deformed. In the case of a pseudofunctor, this will lead us to the notion of a
not necessarily purely pseudofunctorial deformation. However, we
will not proceed in this direction, and we will assume that all
bicategories, the undeformed as well as the deformed ones, are
always 2-categories.
\end{rem}

For example, according to the above definition, an arbitrary
$n^{th}$-order deformation of $(\text\Cg,\otimes,a,\pi)$ amounts to
a new set of structural 2-isomorphism
$\widehat{\otimes}_{\epsilon},(\otimes_0)_{\epsilon},\widehat{a}_{\epsilon},\pi_{\epsilon}$
of the form in Eq.(\ref{tau_n}) with the zero order term equal to
the original 2-isomorphism, i.e.,
\begin{eqnarray*}
&\widehat{\otimes}_{\epsilon}((f',g'),(f,g))=\widehat{\otimes}((f',g'),(f,g))+
\widehat{\otimes}^{(1)}((f',g'),(f,g))\epsilon+\cdots+\widehat{\otimes}^{(n)}((f',g'),(f,g))
\epsilon^n \\ &(\otimes_0)_{\epsilon}(X,Y)=\otimes_0(X,Y)+\otimes_0^{(1)}(X,Y)\epsilon+
\cdots+
\otimes_0^{(n)}(X,Y)\epsilon^n \\
&\widehat{a}_{\epsilon}(f,g,h)=\widehat{a}(f,g,h)+\widehat{a}^{(1)}(f,g,h)\epsilon+\cdots+
\widehat{a}^{(n)}(f,g,h)\epsilon^n \\ &(\pi_{\epsilon})_{X,Y,Z,T}=\pi_{X,Y,Z,T}+
(\pi^{(1)})_{X,Y,Z,T}\epsilon+\cdots+(\pi^{(n)})_{X,Y,Z,T}\epsilon^n
\end{eqnarray*}
where $\widehat{\otimes}^{(i)}((f',g'),(f,g))$,
$\otimes_0^{(i)}(X,Y)$, $\widehat{a}^{(i)}(f,g,h)$ and
$(\pi^{(i)})_{X,Y,Z,T}$, for all $i=1,\ldots,n$, are suitable
2-morphisms in \Cg\ with the same source and target 1-morphisms as
the corresponding undeformed 2-isomorphisms and such that the whole
set of new 2-isomorphisms satisfy all the necessary equations to
define a semigroupal structure on $\text\Cg^0_{(n)}$. Such a
$n^{th}$-order deformation will be denoted by
$(\{\widehat{\otimes}^{(i)}\}_i, \{\otimes_0^{(i)}\}_i,
\{\widehat{a}^{(i)}\}_i,\{\pi^{(i)}\}_i)$. In
particular, when all these 2-morphisms are zero we again recover
the null deformation.

As in the case of pseudofunctors, we are only interested in the
equivalence classes of infinitesimal deformations.

\begin{defn}
Let $(\text\Cg,\otimes,a,\pi)$ be a $K$-linear semigroupal
2-category. Two infinitesimal deformations over $R$
$(\otimes_R,a_R,\pi_R)$ and $(\otimes'_R,a'_R,\pi'_R)$ are called
{\sl equivalent} if the identity 2-functor
$id_{\text\Cg^0_R}:(\text\Cg^0_R,\otimes_R,a_R,\pi_R)\To(\text\Cg^0_R,\otimes'_R,a'_R,\pi'_R)$
admits a semigroupal structure $(\psi,\omega)$ such that:
\begin{enumerate}
\item
$\psi_{X,Y}=id_{X\otimes Y}$ for all objects $X,Y$;
\item
$\widehat{\psi}(f,g)=1_{f\otimes g}\ \ (\text{mod. \mm})$, for all
1-morphisms $f,g$, and
\item
$\omega_{X,Y,Z}=(\otimes_0(X\otimes Y,Z)^{-1}\circ
1_{a_{X,Y,Z}})\cdot(1_{a_{X,Y,Z}}\circ\otimes_0(X,Y\otimes Z))\ \
(\text{mod. \mm})$, for all objects $X,Y,Z$.
\end{enumerate}
The deformations will be called $\omega$-{\sl equivalent} when
there exists a semigroupal structure $(\psi,\omega)$ satisfying the
first and third conditions above and such that
$\widehat{\psi}(f,g)=1_{f\otimes g}$ for all $f,g$ (not only mod.
\mm). Similarly, the deformations will be called $\psi$-{\sl equivalent} when
there exists a semigroupal structure $(\psi,\omega)$ satisfying the
first and second conditions above and such that
$\omega_{X,Y,Z}=(\otimes_0(X\otimes Y,Z)^{-1}\circ
1_{a_{X,Y,Z}})\cdot(1_{a_{X,Y,Z}}\circ\otimes_0(X,Y\otimes Z))$ for
all $X,Y,Z$ (not only mod. \mm).
\end{defn}

Let us also make explicit for its later use what this definition
means in the case of first order deformations. To simplify
equations, however, let us assume, without loss of generality by
Theorem~\ref{strictification}, that the undeformed tensor product
$\otimes$ is unitary, i.e., that all 2-isomorphisms
$\otimes_0(X,Y)$ are identities (see
Definition~\ref{pseudofunctor}). Notice, however, that the deformed
tensor product may no longer be unitary.

\begin{prop} \label{equivalents}
Let $(\text\Cg,\otimes,a,\pi)$ be a $K$-linear semigroupal
2-category, with $\otimes$ a unitary tensor product. Let's consider
two first order deformations defined by 2-morphisms
$(\widehat{\otimes}^{(1)},\otimes_0^{(1)},\widehat{a}^{(1)},\pi^{(1)})$
and
$((\widehat{\otimes}^{(1)})',(\otimes_0^{(1)})',(\widehat{a}^{(1)})',(\pi^{(1)})')$.
Then, they are equivalent if and only there exists 2-morphisms
$\widehat{\psi}^{(1)}(f,g):f\otimes g\Longrightarrow f\otimes g$
and $(\omega^{(1)})_{X,Y,Z}:a_{X,Y,Z}\Longrightarrow a_{X,Y,Z}$,
for all objects $X,Y,Z$ and 1-morphisms $f,g$ of \Cg, such that the
following equations hold:
\begin{description}

\item[E$\widehat{\psi}$1]
For all 2-morphisms
$(\tau,\sigma):(f,g)\Longrightarrow(\tilde{f},\tilde{g})$ of
$\text\Cg^2$
$$
(\tau\otimes\sigma)\cdot\widehat{\psi}^{(1)}(f,g)=
\widehat{\psi}^{(1)}(\tilde{f},\tilde{g})\cdot(\tau\otimes\sigma)
$$

\item[E$\widehat{\psi}$2]
For all composable 1-morphisms
$(X,Y)\stackrel{(f,g)}{\To}(X',Y')\stackrel{(f',g')}{\To}(X'',Y'')$
of $\text\Cg^2$
\begin{align*}
\widehat{\psi}^{(1)}(f'\circ
f,g'\circ g)\cdot&\widehat{\otimes}((f',g'),(f,g)))+
\widehat{\otimes}^{(1)}(f',g'),(f,g))= \\
&=(\widehat{\otimes}^{(1)})'((f',g'),(f,g))+\widehat{\otimes}((f',g'),(f,g))
\cdot(\widehat{\psi}^{(1)}(f',g')\circ 1_{f\otimes g})+ \\ &\ \ \ \ +
\widehat{\otimes}((f',g'),(f,g))\cdot(1_{f'\otimes g'}\circ\widehat{\psi}^{(1)}(f,g))
\end{align*}
\item[E$\widehat{\psi}$3]
For all objects $(X,Y)$ of $\text\Cg\times\text\Cg$
$$
\widehat{\psi}^{(1)}(id_X,id_Y)=\otimes_0^{(1)}(X,Y)-
(\otimes_0^{(1)})'(X,Y)
$$

\item[E$\omega$1]
For all 1-morphisms $(f,g,h):(X,Y,Z)\To(X',Y',Z')$ of
$\text\Cg\times\text\Cg\times\text\Cg$
\begin{align*}
(\widehat{a}^{(1)})'(f,g,h)-&\widehat{a}(f,g,h)\cdot((\widehat{\otimes}^{(1)})'((id_{X'\otimes
Y'},id_{Z'}),(f\otimes g,h))\circ 1_{a_{X,Y,Z}})+
\\ & +\widehat{a}(f,g,h)\cdot((\widehat{\psi}^{(1)}(f,g)\otimes 1_h)\circ
1_{a_{X,Y,Z}})+
\\ &+\widehat{a}(f,g,h)\cdot((\widehat{\otimes}^{(1)})'((f\otimes
g,h),(id_{X\otimes Y},id_Z))\circ 1_{a_{X,Y,Z}})+\
\\ &+\widehat{a}(f,g,h)\cdot(\widehat{\psi}^{(1)}(f\otimes g,h)\circ
1_{a_{X,Y,Z}})+ \\ &+\widehat{a}(f,g,h)\cdot(1_{(f\otimes g)\otimes
h}\circ (\omega^{(1)})_{X,Y,Z})= \\ &\ \ \ \
=((\omega^{(1)})_{X',Y',Z'}\circ 1_{f\otimes(g\otimes h)})\cdot\widehat{a}(f,g,h)- \\
&\ \ \ \ \ \ \
-(1_{a_{X',Y',Z'}}\circ(\widehat{\otimes}^{(1)})'((id_{X',}id_{Y'\otimes
Z'}),(f,g\otimes h)))\cdot\widehat{a}(f,g,h)+
\\ &\ \ \ \ \ \ \ +(1_{a_{X',Y',Z'}}\circ(1_f\otimes\widehat{\psi}^{(1)}(g,h)))
\cdot\widehat{a}(f,g,h)+ \\ &\ \ \ \ \ \ \ +
(1_{a_{X',Y',Z'}}\circ(\widehat{\otimes}^{(1)})'((f,g\otimes h),
(id_X,id_{Y\otimes Z})))\cdot\widehat{a}(f,g,h)+ \\ &\ \ \ \ \ \ \
+ (1_{a_{X',Y',Z'}}\circ\widehat{\psi}^{(1)}(f,g\otimes
h))\cdot\widehat{a}(f,g,h)+\widehat{a}^{(1)}(f,g,h)
\end{align*}

\item[E$\omega$2]
For all objects $(X,Y,Z,T)$ of $\text\Cg^4$
\begin{align*}
(\pi^{(1)})'&_{X,Y,Z,T}+
\\ &+\pi_{X,Y,Z,T}\cdot((\otimes_0^{(1)})'((X\otimes Y)\otimes Z,T)\circ
1_{(a_{X,Y,Z}\otimes id_T)\circ a_{X,Y\otimes Z,T}\circ(id_X\otimes
a_{Y,Z,T})})+ \\
&+\pi_{X,Y,Z,T}\cdot(((\omega^{(1)})_{X,Y,Z}\otimes 1_{id_T})\circ
1_{a_{X,Y\otimes Z,T}\circ(id_X\otimes a_{Y,Z,T})})+ \\ &+
\pi_{X,Y,Z,T}\cdot(1_{a_{X,Y,Z}\otimes id_T}\circ(\otimes_0^{(1)})'(X\otimes(Y\otimes Z),T)
\circ 1_{a_{X,Y\otimes Z,T}\circ(id_X\otimes a_{Y,Z,T})})+ \\ &+\pi_{X,Y,Z,T}\cdot
(\widehat{\psi}^{(1)}(a_{X,Y,Z},id_T)
\circ 1_{a_{X,Y\otimes Z,T}\circ(id_X\otimes a_{Y,Z,T})})- \\ &-\pi_{X,Y,Z,T}\cdot
(1_{a_{X,Y,Z}\otimes
id_T}\circ(\widehat{a}^{(1)})'(id_X,id_{Y\otimes Z},id_T)\circ
1_{id_X\otimes a_{Y,Z,T}})+ \\ &+\pi_{X,Y,Z,T}\cdot
(1_{a_{X,Y,Z}\otimes id_T}\circ(\omega^{(1)})_{X,Y\otimes Z,T}\circ
1_{id_X\otimes a_{Y,Z,T}})+ \\
&+\pi_{X,Y,Z,T}\cdot(1_{(a_{X,Y,Z}\otimes id_T)\circ a_{X,Y\otimes
Z,T}}\circ(\otimes_0^{(1)})'(X,(Y\otimes Z)\otimes T)\circ
1_{id_X\otimes a_{Y,Z,T}})+ \\
&+\pi_{X,Y,Z,T}\cdot(1_{(a_{X,Y,Z}\otimes id_T)\circ a_{X,Y\otimes
Z,T}}\circ(1_{id_X}\otimes(\omega^{(1)})_{Y,Z,T}))+ \\ &+
\pi_{X,Y,Z,T}\cdot(1_{(a_{X,Y,Z}\otimes id_T)\circ a_{X,Y\otimes
Z,T}\circ(id_X\otimes a_{Y,Z,T})}\circ(\otimes_0^{(1)})'(X,Y\otimes
(Z\otimes T)))+  \\ &+\pi_{X,Y,Z,T}\cdot(1_{(a_{X,Y,Z}\otimes
id_T)\circ a_{X,Y\otimes
Z,T}}\circ\widehat{\psi}^{(1)}(id_X,a_{Y,Z,T}))= \\ &\ \ \ \ =
-(\widehat{a}^{(1)}(id_{X\otimes Y},id_Z,id_T)\circ 1_{a_{X,Y,Z\otimes T}})\cdot\pi_{X,Y,Z,T}+
\\ &\ \ \ \ \ \ +((\omega^{(1)})_{X\otimes Y,Z,T}\circ 1_{a_{X,Y,Z\otimes T}})
\cdot\pi_{X,Y,Z,T}- \\ &\ \ \ \ \ \ -(1_{a_{X\otimes Y,Z,T}}\circ
\widehat{\otimes}^{(1)}((id_{X\otimes Y},id_{Z\otimes T}),(id_{X\otimes Y},id_{Z\otimes T}))
\circ 1_{a_{X,Y,Z\otimes T}})\cdot\pi_{X,Y,Z,T}- \\
&\ \ \ \ \ \ -(1_{a_{X\otimes Y,Z,T}}\circ
(\widehat{\psi}^{(1)}(id_X,id_Y)\otimes 1_{id_{Z\otimes T}})
\circ 1_{a_{X,Y,Z\otimes T}})\cdot\pi_{X,Y,Z,T}+ \\
&\ \ \ \ \ \ +(1_{a_{X\otimes Y,Z,T}}\circ (1_{id_{X\otimes
Y}}\otimes\widehat{\psi}^{(1)}(id_Z,id_T))
\circ 1_{a_{X,Y,Z\otimes T}})\cdot\pi_{X,Y,Z,T}+ \\
&\ \ \ \ \ \ +(1_{a_{X\otimes Y,Z,T}}\circ
(\widehat{\otimes}^{(1)})'((id_{X\otimes Y},id_{Z\otimes
T}),(id_{X\otimes Y},id_{Z\otimes T}))
\circ 1_{a_{X,Y,Z\otimes T}})\cdot\pi_{X,Y,Z,T}- \\
&\ \ \ \ \ \ -(1_{a_{X\otimes
Y,Z,T}}\circ(\widehat{a}^{(1)})'(id_X,id_Y,id_{Z\otimes
T}))\cdot\pi_{X,Y,Z,T}+ \\ &\ \ \ \ \ \ + (1_{a_{X\otimes
Y,Z,T}}\circ(\omega^{(1)})_{X,Y,Z\otimes T})\cdot\pi_{X,Y,Z,T}+
\\ &\ \ \ \ \ \ +(\pi^{(1)})_{X,Y,Z,T}
\end{align*}
\end{description}
\end{prop}

\begin{proof}
The proof is a long but straightforward computation of the first
order term in each of the conditions in
Proposition~\ref{estructura_semigrupal} when $\F$ is taken equal to
the identity 2-functor $id_{\text\Cg^0_{(1)}}$, the semigroupal
structures $(\otimes,a,\pi)$ and $(\otimes',a',\pi')$ are those of
the first order deformations, and the $\psi$ and $\omega$ are of
the form
\begin{eqnarray*}
&\psi_{X,Y}=id_{X\otimes Y} \\ &\widehat{\psi}(f,g)=1_{f\otimes g}+\widehat{\psi}^{(1)}(f,g)\epsilon \\ &\omega_{X,Y,Z}=1_{a_{X,Y,Z}}+\omega^{(1)}_{X,Y,Z}\epsilon
\end{eqnarray*}
Notice that the zero order term of $\omega$ is trivial because we
are assuming $\otimes$ is unitary.
\end{proof}

\section{Cohomology of a unitary pseudofunctor}

\subsection{}
This section contains preliminary results that will be used in
Section 8 to construct the cochain complex which describes the
simultaneous deformations of both the tensor product and the
associator in a $K$-linear semigroupal 2-category. More explicitly,
we associate a cohomology to an arbitrary $K$-linear unitary
pseudofunctor and prove that this cohomology describes its purely
pseudofunctorial infinitesimal deformations in the sense of
Gerstenhaber. The main idea is to use the fact mentioned in Section
2 that a pseudofunctor between one object bicategories corresponds
to the notion of a monoidal functor. In this sense, our results
generalize the cohomology theory for monoidal functors described by
Yetter in \cite{dY98}. Let us remark that our restriction to the
case of unitary pseudofunctors implies no loss of generality for
our purposes, because, as indicated before, the results obtained
here will just be used to study the deformations of the tensor
product and the associator in a $K$-linear semigroupal 2-category.
Now, by Theorem~\ref{strictification}, the undeformed semigroupal
2-category may be assumed to be a Gray semigroup and, hence, such
that the original tensor product is indeed a unitary pseudofunctor.

\subsection{}
Let us consider a $K$-linear unitary pseudofunctor
$\F=(|\F|,\F_*,\widehat{\F}_*)$ between $K$-linear 2-categories
\Cg\ and \Dg. Since we assume
\Cg\ is a 2-category, for each $n\geq 2$ and each ordered family
$X_0,\ldots,X_n$ of $n+1$ objects of \Cg, we have a uniquely
induced composition functor
$c^{\text\Cg}_{X_0,\ldots,X_n}:\text\Cg(X_{n-1},X_n)\times\text
\Cg(X_{n-2},X_{n-1})\times\cdots\times\text\Cg(X_0,X_1)\longrightarrow\text\Cg(X_0,X_n)$,
obtained by applying the appropriate elementary composition
functors $c^{\text\Cg}_{X,Y,Z}$ in any order \footnote{Here, we
think of the elementary composition functors $c^{\text\Cg}_{X,Y,Z}$
as defined on the product category
$\text\Cg(Y,Z)\times\text\Cg(X,Y)$, instead of
$\text\Cg(X,Y)\times\text\Cg(Y,Z)$. Hence, they differ from those
appearing in Definition~\ref{bicategory} by a permutation
functor.}. In the same way, we have the induced composition
functors $c^{\text\Dg}_{\F(X_0),\ldots,\F(X_n)}$ for all
$X_0,\ldots,X_n$. Then, given $X_0,\ldots,X_n$, let's consider the
functors $C\F_{X_0,\ldots,X_n},\F
C_{X_0,\ldots,X_n}:\text\Cg(X_{n-1},X_n)\times\text
\Cg(X_{n-2},X_{n-1})\times\cdots\times\text\Cg(X_0,X_1)\To
\text\Dg(\F(X_0),\F(X_n))$ defined by
\begin{align*}
C\F_{X_0,\ldots,X_n}&:=c^{\text\Dg}_{\F(X_0),\ldots,\F(X_n)}
\circ(\F_{X_{n-1},X_n}\times\F_{X_{n-2},X_{n-1}}\times\cdots\times\F_{X_0,X_1})
\\
\F C_{X_0,\ldots,X_n}&:=\F_{X_0,X_n}\circ c^{\text\Cg}_{X_0,\ldots,X_n}
\end{align*}
When $n=1$, let $C\F_{X_0,X_1}:=\F_{X_0,X_1}=:\F C_{X_0,X_1}$.

We now define the vector spaces $X^n(\F)$ of the cochain complex we
are looking for as follows:
$$
X^n(\F):=\left\{
\begin{array}{ll} \prod_{X_0,\ldots,X_n\in|\text\Cg| }
\text{Nat}(C\F_{X_0,\ldots,X_n},\F C_{X_0,\ldots,X_n}) & n\geq 1 \\ 0 & \text{otherwise}
\end{array}\right.
$$
Notice that they are indeed vector spaces over $K$ because we are
assuming that the target 2-category \Dg\ is $K$-linear. According
to this definition, a generic element $\phi\in X^n(\F)$, $n\geq 1$,
is of the form $\phi=(\phi_{X_0,\ldots,X_n})_{X_0,\ldots,X_n}$,
with $\phi_{X_0,\ldots,X_n}=\{\phi(f_0,\ldots,f_{n-1})\ |\
f_i\in|\text\Cg(X_{n-i-1},X_{n-i})|,\ i=0,\ldots,n-1\}$ and
$$
\phi(f_0,\ldots,f_{n-1}):\F(f_0)\circ\F(f_1)\circ\cdots\circ\F(f_{n-1})\Longrightarrow
\F(f_0\circ f_1\circ\cdots\circ f_{n-1})
$$
a 2-morphism natural in $(f_0,\ldots,f_{n-1})$. On the other hand,
the ``padding'' composition operators introduced in Section 3
allows us to define coboundary maps
$\delta:X^{n-1}(\F)\longrightarrow X^n(\F)$, for all $n\geq 2$, in
the usual way. So, if $\phi\in X^{n-1}(\F)$, $\delta\phi\in
X^n(\F)$ is given by
\begin{align*}
(\delta\phi)(f_0,f_1,\ldots,f_{n-1})=&\lceil
1_{\F(f_0)}\circ\phi(f_1,\ldots,f_{n-1})
\rceil_{\F(X_0),\F(X_n)}+
\\ &+\sum_{i=1}^{n-1}(-1)^i\lceil\phi(f_0,\ldots,f_{i-1}\circ
f_i,\ldots,f_{n-1})\rceil_{\F(X_0),\F(X_n)}+
\\ &+(-1)^n\lceil\phi(f_0,\ldots,f_{n-2})\circ
1_{\F(f_{n-1})}\rceil_{\F(X_0),\F(X_n)}
\end{align*}
for all $f_i\in|\text\Cg(X_{n-i-1},X_{n-i})|$, $i=0,\ldots,n-1$.
\begin{prop} \label{complex}
For any $K$-linear unitary pseudofunctor
$\F=(|\F|,\F_*,\widehat{\F}_*)$, the pair
$(X^{\bullet}(\F),\delta)$ is a cochain complex.
\end{prop}
\begin{proof}
The map $\delta$ is clearly linear in $\phi$. On the other hand,
naturality of $\phi(f_1,\ldots,f_{n-1})$ in $(f_1,\ldots,f_{n-1})$
easily implies the naturality of each summand of
$\delta(\phi)(f_0,\ldots,f_{n-1})$ in $(f_0,\ldots,f_{n-1})$. For
example, making explicit the padding, the first term is
$$
\lceil 1_{\F(f_0)}\circ\phi(f_1,\ldots,f_{n-1})
\rceil_{\F(X_0),\F(X_n)}=\widehat{\F}(f_0,f_1\circ\ldots\circ f_{n-1})\cdot
(1_{\F(f_0)}\circ\phi(f_1,\ldots,f_{n-1}))
$$
Now, for any 2-morphism
$(\tau_0,\ldots,\tau_{n-1}):(f_0,\ldots,f_{n-1})\Longrightarrow(f'_0,\ldots,f'_{n-1})$,
the naturality of $\widehat{\F}$ and $\phi$ in its arguments and
the interchange law gives that
\begin{align*}
\F(\tau_0\circ\ldots\circ&\tau_{n-1})\cdot\lceil 1_{\F(f_0)}\circ\phi(f_1,\ldots,f_{n-1})
\rceil_{\F(X_0),\F(X_n)}=
\\ &=\F(\tau_0\circ\ldots\circ\tau_{n-1})\cdot
\widehat{\F}(f_0,f_1\circ\ldots\circ f_{n-1})\cdot
(1_{\F(f_0)}\circ\phi(f_1,\ldots,f_{n-1}))
\\ &=\widehat{\F}(f'_0,f'_1\circ\ldots\circ f'_{n-1})\cdot
(\F(\tau_0)\circ\F(\tau_1\circ\ldots\circ\tau_{n-1}))\cdot
(1_{\F(f_0)}\circ\phi(f_1,\ldots,f_{n-1}))
\\ &=\widehat{\F}(f'_0,f'_1\circ\ldots\circ f'_{n-1})\cdot
(\F(\tau_0)\circ[\F(\tau_1\circ\ldots\circ\tau_{n-1})\cdot
\phi(f_1,\ldots,f_{n-1})])
\\ &=\widehat{\F}(f'_0,f'_1\circ\ldots\circ f'_{n-1})\cdot
(\F(\tau_0)\circ[\phi(f'_1,\ldots,f'_{n-1})\cdot(\F(\tau_1)\circ\ldots\circ\F(\tau_{n-1}))])
\\ &=\widehat{\F}(f'_0,f'_1\circ\ldots\circ f'_{n-1})\cdot
(1_{\F(f'_0)}\circ\phi(f'_1,\ldots,f'_{n-1}))\cdot
(\F(\tau_0)\circ\F(\tau_1)\circ\ldots\circ\F(\tau_{n-1}))
\\ &=\lceil 1_{\F(f_0)}\circ\phi(f_1,\ldots,f_{n-1})
\rceil_{\F(X_0),\F(X_n)}\cdot
(\F(\tau_0)\circ\F(\tau_1)\circ\ldots\circ\F(\tau_{n-1}))
\end{align*}
The other terms are similarly worked. Finally, to prove that
$\delta^2=0$, notice first that, in computing $\delta^2(\phi)$, one
can initially forget the padding operators and take them into
account at the end of the computation. For example, the first term
of $\delta^2(\phi)$ reads
\begin{align*}
\delta(\lceil 1_{\F}\circ\phi\rceil)(f_0,\ldots,f_n)&=
\lceil 1_{\F(f_0)}\circ\lceil 1_{\F(f_1)}\circ\phi(f_2,\ldots,f_n)\rceil\rceil
\\ &\ \ \ \ \ \ \sum_{i=1}^n(-1)^i\lceil\lceil (1_{\F}\circ\phi)
(f_0,\ldots,f_{i-1}\circ f_i,\ldots,f_n)\rceil\rceil
\\ &\ \ \ \ \ \ (-1)^{n+1}\lceil\lceil 1_{\F(f_0)}\circ\phi(f_1,\ldots,f_{n-1})\rceil\circ
1_{\F(f_n)}\rceil
\end{align*}
Now, it is easy to check that the horizontal compositions of the
2-morphisms $1_{\F(f_0)}$ and $1_{\F(f_n)}$ in the first and last
terms commute with the padding. Our assertion follows then from the
obvious fact that taking a padding of a padding is the same as
doing nothing. So, let's provisionally forget the extra padding
operators in the computation of $\delta^2(\phi)$ and use the same
argument which shows the $\delta$ in the bar resolution satisfies
$\delta^2=0$ to deduce that the terms formally cancel out each
other. Reinserting now the padding operators in each summand of
this formal expression, corresponding terms still cancel out each
other because, by the coherence theorem, their paddings will also
coincide.
\end{proof}
This complex will be called the {\sl purely pseudofunctorial
deformation complex of} $\F$, and the corresponding cohomology will
be denoted by $H^{\bullet}(\F)$. Notice that the dependence of this
cohomology on the structural 2-isomorphisms $\widehat{\F}_*$ of
$\F$ is entirely encoded in the padding operators involved in the
definition of $\delta$.

\subsection{}
Let us suppose that both $\Cg$ and \Dg\ have only one object. Let
us denote by $X$ the only object of $\Cg$, so that the $\F(X)$ will
be the only object of \Dg. If we denote the (unique) composition
functor $c_{X,X,X}:\Cg(X,X)\times\Cg(X,X)\To\Cg(X,X)$ in $\Cg$ by
$\otimes^{\Cg}$ and in the same way denote by $\otimes^{\Dg}$ the
(unique) composition functor in \Dg, the purely pseudofunctorial
deformation complex of $\F$ clearly reduces to
$$
X^n(\F):=\left\{
\begin{array}{ll} {\rm Nat}((\otimes^{\Dg})^n\circ\F^n,\F\circ(\otimes^{\Cg})^n) & n\geq 1 \\
0 & {\rm otherwise}
\end{array}\right.
$$
which is exactly the cochain complex associated by Yetter
\cite{dY98} to a semigroupal functor. We have then the following
generalization of Yetter's result \cite{dY98}:

\begin{thm}
The equivalences classes of purely pseudofunctorial first order
deformations of a $K$-linear unitary pseudofunctor
$\F=(|\F|,\F_*,\widehat{\F}_*)$ are in bijection with the elements
of $H^2(\F)$.
\end{thm}
\begin{proof}
Let us consider 2-isomorphisms
$\widehat{\F}_{\epsilon}(g,f)=\widehat{\F}(g,f)+\widehat{\F}^{(1)}(g,f)\epsilon$
and $(\F_0)_{\epsilon}(X)=1_{id_{\F(X)}}+\F_0^{(1)}(X)\epsilon$,
with
$\widehat{\F}^{(1)}(g,f):\F(g)\circ\F(f)\Longrightarrow\F(g\circ
f)$ and $\F_0^{(1)}:\F(id_X)\Longrightarrow id_{\F(X)}$. We want to
find the necessary and sufficient conditions on the
$\widehat{\F}^{(1)}(g,f)$ and the $\F_0^{(1)}(X,Y)$ for these
2-isomorphisms to define a purely pseudofunctorial first order
deformation of $\F$. Let us first observe the following, which in
particular shows that such a first order deformation of a unitary
pseudofunctor $\F$ is completely determined by the 2-morphisms
$\widehat{\F}^{(1)}(g,f)$
\footnote{ This result should be viewed as an analog of Yetter's
result that a semigroupal deformation of a monoidal functor
uniquely extends to a monoidal deformation. See \cite{dY01},
Theorem 17.2.}:
\begin{lem}
Let $\F=(|\F|,\F_*,\widehat{\F}_*)$ be a $K$-linear unitary
pseudofunctor between $K$-linear 2-categories \Cg\ and \Dg, and
let's consider a purely pseudofunctorial first order deformation
given by 2-morphisms $\widehat{\F}^{(1)}(g,f)$ and $\F_0^{(1)}(X)$
(the deformation need not be unitary). Then, for all objects $X$ of
\Cg, we have:

(i) $\widehat{\F}(id_X,id_X)=1_{id_{\F(X)}}$.

(ii) $\F_0^{(1)}(X)=\widehat{\F}^{(1)}(id_X,id_X)$.
\end{lem}
\begin{proof}
(i) For any pseudofunctor between 2-categories, it directly follows
from the axioms that $\widehat{\F}(f,id_X)=1_{\F(f)}\circ \F_0(X)$,
for all 1-morphisms $f$. In particular, this is true when $f=id_X$.
Now, if $\F$ is unitary, we have $\F(id_X)=id_{\F(X)}$, and since
in any 2-category identity 2-morphisms of an identity 1-morphism
are units with respect to horizontal composition, we get
$\widehat{\F}(id_X,id_X)=\F_0(X)=1_{id_{\F(X)}}$.

(ii) The same argument as before shows that
$(\F_0)_{\epsilon}(X)=\widehat{\F}_{\epsilon}(id_X,id_X)$. Notice
that, although $\F_{\epsilon}$ is no longer unitary, it still holds
that $\F_{\epsilon}(id_X)=id_{\F_{\epsilon}(X)}$, which is the only
thing needed to show the previous equality. The desired result
follows then by taking the first order terms in $\epsilon$.
\end{proof}

Let us now prove the proposition. According to the lemma and the
definition of a pseudofunctor, the 2-morphisms
$\widehat{\F}^{(1)}(g,f)$ and $\F_0^{(1)}(X)$ above define a purely
pseudofunctorial first order deformation of $\F$ if and only if:
(1) $\F_0^{(1)}(X)=\widehat{\F}^{(1)}(id_X,id_X)$, and (2) the
$\widehat{\F}^{(1)}(g,f)$ are such that
$\widehat{\F}_{\epsilon}(g,f)$ is natural in $(g,f)$ and satisfies
the hexagonal and triangular axioms in
Definition~\ref{pseudofunctor}. Since $\widehat{\F}(g,f)$ is
natural in $g,f$ by hypothesis, naturality in $g,f$ of
$\widehat{\F}_{\epsilon}(g,f)$ amounts to the naturality of
$\widehat{\F}^{(1)}(g,f)$ in $g,f$. Hence the
$\widehat{\F}^{(1)}(g,f)$ define an element $\widehat{\F}^{(1)}\in
X^2(\F)$. On the other hand, the hexagonal axiom on $\F_{\epsilon}$
gives the following condition on $\widehat{\F}^{(1)}$: for all
composable morphisms
$X\stackrel{f}{\longrightarrow}Y\stackrel{g}{\longrightarrow}Z\stackrel{h}{\longrightarrow}T$
\begin{align*}
\widehat{\F}&(h,g\circ f)\cdot(1_{\F(h)}\circ\widehat{\F}^{(1)}(g,f))+
\widehat{\F}^{(1)}(h,g\circ f)\cdot(1_{\F(h)}\circ \widehat{\F}(g,f))= \\ &=
\widehat{\F}^{(1)}(h\circ g,f)\cdot(\widehat{\F}(h,g)\circ 1_{\F(f)})+
\widehat{\F}(h\circ g,f)\cdot(\widehat{\F}^{(1)}(h,g)\circ 1_{\F(f)})
\end{align*}
It is easily seen that this condition exactly corresponds to the
fact that $\delta(\widehat{\F}^{(1)})=0$. Hence,
$\widehat{\F}^{(1)}$ is a 2-cocycle of the complex
$X^{\bullet}(\F)$. As regards the triangular axioms, notice that
they imply no additional conditions on $\widehat{\F}^{(1)}$. For
example, since $\F_0^{(1)}(X)=\widehat{\F}^{(1)}(id_X,id_X)$, the
first of these triangular axioms gives the condition
$$
\widehat{\F}^{(1)}(f,id_X)=1_{\F(f)}\circ\widehat{\F}^{(1)}(id_X,id_X)
$$
for all 1-morphisms $f:X\To Y$. Now, this condition is nothing more
than the condition $\delta(\widehat{\F}^{(1)})(f,id_X,id_X)=0$, as
the reader may easily check.

Suppose now that the 2-morphisms $(\widehat{\F}^{(1)})'(g,f)$
define another purely pseudofunctorial first order deformation of
$\F$ equivalent to the previous one. We need to show that
$\widehat{\F}^{(1)}$ and $(\widehat{\F}^{(1)})'$ are cohomologous
2-cocycles. Now, from Proposition~\ref{equivalent_pseudo} and by
definition of $\delta$, it follows immediately that both
deformations are equivalent if and only if there exists
$\widehat{\xi}^{(1)}\in X^1(\F)$ such that
$$
(\widehat{\F}^{(1)})'-\widehat{\F}^{(1)}=\delta(\widehat{\xi}^{(1)}),
$$
as required. Let's remark that the third condition in
Proposition~\ref{equivalent_pseudo} is again superfluous. Indeed,
just take $f=g=id_X$ in the second condition and use the previous
lemma to conclude that
$$
(\F_0^{(1)})'(X)-\F_0^{(1)}(X)=(\widehat{\F}^{(1)})'(id_X,id_X)-
\widehat{\F}^{(1)}(id_X,id_X)=\widehat{\xi}^{(1)}(id_X)=\F_0(X)\cdot\widehat{\xi}^{(1)}(id_X)
$$
\end{proof}

%\subsection{}
%Finally, we can consider the question of extending a purely
%pseudofunctorial $n^{th}$-order deformation of $\F$ one higher
%order.
% As Yetter \cite{dY98},\cite{dY01}, we closely follow
%Gerstenhaber's argument in terms of pre-Lie systems \cite{mG64}.
%
%Recall that, according to Gerstenhaber \cite{mG63}, a {\sl (right)
%pre-Lie system} of $K$-vector spaces is any sequence
%$\ldots,V_{-1},V_0,V_1,\ldots$ of vector spaces over $K$ indexed by
%the integers together with a family of $K$-bilinear maps
%$\ast_i:V_m\times V_n\To V_{m+n}$ for all triples $m,n,i\geq 0$
%with $i\leq m$, satisfying the following conditions:
%
%So, suppose we have such a deformation, given by structural
%2-isomorphisms
%\begin{align*}
%\widehat{\F}_{\epsilon}(g,f)&\widehat{\F}(g,f)+\widehat{\F}^{(1)}(g,f)\epsilon+\cdots+
%\widehat{\F}^{(n)}(g,f)\epsilon^n
%\\ (\F_0)_{\epsilon}(X)&=1_{\F(id_X)}+\F_0^{(1)}(X)\epsilon+\cdots+\F_0^{(n)}(X)\epsilon^n
%\end{align*}
%The same argument as before shows that the 2-morphisms
%$\F_0^{(i)}(X)$, for all $i=1,\ldots,n$, are already implicit in
%the corresponding $\widehat{\F}^{(i)}(g,f)$, as the terms
%$\widehat{\F}^{(i)}(id_X,id_X)$. (...)
%
%\begin{thm}
%Let $\F=(|\F|,\F_*,\widehat{\F}_*)$ be a non-unital pseudofunctor.
%Then, the obstruction to extending a purely pseudofunctorial
%$n^{th}$-order deformation of $\F$ to a purely pseudofunctorial
%$(n+1)^{th}$-order deformation is an element of $H^3(\F)$.
%\end{thm}
%\begin{proof}
%???
%\end{proof}

\section{Cohomology theory for the deformations of the pentagonator}

\subsection{}
In this section we initiate the study of the infinitesimal
deformations of a $K$-linear semigroupal 2-category
$(\text\Cg,\otimes,a,\pi)$. Notice first of all that, according to
Definition~\ref{deformacio_2_categoria_semigrupal}, in a generic
infinitesimal deformation of \Cg\ all structural 2-isomorphisms
$\widehat{\otimes}((f',g'),(f,g))$, $\otimes_0(X,Y)$,
$\widehat{a}(f,g,h)$, $\pi_{X,Y,Z,T}$ will be deformed. Now,
instead of treating directly such a generic deformation from the
outset, we will proceed in three steps. So, in this section we
consider those infinitesimal deformations where only the
pentagonator is deformed, the tensor product and the associator
remaining undeformed. They will be called infinitesimal {\sl
pentagonator-deformations}. In the following section we will treat
the case where both the tensor product and the associator are
simultaneously deformed, although under the assumption that the
tensor product remains unitary, even after the deformation. These
deformations will be called infinitesimal {\sl unitary
(tensorator,associator)-deformations}. We will obtain in this way
two different cohomologies that separately describe the
deformations of both parts of the semigroupal structure. Section 9
is devoted to see how both cohomologies fit together in a global
cohomology describing the generic infinitesimal unitary
deformations.

\subsection{}
Let us consider an arbitrary $K$-linear semigroupal 2-category
$(\text\Cg,\otimes,a,\pi)$. Recall from Section 4 that, given the
data $(\text\Cg,\otimes,a)$, a pentagonator $\pi$ is defined as a
modification between two induced pseudonatural transformations
$a^{(4)},^{(4)}a:\otimes^{(4)}\Longrightarrow ^{(4)}\otimes$ which
satisfies the $K_5$ coherence relation.

More generally, given the pseudofunctor
$\otimes:\text\Cg\times\text\Cg\To\text\Cg$, we can consider the
induced pseudofunctors
$\otimes^{(n)},^{(n)}\otimes:\text\Cg\times\stackrel{n)}{\cdots}\times\text\Cg\To\text\Cg$,
$n\geq 3$, defined by
\begin{eqnarray*}
&\otimes^{(n)}=\otimes\circ(id_{\text\Cg}\times\otimes)\circ(id_{\text\Cg}\times
id_{\text\Cg}\times\otimes)\circ\cdots\circ(id_{\text\Cg}\times\stackrel{n-2)}{\cdots}\times
id_{\text\Cg}\times\otimes) \\
&^{(n)}\otimes=\otimes\circ(\otimes\times
id_{\text\Cg})\circ(\otimes\times id_{\text\Cg}\times
id_{\text\Cg})\circ\cdots\circ(\otimes\times
id_{\text\Cg}\times\stackrel{n-2)}{\cdots}\times id_{\text\Cg})
\end{eqnarray*}
These are just two examples of a lot of induced tensor products of
multiplicity $n$. In the same way, we can generalize the induced
pseudonatural transformations $a^{(4)},^{(4)}a$ to suitable
pseudonatural transformations
$a^{(n)},^{(n)}a:\otimes^{(n)}\Longrightarrow
^{(n)}\otimes$, for all $n\geq 4$. Here, we also have many possible
choices, because there are many possibles $a${\it -paths} (i.e.,
paths constructed as compositions of expansions of instances of the
1-isomorphisms $a_{X,Y,Z}$) from the completely right-parenthesized
object $\otimes^{(n)}(X_1,\ldots,X_n)$ to the completely
left-parenthesized one $^{(n)}\otimes(X_1,\ldots,X_n)$. In the case
$n=4$, the 1-isomorphisms of $a^{(4)}$ and $^{(4)}a$ are defined by
taking the {\it extremal paths}, i.e., those characterized by the
fact that, in each step, always the most internal parenthesis or
the most external parenthesis, respectively, is moved. This leads
us to introduce the following generalization.

\begin{defn}
Given 1-isomorphisms $a_{(X,Y,Z)}:X\otimes(Y\otimes Z)\To (X\otimes
Y)\otimes Z$ for all objects $(X,Y,Z)$ of $\text\Cg^3$, let
$a^{(n)}_{(X_1,\ldots,X_n)}$, $^{(n)}a_{(X_1,\ldots,X_n)}$, $n\geq
4$, be the $a$-paths from $\otimes^{(n)}(X_1,\ldots,X_n)$ to
$^{(n)}\otimes(X_1,\ldots,X_n)$ induced by the $a_{X,Y,Z}$ and
corresponding to always moving the most internal parenthesis and
the most external parenthesis, respectively.
\end{defn}

It is possible to give a more explicit description of these
1-isomorphisms. Indeed, the objects $\otimes^{(n)}(X_1,\ldots,X_n)$
and $^{(n)}\otimes(X_1,\ldots,X_n)$ can be graphically represented
as in Fig.\ref{objectes}. Then, the $a$-path $a^{(n)}_{X_1,\dots,X_n}$
corresponds to moving to the left all the legs associated to the
objects $X_2,\ldots,X_{n-1}$, starting with $X_{n-1}$ and so on
until $X_2$, while the path $^{(n)}a_{X_1,\dots,X_n}$ corresponds
to doing the same thing but starting with $X_2$ and so on until
$X_{n-1}$. Using this graphical presentation, we obtain the
following description of both paths:

\begin{figure}
\centering
\input{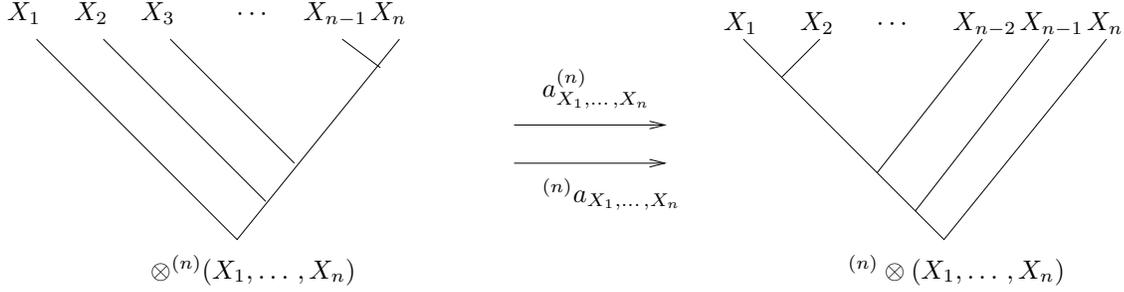}
\caption{A graphical representation of the objects
$\otimes^{(n)}(X_1,\ldots,X_n)$ and $^{(n)}\otimes(X_1,\ldots,X_n)$.}
\label{objectes}
\end{figure}

\begin{lem} \label{na_an}
For any $n\geq 4$ and any objects $X_1,\ldots,X_n$, we have
\begin{align*}
a^{(n)}_{X_1,\dots,X_n}=&\left(\prod_{i=3}^n(\cdots((a_{X_1,^{(i-2)}\otimes
(X_2,\ldots,X_{i-1}),X_i}\otimes id_{X_{i+1}})\otimes
id_{X_{i+2}})\cdots)\otimes id_{X_n}\right)\circ
\\ &\ \ \ \ \ \ \ \ \ \ \ \ \ \ \circ(id_{X_1}\hat{\otimes}
a^{(n-1)}_{X_2,\ldots,X_n})
\\ ^{(n)}a_{X_1,\dots,X_n}&=\prod_{i=2}^{n-1}a_{^{(n-i)}\otimes(X_1,\ldots,X_{n-i}),X_{n-i+1},
\otimes^{(i-1)}(X_{n-i+2},\ldots,X_n)}
\end{align*}
(the product denotes composition of 1-morphisms and the symbol
$\hat{\otimes}$ in the term $id_{X_1}\hat{\otimes}
a^{(n-1)}_{X_2,\ldots,X_n}$ is intended to mean the composition of
the tensor multiplications of $id_{X_1}$ by each one of the
composition factors defining $a^{(n-1)}_{X_2,\ldots,X_n}$).
\end{lem}

For example, when $n=5$, the reader may easily check that one
recovers the $a$-paths that appear in Fig.~\ref{relacioK5} defining
the common boundary of the polytope.
%namely
%\begin{align*}
%a^{(5)}_{X,Y,Z,T,U}&=((a_{X,Y,Z}\otimes id_T)\otimes id_U)\circ
%(a_{X,Y\otimes Z,T}\otimes id_U)\circ a_{X,(Y\otimes Z)\otimes
%T,U}\circ
%\\ &\ \ \ \ \ \ \ \circ(id_X\otimes(a_{Y,Z,T}\otimes id_U))\circ
%(id_X\otimes a_{Y,Z\otimes T,U})\circ
%(id_X\otimes(id_Y\otimes a_{Z,T,U}))
%\\ ^{(5)}a_{X,Y,Z,T,U}&=a_{(X\otimes Y)\otimes Z,T,U}\circ a_{X\otimes Y,Z,T\otimes U}
%\circ a_{X,Y,Z\otimes (T\otimes U)}
%\end{align*}

\begin{figure}
\centering
\input{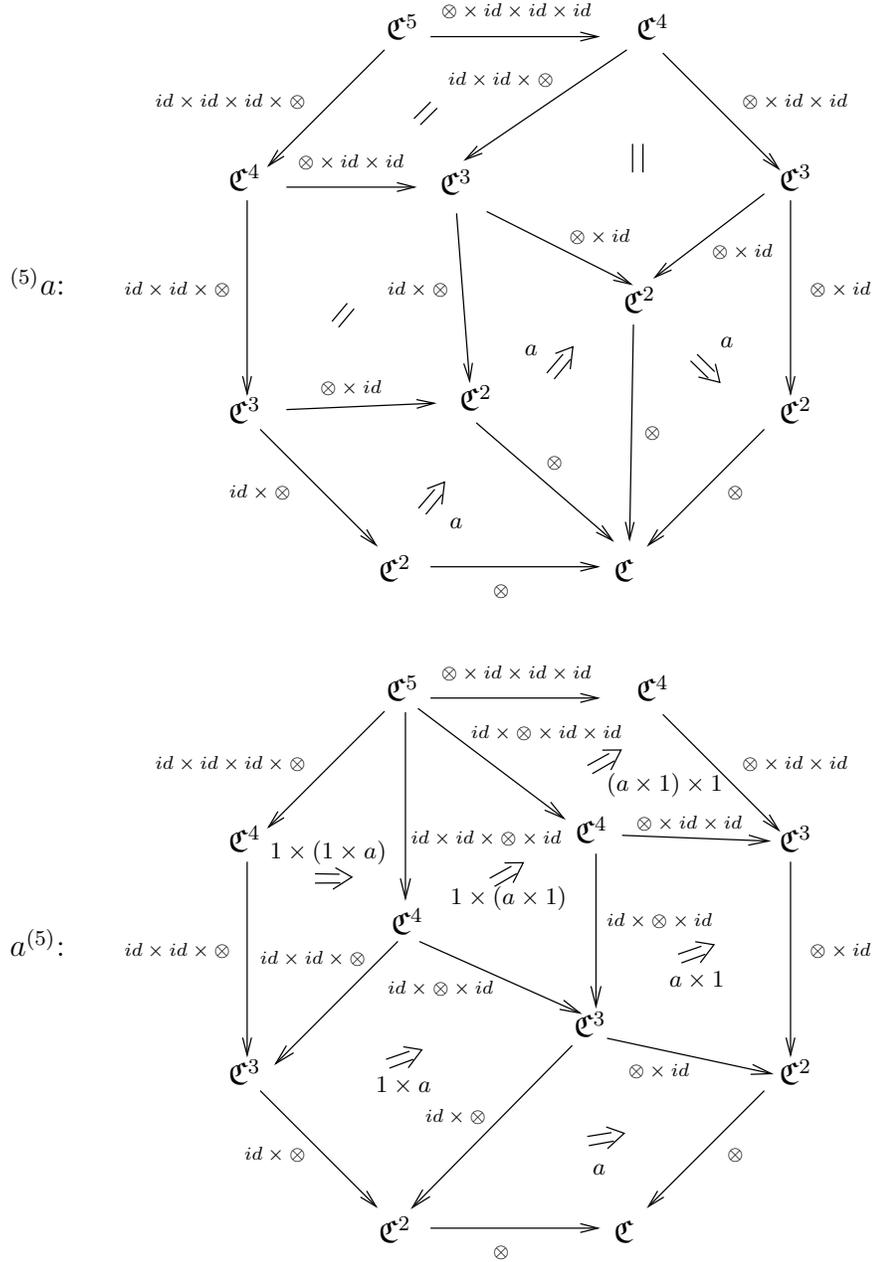}
\caption{Pastings defining the pseudonatural isomorphisms $^{(5)}a$ and $a^{(5)}$}
\label{figura_pastings_a_5}
\end{figure}

By definition, $a^{(n)},^{(n)}a:\otimes^{(n)}\Longrightarrow
^{(n)}\otimes$ are the pseudonatural isomorphisms (induced by $a$) whose 1-isomorphisms
are precisely the above 1-morphisms
$a^{(n)}_{X_1,\ldots,X_n},^{(n)}a_{X_1,\ldots,X_n}$. So, from the
formulas defining the vertical composition of 2-morphisms and the
horizontal compositions of the form $\xi\circ 1_{\F}$ (see Section
2), it is clear that $^{(n)}a$ is the pseudonatural isomorphism
given by the pasting
$$
^{(n)}a=\prod_{i=2}^{n-1}(a\circ 1_{^{(n-i)}\otimes\times id_{\Cg}\times\otimes^{(i-1)}})
$$
the product here denoting vertical composition of pseudonatural
transformations. The formula giving the pasting that defines
$a^{(n)}$ is a bit more complicated and is omitted because it is
not relevant in what follows.
% where
%it will appear the term $1_{\otimes}\circ(1_{id_{\Cg}}\times
%a^{(n-1)})$, corresponding to the tensor product $X_1\otimes
%a^{(n-1)}_{(X_2,\ldots,X_n)}$, vertically composed with similar
%terms. We leave to the reader to find the explicit formula.
In Fig.~\ref{figura_pastings_a_5}, however, both pastings are
explicitly represented in the case $n=5$. This defines the
pseudonatural isomorphisms $a^{(n)},^{(n)}a$ for all $n\geq 4$.
%\begin{lem} \label{na_an}
%For any $n\geq 4$ and any objects $X_1,\ldots,X_n$, we have
%\begin{eqnarray*}
%&a^{(n)}_{X_1,\dots,X_n}=\left(\prod_{i=3}^n(\cdots((a_{X_1,^{(i-2)}\otimes(X_2,\ldots,X_{i-1}),X_i}\otimes
%id_{X_{i+1}})\otimes id_{X_{i+2}})\cdots)\otimes id_{X_n}\right)(id_{X_1}\otimes
%a^{(n-1)}_{X_2,\ldots,X_n}) \\
%&^{(n)}a_{X_1,\dots,X_n}=\prod_{i=2}^{n-1}a_{^{(n-i)}\otimes(X_1,\ldots,X_{n-i}),X_{n-i+1},
%\otimes^{(i-1)}(X_{n-i+2},\ldots,X_n)}
%\end{eqnarray*}
%where the product denotes composition of 1-morphisms.
%\end{lem}
%
%These 1-isomorphisms correspond to induced pseudonatural
%isomorphisms $a^{(n)},^{(n)}a:\otimes^{(n)}\Longrightarrow
%^{(n)}\otimes$. For example, from the formulas defining
%the vertical composition of 2-morphisms and the horizontal
%compositions of the form $\xi\circ 1_{\F}$ (see Section 2), it is
%clear that $^{(n)}a$ is the pseudonatural isomorphism given by the pasting
%$$
%^{(n)}a=\prod_{i=2}^{n-1}(a\circ 1_{^{(n-i)}\otimes\times id_{\text\Cg}\times\otimes^{(i-1)}})
%$$
%the product here denoting vertical composition of pseudonatural
%transformations. We leave to the reader to find the pasting which
%gives $a^{(n)}$.
% where
%it will appear the term $1_{\otimes}\circ(1_{id_{\text\Cg}}\times
%a^{(n-1)})$, corresponding to the tensor product $X_1\otimes
%a^{(n-1)}_{(X_2,\ldots,X_n)}$, vertically composed with similar
%terms. We leave to the reader to find the explicit formula.
%In Fig.~\ref{figura_pastings_a_5}, both pastings are explicitly represented in the
%case $n=5$.
When $n=1,2,3$, let us take $\otimes^{(2)}=^{(2)}\otimes=\otimes$,
$\otimes^{(1)}=^{(1)}\otimes=id_{\text\Cg}$ and define
\begin{eqnarray*}
&a^{(3)}=^{(3)}a=a \\ &a^{(2)}=^{(2)}a=1_{\otimes}
\\ &a^{(1)}=^{(1)}a=1_{id_{\text\Cg}}
\end{eqnarray*}

\subsection{}
\begin{figure}
\centering
\input{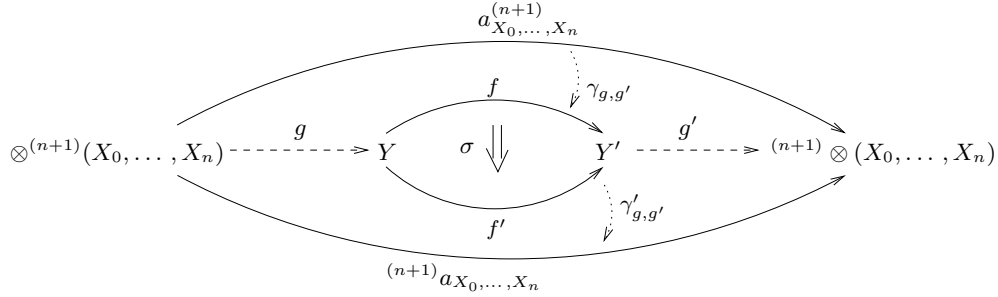}
\caption{Extending by canonical 2-isomorphisms.}
\label{figura_extensio}
\end{figure}

We can now define the cochain complex we are looking for. So, for
all $n\in{\mathbf N}$, let's denote by
$\widetilde{X}^n_{pent}(\text\Cg)$ the following vector spaces over
$K$:
$$
\widetilde{X}^n_{pent}(\text\Cg)=\left\{ \begin{array}{ll}
\text{PseudMod}(a^{(n+1)},^{(n+1)}a) & \text{if}\ n\geq 0 \\ 0 &
\text{otherwise}
\end{array} \right.
$$
where $\text{PseudMod}(a^{(n+1)},^{(n+1)}a)$ denotes the set of
pseudomodifications from $a^{(n+1)}$ to $^{(n+1)}a$ (see
Definition~\ref{modification}). They are indeed vector spaces over
$K$ because of Proposition~\ref{prop_K_linear}. As for the
coboundary operator
$\delta_{pent}:\widetilde{X}^{n-1}_{pent}(\text\Cg)\To\widetilde{X}^n_{pent}(\text\Cg)$,
we would like to take the usual formula, i.e.,
$$
(\delta_{pent}(\text\nn))_{X_0,\ldots,X_n}\approx
1_{id_{X_0}}\otimes\text\nn_{X_1,\ldots,X_n}+
\sum_{i=1}^n(-1)^i\text\nn_{X_0,\ldots,X_{i-1}\otimes X_i,\ldots,X_n}+
(-1)^{n+1}\text\nn_{X_0,\ldots,X_{n-1}}\otimes 1_{id_{X_n}}
$$
But the 2-morphisms corresponding to each of the terms in this sum
are not 2-cells from $a^{(n+1)}_{X_0,\ldots,X_n}$ to
$^{(n+1)}a_{X_0,\ldots,X_n}$, as required. If
$\sigma:f\Longrightarrow f'$ denotes any one of these 2-morphisms,
the situation is like that in Fig.~\ref{figura_extensio}.
Note, however, that $f,f'$ are always $a$-paths, because
\nn\ is a modification from $a^{(n)}$ to $^{(n)}a$.
The claim is that, once more, there exist suitable analogs of
Crane-Yetter's ``padding'' operators in this 2-dimensional setting
of pastings that give sense to the previous definition. Behind
these new ``padding'' operators there is again a coherence theorem,
which in this case can be stated as follows:

\begin{thm}
Let $(\text\Cg,\otimes,a,\pi)$ be a semigroupal 2-category, and let
$U,V$ be any two objects of \Cg\, both obtained as a certain tensor
product of the objects $X_1,\ldots,X_n$. Then, given any two
$a$-paths from $U$ to $V$, there is a unique 2-isomorphism between
them constructed as a pasting of instances of the structural
2-isomorphisms of \Cg\ and identity 2-morphisms (of expansions of
instances) of the structural 1-isomorphisms.
\end{thm}
\begin{proof}
Although this is a particular consequence of the strictification
theorem for semigroupal 2-categories (see
Theorem~\ref{strictification}), let's give a direct and somewhat
more appealing argument using the Stasheff polyhedra \cite{jS63}.
Recall that in \cite{jS63}, the author introduces, for each $n\geq
2$, a polyhedron $K_n$ whose vertices are in bijection with all
possible parenthesizations of a word $x_1x_2\cdots x_n$ of length
$n$ and whose edges all correspond to moves of the type $-(--)\To
(--)-$, where $-$ stands for a letter or a block of letters.
Stasheff shows that $K_n$ is homeomorphic to the
$(n-2)$-dimensional ball $D^{n-2}$. In particular, $K_5$ is a
homeomorphic image of the 3-ball whose faces (six pentagons and
three quadrilaterals) are those represented in
Fig.~\ref{relacioK5}. On the other hand, notice that the
$(n-3)$-dimensional faces of $K_n$ constituting its boundary
$\partial K_n$ correspond to all meaningful ways of inserting one
pair of parentheses $x_1x_2\cdots(x_k\cdots x_{k+s-1})\cdots x_n$,
where $2\leq s\leq n-1$ and $1\leq k\leq n-s+1$ (in particular,
$K_n$ has $n(n-1)/2-1$ such faces). Since the next insertion of
parentheses must be either within the block $(x_k\cdots x_{k+s-1})$
or treating this block as a unit, this face can be thought of as a
homeomorphic image of $K_{n-s+1}\times K_s$. All these faces are
not disjoint, but intersect along their boundaries in such a way
that the ``edges'' so formed correspond to inserting two pairs of
parentheses in the word $x_1\cdots x_n$. This allows one to
construct the $K_n$, for all $n\geq 2$, by induction: $K_2$ is a
point, and given $K_2,\ldots,K_{n-1}$, the next one $K_n$ is
defined as the cone on $\partial K_n$, where $\partial K_n$ is a
quotient of the form
$$
\partial K_n:=\left(\coprod_{s,k}(K_{n+s-1}\times K_s)_k\right)/\sim
$$
with $2\leq s\leq n-1$ and $1\leq k\leq n-s+1$.

MacLane's classical coherence theorem for semigroupal categories
(see, for example, \cite{sM98}) is nothing more that an algebraic
interpretation of the fact that the 2-dimensional skeleton of
$K_n$, for all $n\geq 4$, is a union of homeomorphic copies of
$K_4$ or $K_3\times K_3$ (see \cite{SS93}), a copy of $K_4$
corresponding to an instance of the Stasheff pentagon axiom and a
copy of $K_3\times K_3$ corresponding to an instance of a
naturality square of $a_{X,Y,Z}$ applied to a morphism which is
itself some $a_{X',Y',Z'}$. In the same way, the above coherence
result we want to prove is a consequence of the following fact
about these polyhedra:

\begin{lem}
For all $n\geq 5$, the 3-dimensional skeleton of $K_n$ is a union
of homeomorphic copies of $K_5$, $K_3\times K_4$ or $K_3\times
K_3\times K_3$.
\end{lem}
\begin{proof}
Indeed, the 3-cells of $K_n$ correspond to all ways of inserting
$n-5$ pairs of parenthesis in the word $x_1\cdots x_n$. Now, as we
have seen before, the insertion of the first pair gives an
$(n-3)$-cell of $K_n$ homeomorphic to a suitable product
$K_{s'_1}\times K_{s'_2}$. Similarly, the insertion of the second
pair corresponds to an $(n-4)$-cell of $K_n$ homeomorphic to some
product $K_{s''_1}\times K_{s''_2}\times K_{s''_3}$, because it is
obtained by substituting one of the previous factors $K_{s'_i}$ for
one of its faces, etc. We conclude that all 3-cells of $K_n$ will
be homeomorphic images of a suitable product
$K_{s_1}\times\cdots\times K_{s_{n-4}}$. Furthermore, since
$K_{s_i}$ is of dimension $s_i-2$, it must be
$(s_1-2)+\cdots+(s_{n-4}-2)=3$. This, together with the fact that
$s_i\geq 2$ for all $i=1,\ldots,n-3$, implies that at most three of
the $s_i$ can be greater than 2. In other words, any 3-cell of
$K_n$ is homeomorphic to a product of $n-7$ copies of $K_2$ by
$K_3\times K_3\times K_3$, by $K_2\times K_3\times K_4$ or by
$K_2\times K_2\times K_5$.
\end{proof}
To prove the proposition using this lemma, let's consider the
$(\otimes,a,\pi)$-{\it realization} of $K_n$ associated to the
objects $X_1,\ldots,X_n$, defined as follows: (1) as vertices, it
has all possible tensor products of $X_1,\ldots,X_n$, with all
possible parenthesizations, (2) as edges, it has expansions of
instances of the structural 1-isomorphisms $a_{X,Y,Z}$, and (3) as
2-faces, instances of the structural 2-isomorphisms
$\widehat{\otimes}((f',g'),(f,g))$, $\otimes_0(X,Y)$,
$\widehat{a}(f,g,h)$ and $\pi_{X,Y,Z,T}$. Observe that, here, the
$X,Y,Z,T$ are all objects obtained as tensor products of the
$X_1,\ldots,X_n$, and that the $f,g,h,f',g'$ are all identity
1-morphisms or instances of the 1-isomorphisms $a_{X,Y,Z}$. Any
a-path from $U$ to $V$ is then a path in this realization of $K_n$,
and the 2-isomorphisms mentioned in the proposition between two
such paths correspond to 2-faces between them in this realization.
Now, two such 2-faces are equal whenever the 3-cell diagram they
define is commutative. But, according to the previous lemma, any
3-cell in $K_n$ is a union of 3-cells of the types $K_5$,
$K_3\times K_4$ and $K_3\times K_3\times K_3$. The proof of the
proposition then finishes by checking that the 3-dimensional
diagrams corresponding to these three possible types of 3-cells are
just realizations of $K_5$, pentagonal prisms corresponding to the
naturality of the pentagonator in any one of the variables and
instances of the cube in Fig.~\ref{cub}, all of them commutative by
hypothesis.
\end{proof}

This unique 2-isomorphism will be called the {\it canonical}
2-isomorphism between both $a$-paths. Using them, we can extend any
$\sigma:f\Longrightarrow f'$ above to a 2-morphism from
$a^{(n+1)}_{X_0,\ldots,X_n}$ to $^{(n+1)}a_{X_0,\ldots,X_n}$ as
follows. Since the source and target objects $Y$, $Y'$ of $f$ and
$f'$ are canonically isomorphic to the reference objects
$\otimes^{(n+1)}(X_0,\ldots,X_n)$ and
$^{(n+1)}\otimes(X_0,\ldots,X_n)$, we can choose $a$-paths
$g:\otimes^{(n+1)}(X_0,\ldots,X_n)\To Y$ and $g':Y'\To
^{(n+1)}\otimes(X_0,\ldots,X_n)$, represented in Fig.~\ref{figura_extensio} by dashed
arrows. Now, by the previous coherence theorem, there are unique
canonical 2-isomorphisms
$\gamma_{g,g'}:a^{(n+1)}_{X_0,\ldots,X_n}\Longrightarrow g'\circ
f\circ g$ and $\gamma'_{g,g'}:g'\circ f\circ g\Longrightarrow
^{(n+1)}a_{X_0,\ldots,X_n}$. The desired extension of $\sigma$ is
then the pasting
$$
\lceil\lceil\sigma\rceil\rceil_{g,g'}:=
\gamma'_{g,g'}\cdot(1_{g'}\circ\sigma\circ 1_g)\cdot\gamma_{g,g'}.
$$
Since the 2-morphism $\sigma$ will not generally be a pasting of
the structural 2-isomorphisms of \Cg, this extension may a priori
depend on the paths $g,g'$. The next result shows that this is not
the case.

\begin{prop}
In the above notations, the extension
$\lceil\lceil\sigma\rceil\rceil_{g,g'}$ is independent of the
chosen canonical 1-isomorphisms $g,g'$. In particular, there is a
unique extension of $\sigma$ by canonical 2-isomorphisms.
\end{prop}
\begin{proof}
Let $\hat{g},\hat{g}'$ be any other choice. Then, we have the two extensions
\begin{align*}
\lceil\lceil\sigma\rceil\rceil_{g,g'}&=
\gamma'_{g,g'}\cdot(1_{g'}\circ\sigma\circ 1_g)\cdot\gamma_{g,g'}
\\ \lceil\lceil\sigma\rceil\rceil_{\hat{g},\hat{g}'}&=
\gamma'_{\hat{g},\hat{g}'}\cdot(1_{\hat{g}'}\circ\sigma\circ 1_{\hat{g}})
\cdot\gamma_{\hat{g},\hat{g}'}
\end{align*}
Now, since $g,\hat{g}$ are both canonical 1-isomorphisms between
the same vertices, coherence theorem implies that there exists a
unique canonical 2-isomorphism $\tau:g\Longrightarrow \hat{g}$. By
the same reason, we also have a unique canonical 2-isomorphism
$\tau':g'\Longrightarrow \hat{g}'$. Hence, the pastings
\begin{eqnarray*}
&(\tau'\circ 1_f\circ\tau)\cdot\gamma_{g,g'} \\
&\gamma'_{g,g'}\cdot((\tau')^{-1}\circ 1_{f'}\circ\tau^{-1})
\end{eqnarray*}
define canonical 2-isomorphisms from $a^{(n+1)}_{X_0,\ldots,X_n}$
to $\hat{g}'\circ f\circ\hat{g}$ and from $\hat{g}'\circ
f'\circ\hat{g}$ to $^{(n+1)}a_{X_0,\ldots,X_n}$, respectively. By
unicity, we must have
\begin{align*}
\gamma_{\hat{g},\hat{g}'}&=(\tau'\circ 1_f\circ\tau)\cdot\gamma_{g,g'} \\
\gamma'_{\hat{g},\hat{g}'}&=\gamma'_{g,g'}\cdot((\tau')^{-1}\circ 1_{f'}\circ\tau^{-1})
\end{align*}
Hence, applying the interchange law, we obtain that
\begin{align*}
\lceil\lceil\sigma\rceil\rceil_{\hat{g},\hat{g}'}&=
\gamma'_{g,g'}\cdot((\tau')^{-1}\circ 1_{f'}\circ\tau^{-1})
\cdot(1_{\hat{g}'}\circ\sigma\circ 1_{\hat{g}})\cdot(\tau'\circ 1_f\circ\tau)
\cdot\gamma_{g,g'} \\
&=\gamma'_{g,g'}\cdot(1_{g'}\circ\sigma\circ 1_g)\cdot\gamma_{g,g'} \\
&=\lceil\lceil\sigma\rceil\rceil_{g,g'}.
\end{align*}
\end{proof}
%
%Let us denote by $\lceil\lceil\sigma\rceil\rceil$ this unique
%extension of $\sigma$ by canonical 2-isomorphisms.

\begin{cor}
Let us consider a 2-morphism $\sigma:f\Longrightarrow f'$, where
$f,f'$ are some $a$-paths between suitable parenthesizations of the
tensor product of $X_0,\ldots,X_n$. Then, there exists a unique
extension of $\sigma$ by canonical 2-isomorphisms to a 2-morphism
between the reference $a$-paths $a^{(n+1)}_{X_0,\ldots,X_n}$ and
$^{(n+1)}a_{X_0,\ldots,X_n}$.
\end{cor}
Let us denote by $\lceil\lceil\sigma\rceil\rceil$ this unique
extension of $\sigma$ by canonical 2-isomorphisms. The
$\lceil\lceil-\rceil\rceil$ are, then, the analogs of the
``padding'' operators in this 2-dimensional setting (for the chosen
reference 1-morphisms). Notice that they should be strictly denoted
by $\lceil\lceil-\rceil\rceil_{X_0,\ldots,X_n}$, because there is
such an operator for every ordered set of objects
$(X_0,\ldots,X_n)$.

We can now define the coboundary operator
$\delta_{pent}:\widetilde{X}^{n-1}_{pent}(\Cg)\To\widetilde{X}^n_{pent}(\Cg)$
by
\begin{align*}
(\delta_{pent}(\text\nn))_{X_0,\ldots,X_n}=\lceil\lceil
&1_{id_{X_0}}\otimes\text\nn_{X_1,\ldots,X_n}\rceil\rceil
\\ &+\sum_{i=1}^n(-1)^{i}\lceil\lceil\text\nn_{X_0,\ldots,X_{i-1}\otimes
X_i,\ldots,X_n}\rceil\rceil+(-1)^{n+1}\lceil\lceil\text\nn_{X_0,\ldots,X_{n-1}}\otimes
1_{id_{X_n}}\rceil\rceil.
\end{align*}
Using similar arguments to those made to prove
Proposition~\ref{complex}, it can be shown that

\begin{prop}
For any $K$-linear semigroupal 2-category
$(\text\Cg,\otimes,a,\pi)$, the pair
$(\widetilde{X}^{\bullet}_{pent}(\text\Cg),\delta_{pent})$ is a
cochain complex.
\end{prop}

This complex will be called the {\sl general
pentagonator-deformation complex of} $(\text\Cg,\otimes,a,\pi)$,
and the corresponding cohomology groups will be denoted by
$\widetilde{H}^{\bullet}_{pent}(\text\Cg)$, the semigroupal
structure $(\otimes,a,\pi)$ being omitted for the sake of
simplicity. Note that the dependence on the pentagonator $\pi$
comes exclusively through the ``padding'' operators
$\lceil\lceil-\rceil\rceil$. Although this complex and its
cohomology will be relevant in the sequel, it is not the right
complex describing the infinitesimal pentagonator-deformations.
Indeed, we need to take the following subcomplex:

\begin{prop}
The vector subspaces ${\rm
Mod}(a^{(n)},^{(n)}a)\subset\widetilde{X}^{n-1}_{pent}(\text\Cg)$
define a subcomplex of the general pentagonator-deformation complex
of \Cg.
\end{prop}
\begin{proof}
We only need to see that the naturality of the
$\text\nn_{X_1,\ldots,X_n}$ in $(X_1,\ldots,X_n)$ implies that of
the $(\delta_{pent}(\text\nn))_{X_0,\ldots,X_n}$ in
$(X_0,\ldots,X_n)$. Let us consider for example the first term
$\lceil\lceil
1_{id_{X_0}}\otimes\text\nn_{X_1,\ldots,X_n}\rceil\rceil$. The
naturality in $(X_1,\ldots,X_n)$ of $\text\nn_{X_1,\ldots,X_n}$
implies that of $1_{id_{X_0}}\otimes\text\nn_{X_1,\ldots,X_n}$ in
$(X_0,\ldots,X_n)$. So, the situation is like that in
Fig.~\ref{naturalitat}, with two cylinders, one inside the other.
\begin{figure}
\centering
\input{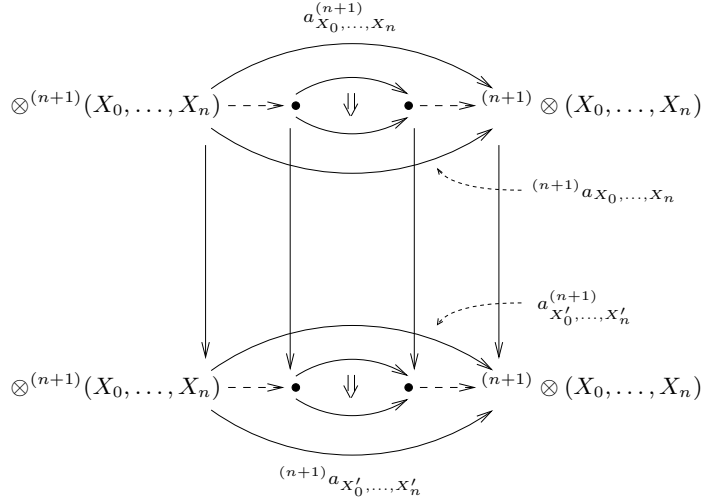}
\caption{Naturality of the first term in $\delta_{pent}(\text\nn)$.}
\label{naturalitat}
\end{figure}
We already know that the
inner one commutes, and we want to see that the same is true for
the outer one. Let's think of this outer cylinder without the inner
one as being decomposed in its upper and lower halves. Each one of
these halves is itself a cylinder. Now, both bases of any one of
these cylinders will correspond to {\it modifications} between
$a^{(n+1)}$ or $^{(n+1)}a$ and some other induced pseudonatural
isomorphism. Indeed, they are nothing more than the canonical
2-isomorphisms of the previous coherence theorem between the
corresponding $a$-paths. But these 2-isomorphisms are pastings of
2-isomorphisms all natural in $(X_0,\ldots,X_n)$. It then follows
that both halves also commute.
\end{proof}

This subcomplex will be denoted by $X^{\bullet}_{pent}(\text\Cg)$
and called the {\sl pentagonator-deformation complex of} \Cg. If we
denote its cohomology by $H_{pent}^{\bullet}(\text\Cg)$, we have
the following:

\begin{thm}
For any $K$-linear semigroupal 2-category
$(\text\Cg,\otimes,a,\pi)$, the $\omega$-equivalence classes of its
first order pentagonator-deformations are in bijection with the
elements of $H^3_{pent}(\text\Cg)$.
\end{thm}
\begin{proof}
Let's consider 2-isomorphisms of the form
\begin{align*}
\widehat{\otimes}_{\epsilon}((f',g'),(f,g))&=\widehat{\otimes}((f',g'),(f,g))
\\ (\otimes_0)_{\epsilon}(X,Y)&=(\otimes_0)(X,Y)
\\ \widehat{a}_{\epsilon}(f,g,h)&=\widehat{a}(f,g,h)
\\ (\pi_{\epsilon})_{X,Y,Z,T}&=\pi_{X,Y,Z,T}+(\pi^{(1)})_{X,Y,Z,T}\epsilon
\end{align*}
From Proposition~\ref{semigroupal_2_category}, it is easy to check
that they define a semigroupal structure on $\text\Cg^0_{(1)}$
(hence, a first order pentagonator-deformation of $\text\Cg$) if
and only if the following conditions are satisfied (each condition
is identified by the structural equation it comes from):

\begin{description}

\item[A$\pi$1]
The $\pi^{(1)}_{X,Y,Z,T}$ are natural in $(X,Y,Z,T)$, i.e., they
define an element $\pi^{(1)}\in X^3_{pent}(\text\Cg)$.

\item[A$\pi$2]
$\delta_{pent}(\pi^{(1)})=0$.
%\begin{align*}
%\lceil\lceil(\pi^{(1)})_{X\otimes Y,Z,T,U}\rceil\rceil+
%\lceil\lceil(\pi^{(1)})_{X,Y,Z\otimes T,U}\rceil\rceil+&
%\lceil\lceil(\pi^{(1)})_{X,Y,Z,T}\otimes 1_{id_U}\rceil\rceil= \\
%&=\lceil\lceil(\pi^{(1)})_{X,Y,Z,T\otimes U}\rceil\rceil+
%\lceil\lceil 1_{id_X}\otimes(\pi^{(1)})_{Y,Z,T,U}\rceil\rceil
%\end{align*}
\end{description}
The remaining structural equations are clearly superfluous in this
case. This proves that first order pentagonator-deformations of
\Cg\ indeed correspond to 3-cocylces of $X^{\bullet}_{pen}(\text\Cg)$.

Let's suppose now that two 3-cocycles $\pi^{(1)}$ and
$(\pi^{(1)})'$ define $\omega$-equivalent first order
pentagonator-deformations. We need to see that they are
cohomologous cocycles. Indeed, from Proposition~\ref{equivalents},
it easily follows that both deformations are $\omega$-equivalent if
and only if there exists 2-morphisms
$(\omega^{(1)})_{X,Y,Z}:a_{X,Y,Z}\Longrightarrow a_{X,Y,Z}$, hence,
and element $\omega^{(1)}\in\widetilde{X}^2_{pent}(\text\Cg)$, such
that (each condition is again identified with the correspondiong
condition in Proposition~\ref{equivalents} it comes from)
\begin{description}

\item[E$\omega$1]
The $(\omega^{(1)})_{X,Y,Z}$ are natural in $(X,Y,Z)$, i.e.,
$\omega^{(1)}\in X^2_{pent}(\text\Cg)$.

\item[E$\omega$2]
$(\pi^{(1)})'-\pi^{(1)}=-\delta_{pent}(\omega^{(1)})$
\end{description}
(the remaining conditions $(E\widehat{\psi}1)-(E\widehat{\psi}3)$
in Proposition~\ref{equivalents} are clearly empty in this case).
Hence, both 3-cocycles are indeed cohomologous.
\end{proof}

%(...)
%\begin{thm}
%Let $(\text\Cg,\otimes,\gamma,a,\pi)$ be a $K$-linear semigroupal
%2-category. Then, the obstruction to extending an $n^{th}$-order
%pentagonator-deformation of \Cg\ to an $(n+1)^{th}$-order
%pentagonator-deformation is an element of $H^?_{pent}(\text\Cg)$.
%\end{thm}
%\begin{proof}
%?????
%\end{proof}

\section{Cohomology theory for the unitary deformations of the tensor product and the associator}

\subsection{}
As already indicated in the previous section, in this section we
give a cohomological description of the infinitesimal
unitary\footnote{Recall that the term {\it unitary} applied to an
infinitesimal deformation always means that the deformed tensor
product is supposed to be still unitary, i.e., such that the
2-isomorphisms $\otimes_0(X,Y)$ remain trivial even after the
deformation. At the moment of writing, the author doesn't know how
to take into account the non trivial deformations of these
2-isomorphisms.} (tensorator,associator)-deformations. To do that,
we will make the simplifying assumption that the undeformed
semigroupal 2-category is actually a Gray semigroup, since
otherwise the theory becomes extremely cumbersome. This means,
however, no loss of generality because of
Theorem~\ref{strictification}

The situation we will encounter for these deformations closely
resembles the cohomology theory discovered by Gerstenhaber and
Schack\cite{GS90} to describe the infinitesimal deformations of a
bialgebra, and later extended by Crane and Yetter \cite{CY981} to
the case of a bitensor category (i.e., the categorification of a
bialgebra). So, we associate a double complex to any $K$-linear
Gray semigroup and prove that the second cohomology group of the
corresponding total complex provides us with the desired
description of the simultaneous first order unitary deformations of
both the tensor product and the associator. As we will see, the
role played by the multiplication and comultiplication in the
bialgebra case corresponds in our case to the tensor product and
composition of 1-morphisms. Furthermore, from this double complex
we will easily get cohomologies describing the (unitary)
deformations of the tensor product and the associator separately.
Roughly, they are respectively related to the rows and the columns
of the double complex of $(\text\Cg,\otimes)$, in much the same way
as in the classical bialgebra case.

\subsection{}
Let $(\text\Cg,\otimes)$ be a $K$-linear Gray semigroup. In
particular, $\otimes:\text\Cg\times\text\Cg\To\text\Cg$ is a
cubical pseudofunctor. Recall, however, that not all cubical
pseudofunctors $\otimes$ will provide the 2-category \Cg\ with the
structure of a Gray semigroup. More explicitly, the non-trivial
2-isomorphisms $\widehat{\otimes}((f',g'),(f,g))$ must additionally
satisfy the equation
\begin{align*}
(\widehat{\otimes}((f',g'),(f,g))\otimes
1_{h'\circ h})\cdot\widehat{\otimes}((f'\otimes g',h')&,(f\otimes
g,h))= \\ &=(1_{f'\circ
f}\otimes\widehat{\otimes}((g',h'),(g,h)))\cdot\widehat{\otimes}((f',g'\otimes
h'),(f,g\otimes h))
\end{align*}
coming from the structural condition $(A\widehat{a}2)$ in
Proposition~\ref{semigroupal_2_category} when the associator is
trivial (the reader may easily check that the remaining structural
equations $(A\widehat{a}3)$ and $(A\pi 1)-(A\pi 2)$ give no
additional conditions on $\otimes$).

Recall from the previous section that, for all $n\geq 1$, we
introduced pseudofunctors
$\otimes^{(n)},^{(n)}\otimes:\text\Cg\times\stackrel{n)}{\cdots}\times\text\Cg\To\text\Cg$.
Then, because of the above additional equation on $\otimes$, in a
Gray semigroup we have the following:

\begin{prop} \label{o_n}
Let $(\text\Cg,\otimes)$ be a Gray semigroup. Then,
for all $n\geq 1$, we have the equality of pseudofunctors
$$
^{(n)}\otimes=\otimes^{(n)}:=\otimes(n).
$$
Moreover, $\otimes(n)$ is unitary
% and the following equalities hold for all
%$(f'_0,\ldots,f'_n),(f_0,\ldots,f_n)$:
%\begin{align*}
%\widehat{\otimes(n+1)}((f'_0,\ldots,f'_n),(f_0,\ldots,f_n))&=
%(\widehat{\otimes}((f'_0,f'_1),(f_0,f_1))\otimes
%1_{(f'_2\circ f_2)\otimes\cdots\otimes(f'_n\circ f_n)})\cdot \\ &\
%\ \ \ \ \ \cdot\widehat{\otimes(n)}
%((f'_0\otimes f'_1,\ldots,f'_n),(f_0\otimes f_1,\ldots,f_n)) \\
%&=(1_{(f'_0\circ f_0)\otimes\cdots\otimes(f'_{n-2}\circ f_{n-2})}\otimes
%\widehat{\otimes}((f'_{n-1},f'_n),(f_{n-1},f_n)))\cdot \\
%&\ \ \ \ \ \ \cdot\widehat{\otimes(n)}
%((f'_0,\ldots,f'_{n-1}\otimes f'_n),(f_0,\ldots,f_{n-1}\otimes f_n)) \\
%&=(1_{f'_0\circ f_0}\otimes\widehat{\otimes(n)}
%((f'_1,\ldots,f'_n),(f_1,\ldots,f_n)))\cdot \\ &\ \ \ \ \ \ \cdot
%\widehat{\otimes}((f'_0,f'_1\otimes\cdots\otimes f'_n),
%(f_0,f_1\otimes\cdots\otimes f_n)) \\
%&=(\widehat{\otimes(n)}
%((f'_0,\ldots,f'_{n-1}),(f_0,\ldots,f_{n-1}))\otimes 1_{f'_n\circ f_n})\cdot
%\\ &\ \ \ \ \ \ \cdot\widehat{\otimes}
%((f'_0\otimes\cdots\otimes f'_{n-1},f'_n),
%(f_0\otimes\cdots\otimes f_{n-1},f_n))
%\end{align*}
\end{prop}

\begin{proof}
For $n=1,2$ it is obvious. Let's consider the case $n\geq 3$. Since
the structural isomorphisms $a_{X,Y,Z}$, $\widehat{a}(f,g,h)$ and
$\otimes_0(X,Y)$ are identities, it is clear that we only need to
prove that
$$
\widehat{^{(n)}\otimes}((f'_1,\ldots,f'_n),(f_1,\ldots,f_n))=
\widehat{\otimes^{(n)}}((f'_1,\ldots,f'_n),(f_1,\ldots,f_n))
$$
The proof is by induction on $n$. The case $n=3$ is nothing more
than the previously mentioned additional equation, as the reader
may easily check. Let $n>3$. By definition of $^{(n)}\otimes$ and
using the induction hypothesis, we have
\begin{align*}
\widehat{^{(n)}\otimes}((f'_1,\ldots,f'_n),(f_1,\ldots,f_n))&=
(\widehat{\otimes}((f'_1,f'_2),(f_1,f_2))\otimes 1_{(f'_3\circ
f_3)\otimes\cdots\otimes(f'_n\circ f_n)})\cdot \\ &\ \ \
\ \ \ \cdot\widehat{^{(n-1)}\otimes} ((f'_1\otimes
f'_2,\ldots,f'_n),(f_1\otimes f_2,\ldots,f_n)) \\
&=(\widehat{\otimes}((f'_1,f'_2),(f_1,f_2))\otimes 1_{(f'_3\circ
f_3)\otimes\cdots\otimes(f'_n\circ f_n)})\cdot \\ &\ \ \ \ \ \
\cdot\widehat{\otimes^{(n-1)}}
((f'_1\otimes f'_2,\ldots,f'_n),(f_1\otimes f_2,\ldots,f_n))
\end{align*}
Now, from the definition of $\otimes^{(n-1)}$ and using the equality
$(\tau'\cdot\tau)\otimes(\sigma'\cdot\sigma)=
(\tau'\otimes\sigma')\cdot(\tau\otimes\sigma)$, it follows that
\begin{align*}
\widehat{\otimes^{(n-1)}}
((f'_1\otimes f'_2,\ldots,f'_n),(f_1\otimes f_2,\ldots,f_n))&=
(1_{(f'_1\otimes f'_2)\circ (f_1\otimes f_2)}\otimes
\widehat{\otimes^{(n-2)}}((f'_3,\ldots,f'_n),(f_3,\ldots,f_n)))\cdot \\
&\ \ \ \cdot\widehat{\otimes}
((f'_1\otimes f'_2,f'_3\otimes\cdots\otimes f'_n),
(f_1\otimes f_2,f_3\otimes\cdots\otimes f_n))
\end{align*}
Therefore, we have
\begin{align*}
\widehat{^{(n)}\otimes}((f'_1\otimes f'_2,\ldots,f'_n),(f_1\otimes f_2,\ldots,f_n))&=
(\widehat{\otimes}((f'_1,f'_2),(f_1,f_2))\otimes 1_{(f'_3\circ
f_3)\otimes\cdots\otimes(f'_n\circ f_n)})\cdot \\ &\ \ \
\cdot(1_{(f'_1\otimes f'_2)\circ (f_1\otimes f_2)}\otimes
\widehat{\otimes^{(n-2)}}((f'_3,\ldots,f'_n),(f_3,\ldots,f_n)))\cdot \\
&\ \ \ \cdot\widehat{\otimes}
((f'_1\otimes f'_2,f'_3\otimes\cdots\otimes f'_n),
(f_1\otimes f_2,f_3\otimes\cdots\otimes f_n)) \\
&=(1_{(f'_1\circ f_1)\otimes(f'_2\circ f_2)}\otimes
\widehat{\otimes^{(n-2)}}((f'_3,\ldots,f'_n),(f_3,\ldots,f_n)))\cdot \\ &\ \ \
\cdot(\widehat{\otimes}((f'_1,f'_2),(f_1,f_2))\otimes
1_{(f'_3\otimes\cdots\otimes f'_n)\circ(f_3\otimes\cdots\otimes f_n)})
\cdot \\ &\
\ \ \cdot\widehat{\otimes}
((f'_1\otimes f'_2,f'_3\otimes\cdots\otimes f'_n),
(f_1\otimes f_2,f_3\otimes\cdots\otimes f_n))
\end{align*}
The proof finishes by applying $(A\widehat{a}2)$ to the last two factors.
\end{proof}

\subsection{}
Since all the pseudofunctors $\otimes(n)$ are unitary, we have for
each of them the corresponding cochain complex
$X^{\bullet}(\otimes(n))$, $n\geq 1$, describing their purely
pseudofunctorial deformations (see Section 6). More precisely, if
$m\leq 0$, $X^m(\otimes(n))=0$, while for all $m\geq 1$, it is
$$
X^m(\otimes(n))=\prod_{(X^0_1,\ldots,X^0_n),\ldots,(X^m_1,\ldots,X^m_n)\in|\text\Cg^n|}
\text{Nat}(C\otimes(n)_{(X^0_i),\ldots,(X^m_i)},\otimes(n)C_{(X^0_i),\ldots,(X^m_i)})
$$
To simplify, we write $(X^j_i)$ instead of
$(X^j_1,\ldots,X^j_n)$. Here,
$C\otimes(n)_{(X^0_i),\ldots,(X^m_i)}$ and
$\otimes(n)C_{(X^0_i),\ldots,(X^m_i)}$ denote the functors
\begin{align*}
\text\Cg^n((X^{m-1}_1,\ldots,X^{m-1}_n),(X^m_1,\ldots,X^m_n))\times
\cdots\times\text\Cg^n&((X^0_1,\ldots,X^0_n),(X^1_1,\ldots,X^1_n))\To \\
&\To \text\Cg(X^0_1\otimes\cdots\otimes
X^0_n,X^m_1\otimes\cdots\otimes X^m_n)
\end{align*}
which apply the composable 1-morphisms
$$
(f^1_i,\ldots,f^n_i):(X^{m-1-i}_1,\ldots,X^{m-1-i}_n)\To(X^{m-i}_1,\ldots,X^{m-i}_n),\
\ \ \ \ \ i=0,\ldots,m-1
$$
to $(f^1_0\otimes\cdots\otimes
f^m_0)\circ\cdots\circ(f_{m-1}^1\otimes\cdots\otimes f_{m-1}^n)$
and $(f_0^1\circ\cdots\circ
f_{m-1}^1)\otimes\cdots\otimes(f_0^n\circ\cdots\circ f_{m-1}^n)$,
respectively. Hence, a generic element $\phi\in X^m(\otimes(n))$ is
a collection of 2-morphisms
\begin{align*}
\phi((f_0^1,\ldots,f_0^n),\ldots,(f_{m-1}^1,\ldots,f_{m-1}^n)):
(f^1_0\otimes\cdots&\otimes
f^n_0)\circ\cdots\circ(f_{m-1}^1\otimes\cdots\otimes
f_{m-1}^n)\Longrightarrow \\ &\Longrightarrow(f_0^1\circ\cdots\circ
f_{m-1}^1)\otimes\cdots\otimes(f_0^n\circ\cdots\circ f_{m-1}^n)
\end{align*}
natural in the $(f_i^1,\ldots,f_i^n)$. In particular, notice that
the structural 2-isomorphisms $\widehat{\otimes}((f',g'),(f,g))$
define an element $\widehat{\otimes}\in X^2(\otimes(2))$, while the
$\widehat{a}(f,g,h)$ define an element $\widehat{a}\in
X^1(\otimes(3))$ (the trivial one in the case of a Gray semigroup).

Instead of $X^m(\otimes(n))$, let's use the more suggestive
notation $X^{m-1,n-1}(\text\Cg,\otimes)$, or just
$X^{m-1,n-1}(\text\Cg)$, for these $K$-vector spaces (the change of
indices is for later convenience). They can be arranged as in
Figure~\ref{figura_complex_doble}, with $m-1\geq 0$ and $n-1\geq 0$
being the row and column index, respectively. Since the elements
$\phi\in X^{m-1,0}(\text\Cg)$ are of the form
$\phi(f_0,\ldots,f_{m-1}):f_0\circ\cdots\circ
f_{m-1}\Longrightarrow f_0\circ\cdots\circ f_{m-1}$, while those
$\phi\in X^{0,n}(\text\Cg)$ are of the form
$\phi(f^0,\ldots,f^n):f^0\otimes\cdots\otimes f^n\Longrightarrow
f^0\otimes\cdots\otimes f^n$, we can think of the rows and columns
as related to the composition and the tensor product, respectively,
of 1-morphisms.

\begin{figure}
\centering
\input{double.pstex_t}
\caption{Arrangement of the vector spaces $X^m(\otimes(n))=X^{m-1,n-1}(\text\Cg)$}
\label{figura_complex_doble}
\end{figure}

Arranged in this way, each row corresponds to the cochain complexes
$X^{\bullet}(\otimes(n))$, the coboundary operators
$\delta_h:X^{m-1,n-1}(\text\Cg)\To X^{m,n-1}(\text\Cg)$, $m\geq 1$,
being those defined in the previous section. Namely, if $\phi\in
X^{m-1,n-1}(\text\Cg)$, then
\begin{align*}
(\delta_h(\phi))((f^1_0,\ldots,&f_0^n),\ldots,(f^1_m,\ldots,f_m^n))=
\lceil 1_{f_0^1\otimes\cdots\otimes f_0^n}\circ
\phi((f^1_1,\ldots,f_1^n),\ldots,(f^1_m,\ldots,f^n_m)\rceil
\\
&+\sum_{i=1}^{m}(-1)^i\lceil\phi((f^1_0,\ldots,f_0^n),\ldots,(f^1_{i-1}\circ
f^1_i,\ldots,f_{i-1}^n\circ
f_i^n),\ldots,(f^1_{m},\ldots,f^n_{m}))\rceil+
\\ &+(-1)^{m+1}\lceil\phi((f^1_0,\ldots,f_0^n),\ldots,(f^1_{m-1},\ldots,f_{m-1}^n))\circ
1_{f^1_{m}\otimes\cdots\otimes f^n_{m}}\rceil
\end{align*}
The claim is that it is possible to define vertical coboundary maps
$\delta_v:X^{m-1,n-1}(\text\Cg)\To X^{m-1,n}(\text\Cg)$, for all
$m\geq 1$, making each column a cochain complex and in such a way
that the whole set of vector spaces and maps define a double
complex. Indeed, if $\phi\in X^{m-1,n-1}(\text\Cg)$, let's define
\begin{align*}
(\delta_v(\phi))((f^0_0,\ldots,&f_0^n),\ldots,(f^0_{m-1},\ldots,f_{m-1}^n))=
\lceil 1_{f_0^0\circ\cdots\circ f_{m-1}^0}\otimes
\phi((f^1_0,\ldots,f_0^n),\ldots,(f^1_{m-1},\ldots,f^n_{m-1})\rceil
\\
&+\sum_{i=1}^n(-1)^i\lceil\phi((f^0_0,\ldots,f_0^{i-1}\otimes
f_0^i,\ldots,f_0^n),\ldots,(f^0_{m-1},\ldots,f_{m-1}^{i-1}\otimes
f_{m-1}^i,\ldots,f^n_{m-1}))\rceil+
\\ &+(-1)^{n+1}\lceil\phi((f^0_0,\ldots,f_0^{n-1}),\ldots,
(f^0_{m-1},\ldots,f_{m-1}^{n-1}))\otimes
1_{f^n_0\circ\cdots\circ f^n_{m-1}}\rceil
\end{align*}
Using arguments similar to those used to prove previous results of
the same kind, one shows the following (once more, the coherence
theorem for unitary pseudofunctors takes account of the padding
operators):

\begin{prop}
For all $m\geq 1$, the pair $(X^{m,\bullet}(\text\Cg),\delta_v)$ is
a cochain complex.
\end{prop}

Actually, as indicated before, we have the following stronger
result, which is fundamental in our theory:

\begin{thm}
The $K$-vector spaces $X^{\bullet,\bullet}(\text\Cg)$ together with
the above maps $\delta_h:X^{\bullet,\bullet}(\text\Cg)\To
X^{\bullet+1,\bullet}(\text\Cg)$ and
$\delta_v:X^{\bullet,\bullet}(\text\Cg)\To
X^{\bullet,\bullet+1}(\text\Cg)$ define a double complex.
\end{thm}
\begin{proof}
It remains to prove that both coboundary maps $\delta_h$ and
$\delta_v$ commute. Let's consider an element $\phi\in
X^{m-2,n-1}(\text\Cg)$, $m\geq 2,n\geq 1$, with
$\phi=\{\phi((f_1^1,\ldots,f_1^n),\ldots,(f_{m-1}^1,\ldots,f_{m-1}^n))\}$.
Then, the reader may easily check that

\begin{align*}
(\delta_h&(\delta_v(\phi))((f_0^0,\ldots,f_0^n),\ldots,(f_{m-1}^0,\ldots,f_{m-1}^n))=
\\ &=\lceil 1_{f_0^0\otimes\cdots\otimes f_0^n}\circ\lceil 1_{f_1^0\circ\cdots\circ f_{m-1}^0}
\otimes\phi((f_1^1,\ldots,f_1^n),\ldots,(f_{m-1}^1,\ldots,f_{m-1}^n))\rceil\rceil
\\ &\ \ +\sum_{i=1}^n(-1)^i\lceil 1_{f_0^0\otimes\cdots\otimes
f_0^n}\circ\lceil\phi((f_1^0,\ldots,f_1^{i-1}\otimes
f_1^i,\ldots,f_1^n),\ldots,(f_{m-1}^0,\ldots,f_{m-1}^{i-1}\otimes
f_{m-1}^i,\ldots,f_{m-1}^n))\rceil\rceil
\\ &\ \ +(-1)^{n+1}\lceil
1_{f_0^0\otimes\cdots\otimes
f_0^n}\circ\lceil\phi((f_1^0,\ldots,f_1^{n-1}),\ldots,(f_{m-1}^0,\ldots,f_{m-1}^{n-1}))
\otimes 1_{f_1^n\circ\cdots\circ f_{m-1}^n}\rceil\rceil
\\ &\ \ +\sum_{i=1}^{m-1}(-1)^i
\lceil 1_{f_0^0\circ\cdots\circ f_{m-1}^0}\otimes
\phi((f_0^1,\ldots,f_o^n),\ldots,(f_{i-1}^1\circ f_i^1,\ldots,f_{i-1}^n\circ f_i^n),\ldots,
(f_{m-1}^1,\ldots,f_{m-1}^n))\rceil
\\ &\ \ +\sum_{i=1}^{m-1}\sum_{j=1}^n(-1)^{i+j}\lceil
\phi((f_0^0,\ldots,f_0^{j-1}\otimes f_0^j,\ldots,f_0^n),\ldots,(f_{i-1}^0\circ f_i^0,\ldots
\\ &\ \ \ \ \ \ \ \ \ \ \ \ \
\ \ \ldots,(f_{i-1}^{j-1}\circ f_i^{j-1})
\otimes(f_{i-1}^j\circ f_i^j),
\ldots,f_{i-1}^n\circ f_i^n),\ldots,(f_{m-1}^0,\ldots,f_{m-1}^{j-1}
\otimes f_{m-1}^j,\ldots,f_{m-1}^n))\rceil
\\ &\ \ +\sum_{i=1}^{m-1}(-1)^{i+n+1}
\lceil\phi((f_0^0,\ldots,f_0^{n-1}),\ldots,(f_{i-1}^0\circ f_i^0,\ldots,
f_{i-1}^{n-1}\circ
f_i^{n-1}),\ldots,(f_{m-1}^0,\ldots,f_{m-1}^{n-1}))\otimes
1_{f_0^n\circ\cdots\circ f_{m-1}^n}\rceil
\\ &\ \ +(-1)^m\lceil\lceil 1_{f_0^0\circ\cdots\circ f_{m-2}^0}\otimes
\phi((f_0^1,\ldots,f_0^n),\ldots,(f_{m-2}^1,\ldots,f_{m-2}^n))\rceil
\circ 1_{f_{m-1}^0\otimes\cdots\otimes f_{m-1}^n}\rceil
\\ &\ \ +\sum_{i=1}^n(-1)^{m+i}
\lceil\lceil\phi((f_0^0,\ldots,f_0^{i-1}\otimes f_0^i,\ldots,f_o^n),\ldots,
(f_{m-2}^0,\ldots,f_{m-2}^{i-1}\otimes
f_{m-2}^i,\ldots,f_{m-2}^n))\rceil\circ
1_{f_{m-1}^0\otimes\cdots\otimes f_{m-1}^n}\rceil
\\ &\ \ +(-1)^{m+n+1}\lceil\lceil
\phi((f_0^0,\ldots,f_0^{n-1}),\ldots,(f_{m-2}^0,\ldots,f_{m-2}^{n-1}))
\otimes 1_{f_0^n\circ\cdots\circ f_{m-2}^n}\rceil
\circ 1_{f_{m-1}^0\otimes\cdots\otimes f_{m-1}^n}\rceil
\end{align*}
while

\begin{align*}
(\delta_v&(\delta_h(\phi))((f_0^0,\ldots,f_0^n),\ldots,(f_{m-1}^0,\ldots,f_{m-1}^n))=
\\ &=\lceil 1_{f_0^0\circ\cdots\circ f_{m-1}^0}\circ\lceil 1_{f_0^1\otimes\cdots
\otimes f_0^n}
\otimes\phi((f_1^1,\ldots,f_1^n),\ldots,(f_{m-1}^1,\ldots,f_{m-1}^n))\rceil\rceil
\\ &\ \ +\sum_{i=1}^{m-1}(-1)^i
\lceil 1_{f_0^0\circ\cdots\circ f_{m-1}^0}\otimes\lceil
\phi((f_0^1,\ldots,f_o^n),\ldots,(f_{i-1}^1\circ f_i^1,\ldots,f_{i-1}^n\circ f_i^n),\ldots,
(f_{m-1}^1,\ldots,f_{m-1}^n))\rceil\rceil
\\ &\ \ +(-1)^m\lceil 1_{f_0^0\circ\cdots\circ f_{m-1}^0}\otimes\lceil
\phi((f_0^1,\ldots,f_0^n),\ldots,(f_{m-2}^1,\ldots,f_{m-2}^n))
\circ 1_{f_{m-1}^1\otimes\cdots\otimes f_{m-1}^n}\rceil\rceil
\\ &\ \ +\sum_{i=1}^n(-1)^i\lceil 1_{f_0^0\otimes\cdots\otimes
f_0^n}\circ\phi(f_1^0,\ldots,f_1^{i-1}\otimes
f_1^i,\ldots,f_1^n),\ldots,(f_{m-1}^0,\ldots,f_{m-1}^{i-1}\otimes
f_{m-1}^i,\ldots,f_{m-1}^n))\rceil
\\ &\ \ +\sum_{i=1}^n\sum_{j=1}^{m-1}(-1)^{i+j}\lceil
\phi((f_0^0,\ldots,f_0^{i-1}\otimes f_0^i,\ldots,f_0^n),\ldots,(f_{j-1}^0\circ f_j^0,\ldots
\\ &\ \ \ \ \ \ \ \ \ \ \ \ \
\ \ \ldots,\ldots,(f_{j-1}^{i-1}\otimes f_{j-1}^i)
\circ(f_j^{i-1}\otimes f_j^i),
\ldots,f_{j-1}^n\circ f_j^n),\ldots,(f_{m-1}^0,\ldots,f_{m-1}^{i-1}
\otimes f_{m-1}^i,\ldots,f_{m-1}^n))\rceil
\\ &\ \ +\sum_{i=1}^n(-1)^{m+i}
\lceil\phi((f_0^0,\ldots,f_0^{i-1}\otimes f_0^i,\ldots,f_o^n),\ldots,
(f_{m-2}^0,\ldots,f_{m-2}^{i-1}\otimes
f_{m-2}^i,\ldots,f_{m-2}^n))\circ 1_{f_{m-1}^0\otimes\cdots\otimes
f_{m-1}^n}\rceil
\\ &\ \ +(-1)^{n+1}\lceil\lceil
1_{f_0^0\otimes\cdots\otimes
f_0^{n-1}}\circ\phi((f_1^0,\ldots,f_1^{n-1}),\ldots,(f_{m-1}^0,\ldots,f_{m-1}^{n-1}))\rceil
\otimes 1_{f_0^n\circ\cdots\circ f_{m-1}^n}\rceil
\\ &\ \ +\sum_{i=1}^{m-1}(-1)^{i+n+1}
\lceil\lceil\phi((f_0^0,\ldots,f_0^{n-1}),\ldots,(f_{i-1}^0\circ f_i^0,\ldots,
f_{i-1}^{n-1}\circ
f_i^{n-1}),\ldots,(f_{m-1}^0,\ldots,f_{m-1}^{n-1}))\rceil\otimes
1_{f_0^n\circ\cdots\circ f_{m-1}^n}\rceil
\\ &\ \ +(-1)^{m+n+1}\lceil\lceil
\phi((f_0^0,\ldots,f_0^{n-1}),\ldots,(f_{m-2}^0,\ldots,f_{m-2}^{n-1}))
\circ 1_{f_{m-1}^0\otimes\cdots\otimes f_{m-1}^{n-1}}\rceil
\otimes 1_{f_0^n\circ\cdots\circ f_{m-1}^n}\rceil
\end{align*}
Notice that in both expressions there are nine terms. Now, recall
that taking the padding $\lceil-\rceil$ of some 2-morphism $(-)$
simply means to take a vertical composition of $(-)$ with the
appropriate expansions of the 2-isomorphisms
$\widehat{\otimes}(-,-)$. It then follows by the interchange law
that
$$
\lceil 1_{f_0^0\otimes\cdots\otimes
f_0^n}\circ\lceil-\rceil\rceil=\lceil\lceil
1_{f_0^0\otimes\cdots\otimes f_0^n}\circ(-)\rceil\rceil=\lceil
1_{f_0^0\otimes\cdots\otimes f_0^n}\circ(-)\rceil
$$
This proves the equality between the second term in the expression
of $\delta_h(\delta_v(\phi))$ and the fourth term in the expression
of $\delta_v(\delta_h(\phi))$. The same argument shows the equality
between the eighth and sixth terms in the first and second
expression, respectively. On the other hand, we can also conclude
that the first term in $\delta_h(\delta_v(\phi))$ is the padding of
$$
1_{f_0^0\otimes\cdots\otimes f_0^n}
\circ(1_{f_1^0\circ\cdots\circ f_{m-1}^0}
\otimes\phi((f_1^1,\ldots,f_1^n),\ldots,(f_{m-1}^1,\ldots,f_{m-1}^n)))
$$
But by the naturality of the $\widehat{\otimes}(-,-)$ (equation
$(A\widehat{\otimes}1)$), this is the same as
$$
\widehat{\otimes}(-,-)^{-1}\cdot(1_{f_0^0\circ\cdots\circ f_{m-1}^0}\otimes
(1_{f_0^1\otimes\cdots\otimes
f_0^n}\circ\phi((f_1^1,\ldots,f_1^n),\ldots,(f_{m-1}^1,\ldots,f_{m-1}^n))))
\cdot\widehat{\otimes(-,-)}
$$
whose padding clearly coincides with that of
$$
1_{f_0^0\circ\cdots\circ f_{m-1}^0}\otimes
(1_{f_0^1\otimes\cdots\otimes
f_0^n}\circ\phi((f_1^1,\ldots,f_1^n),\ldots,(f_{m-1}^1,\ldots,f_{m-1}^n)))
$$
Using now that
$(\tau'\cdot\tau)\otimes(\sigma'\cdot\sigma)=(\tau'\otimes\sigma')\cdot(\tau\otimes\sigma)$,
we obtain that both first terms also coincide. Similar arguments
can be made to show the equality between both last terms and
between the terms: third and seventh, fourth and second, sixth and
eighth and seventh and third in the first and second expressions,
respectively. Hence, it only remains to prove the equality between
both fifth terms in each expression, and this easily follows from
the naturality of $\phi$ in its arguments applied to the 2-morphism
$$
1_{f_0^0\otimes\cdots\otimes
f_0^n}\circ\cdots\circ(1_{f_{i-1}^0\circ
f_i^0}\otimes\cdots\otimes\widehat{\otimes}((f_{j-1}^{i-1},f_{j-1}^i),(f_j^{i-1},f_j^i))
\otimes\cdots\otimes 1_{f_{j-1}^n\circ f_j^n})\circ\cdots\circ
1_{f_{m-1}^0\otimes\cdots\otimes f_{m-1}^n}
$$
\end{proof}

This double complex
$(X^{\bullet,\bullet}(\text\Cg),\delta_h,\delta_v)$ will be called
the {\sl extended double complex} of the Gray semigroup
$(\text\Cg,\otimes)$. We are actually interested in the double
complex obtained after deleting the bottom row $X^{m,0}(\text\Cg)$.
It will be called the {\sl double complex} of $(\text\Cg,\otimes)$.
Furthermore, for our purposes, we need to take a subcomplex of this
double complex. This is related to the fact that we only consider
infinitesimal {\it unitary} deformations.

\begin{defn}
An element $\phi\in X^{m-1,n-1}(\text\Cg)$ will be called {\sl
special} if whenever $(f_i^1,\ldots,f_i^n)=(id,\ldots,id)$ for some
$i\in\{0,\ldots,m-1\}$, it holds
$$
\phi((f_0^1,\ldots,f_0^n),\ldots,(f_{m-1}^1,\ldots,f_{m-1}^n))=0
$$
\end{defn}
The set of $\phi\in X^{m-1,n-1}(\text\Cg)$, $m,n\geq 1$, which are
special clearly define a vector subspace, which will be denoted by
$X_s^{m-1,n-1}(\text\Cg)$. We have then the following:

\begin{prop}
The vector subspaces $X_s^{m,n}(\text\Cg)$, $m,n\geq 0$, define a
subcomplex $X_s^{\bullet,\bullet}(\text\Cg)$ of the extended double
complex of $(\text\Cg,\otimes)$.
\end{prop}
\begin{proof}
We need to see that both coboundary operators $\delta_h$ and
$\delta_v$ preserve the special elements. Indeed, let $\phi\in
X_s^{m-1,n-1}(\text\Cg)$. Then, from the above expression of
$(\delta_h(\phi))$, it is clear that when
$((f_i^1,\ldots,f_i^n)=(id,\ldots,id)$ for some
$i\in\{0,\ldots,m-1\}$, all terms are zero except the $(i+1)^{th}$
and the $(i+2)^{th}$ terms, which are equal but of opposite sign
(recall that, $\otimes$ being unitary, the tensor product of
identity 1-morphisms is always an identity 1-morphism, and that the
identity 2-morphism of an identity 1-morphism is a unit with
respect to horizontal composition). Hence, $\delta_h(\phi)$ is
special. On the other hand, if $\phi$ is special and
$((f_i^1,\ldots,f_i^n)=(id,\ldots,id)$ for some
$i\in\{0,\ldots,m-1\}$, all terms in $\delta_v(\phi)$ are clearly
zero, so that $\delta_v(\phi)$ is also special.
\end{proof}

The double complex defined by the vector subspaces
$X_s^{m,n}(\text\Cg)$, $m,n\geq 0$, and the corresponding
restrictions of both $\delta_h$ and $\delta_v$ will be called the
{\sl special extended double complex} of $(\text\Cg,\otimes)$, or
just the {\sl special double complex} of $(\text\Cg,\otimes)$, when
the bottom row is deleted.

\subsection{}
Let $X^{\bullet}_{tens,ass}(\text\Cg)$ denote the total complex
associated to the special double complex
$(X_s^{\bullet,\bullet}(\text\Cg),\delta_h,\delta_v)$ of
$(\text\Cg,\otimes)$. It will be called the {\sl unitary
(tensorator,associator)-deformation complex} of
$(\text\Cg,\otimes)$. By definition, it is the cochain complex with
vector spaces
$$
X^q_{tens,ass}(\text\Cg)=\bigoplus_{\begin{array}{c}m+n=q \\ m\geq
0,n\geq 1 \end{array}} X_s^{m,n}(\text\Cg),\
\
\
\ \ q\geq 1
$$
and coboundary operators
$\delta_{tens,ass}:X^q_{tens,ass}(\text\Cg)\To
X^{q+1}_{tens,ass}(\text\Cg)$ given by
$$
\delta_{tens,ass}=\bigoplus_{\begin{array}{c} m+n=q \\ m\geq  0,n\geq 1 \end{array}}
((-1)^n\delta_h+\delta_v)\ \ \ \ \ \ q\geq 1
$$
The corresponding cohomology groups will be denoted by
$H^{\bullet}_{tens,ass}(\text\Cg)$. The reason we choose the above
name for this total complex is the following theorem:

\begin{thm} \label{tensorator,associator-deformations}
If $(\text\Cg,\otimes)$ is a $K$-linear Gray semigroup, the
$\psi$-equivalence classes of its first order unitary
(tensorator,associator)-deformations are in bijection with the
elements of the group $H^2_{tens,ass}(\text\Cg)$.
\end{thm}

\begin{proof}
Let's consider 2-isomorphisms
\begin{align*}
\widehat{\otimes}_{\epsilon}((f',g'),(f,g))&=\widehat{\otimes}((f',g'),(f,g))+
\widehat{\otimes}^{(1)}((f',g'),(f,g))\epsilon
\\ (\otimes_0)_{\epsilon}(X,Y)&=1_{id_{X\otimes Y}}%+\otimes_0^{(1)}(X,Y)\epsilon
\\ \widehat{a}_{\epsilon}(f,g,h)&=1_{f\otimes g\otimes h}+\widehat{a}^{(1)}(f,g,h)\epsilon
\\ (\pi_{\epsilon})_{X,Y,Z,T}&=1_{id_{X\otimes Y\otimes Z\otimes T}}
\end{align*}
Substituting these 2-isomorphisms in the structural equations in
Proposition~\ref{semigroupal_2_category} and computing the first
order term in $\epsilon$, it turns out that they define a first
order unitary (tensorator,associator)-deformation of \Cg\ if and
only if the 2-morphisms $\widehat{\otimes}^{(1)}((f',g'),(f,g))$,
$\widehat{a}^{(1)}(f,g,h)$ satisfy the following conditions:
\begin{description}
\item[A$\widehat{\otimes}$1]
The $\widehat{\otimes}^{(1)}((f',g'),(f,g))$ are natural in
$(f',g'),(f,g)$, i.e., they define an element
$$
\widehat{\otimes}^{(1)}\in X^{1,1}(\text\Cg)=X^2(\otimes(2))
$$
\item[A$\widehat{\otimes}$2]
For all composable 1-morphisms $(f'',g''),(f',g'),(f,g)$
\begin{align*}
\widehat{\otimes}&^{(1)}((f'',g''),(f'\circ f,g'\circ g))
\cdot(1_{f''\otimes g''}\circ\widehat{\otimes}((f',g'),(f,g)))+
\\ \widehat{\otimes}&((f'',g''),(f'\circ f,g'\circ g))
\cdot(1_{f''\otimes g''}\circ\widehat{\otimes}^{(1)}((f',g'),(f,g)))=
\\ &\ \ \ \ \ \ \ \ \ \ \ \ =\widehat{\otimes}^{(1)}((f''\circ f',g''\circ g'),(f,g))
\cdot(\widehat{\otimes}((f'',g''),(f',g'))\circ 1_{f\otimes g})+
\\ &\ \ \ \ \ \ \ \ \ \ \ \ \ \ \ \widehat{\otimes}((f''\circ f',g''\circ g'),(f,g))
\cdot(\widehat{\otimes}^{(1)}((f'',g''),(f',g'))\circ 1_{f\otimes g})
\end{align*}
It is easy to check that this is exactly the condition
$\delta_h(\widehat{\otimes}^{(1)})=0$.

\item[A$\widehat{\otimes}$3]
For all 1-morphisms $(f,g):(X,Y)\To (X',Y')$, it holds
$$
\widehat{\otimes}^{(1)}((id_{X'},id_{Y'}),(f,g))=
\widehat{\otimes}^{(1)}((f,g),(id_X,id_Y))=0
$$
i.e., $\widehat{\otimes}^{(1)}\in X^{1,1}_s(\text\Cg)\subset
X^2_{tens,ass}(\text\Cg)$.
% In particular, we obtain that
%$\otimes_0^{(1)}(X,Y)=\widehat{\otimes}^{(1)}((id_X,id_Y),(id_X,id_Y))$.
%
\item[A$\widehat{a}$1]
The $\widehat{a}^{(1)}(f,g,h)$ are natural in $(f,g,h)$, i.e., they
define an element
$$
\widehat{a}^{(1)}\in X^{0,2}(\text\Cg)=X^1(\otimes(3))
$$
\item[A$\widehat{a}$2]
For all 1-morphisms $(f,g,h)$
\begin{align*}
(\widehat{\otimes}&^{(1)}((f',g'),(f,g))\otimes 1_{h'\circ
h})\cdot\widehat{\otimes}((f'\otimes g',h'),(f\otimes g,h))+
\\ &(\widehat{\otimes}((f',g'),(f,g))\otimes 1_{h'\circ
h})\cdot\widehat{\otimes}^{(1)}((f'\otimes g',h'),(f\otimes g,h))
\\ &\widehat{a}^{(1)}(f'\circ f,g'\circ g,h'\circ
h)\cdot(\widehat{\otimes}((f',g'),(f,g))\otimes 1_{h'\circ
h})\cdot\widehat{\otimes}((f'\otimes g',h'),(f\otimes g,h))=
\\ &\ \ \ \ \ \ \ \ \ \ =(1_{f'\circ
f}\otimes\widehat{\otimes}^{(1)}((g',h')(g,h)))\cdot
\widehat{\otimes}((f',g'\otimes h'),(f,g\otimes h))+
\\ &\ \ \ \ \ \ \ \ \ \ \ \ \ \ (1_{f'\circ
f}\otimes\widehat{\otimes}((g',h')(g,h)))\cdot
\widehat{\otimes}^{(1)}((f',g'\otimes h'),(f,g\otimes h))+
\\ &\ \ \ \ \ \ \ \ \ \ \ \ \ \ (1_{f'\circ
f}\otimes\widehat{\otimes}((g',h')(g,h)))\cdot
\widehat{\otimes}((f',g'\otimes h'),(f,g\otimes h))\cdot
(\widehat{a}^{(1)}(f',g',h')\circ 1_{f\otimes(g\otimes h)})+
\\ &\ \ \ \ \ \ \ \ \ \ \ \ \ \ (1_{f'\circ
f}\otimes\widehat{\otimes}((g',h')(g,h)))\cdot
\widehat{\otimes}((f',g'\otimes h'),(f,g\otimes h))\cdot
(1_{(f'\otimes g')\otimes h'}\circ\widehat{a}^{(1)}(f,g,h))
\end{align*}
It is easy to check that this is exactly the condition
$\delta_v(\widehat{\otimes}^{(1)})+\delta_h(\widehat{a}^{(1)})=0$
\item[A$\widehat{a}$3]
For all objects $(X,Y,Z)$, it holds
$$
\widehat{a}^{(1)}(id_X,id_Y,id_Z)=0
$$
i.e., $\widehat{a}^{(1)}\in X^{0,2}_s(\text\Cg)\subset
X^2_{tens,ass}(\text\Cg)$.
\item[A$\pi$1]
For all 1-morphisms $(f,g,h,k)$
\begin{align*}
\widehat{a}&^{(1)}(f,g,h\otimes k)+\widehat{a}^{(1)}(f\otimes g,h,k)=
\\ &\ \ \ \ \ \ =-\widehat{\otimes}^{(1)}
((id_{X'},id_{Y'\otimes Z'\otimes T'}),(f,g\otimes h\otimes k))
+1_f\otimes \widehat{a}^{(1)}(g,h,k)
\\ &\ \ \ \ \ \ \ \ \ \ +\widehat{\otimes}^{(1)}
((f,g\otimes h\otimes k),(id_{X},id_{Y\otimes Z\otimes T}))
+\widehat{a}^{(1)}(f,g\otimes h,k)
\\ &\ \ \ \ \ \ \ \ \ \ -\widehat{\otimes}^{(1)}
((id_{X'\otimes Y'\otimes Z'},id_{T'}),(f\otimes g\otimes h,k))
+\widehat{a}^{(1)}(f,g,h)\otimes 1_k
\\ &\ \ \ \ \ \ \ \ \ \ +\widehat{\otimes}^{(1)}
((f\otimes g\otimes h,k),(id_{X\otimes Y\otimes Z},id_T))
\end{align*}
Since the terms in $\widehat{\otimes}^{(1)}$ are all zero by
condition $(A\widehat{\otimes}3)$, this exactly corresponds to the
condition $\delta_v(\widehat{a}^{(1)})=0$.
\end{description}
As the reader may easily check, the structural equation $(A\pi 2)$
gives no additional conditions for this kind of deformation. Now,
$(A\widehat{\otimes}1),(A\widehat{\otimes}3)$ together with
$(A\widehat{a}1),(A\widehat{a}3)$ say that
$(\widehat{a}^{(1)},\widehat{\otimes}^{(1)})\in
X^2_{tens,ass}(\text\Cg)$. On the other hand, we have
$$
\delta_{tens,ass}(\widehat{a}^{(1)},\widehat{\otimes}^{(1)})=(\delta_v(\widehat{a}^{(1)}),
\delta_h(\widehat{a}^{(1)})+\delta_v(\widehat{\otimes}^{(1)}),
-\delta_h(\widehat{\otimes}^{(1)}))
$$
Hence, $(A\widehat{\otimes}2)$, $(A\widehat{a}2)$ and $(A\pi 1)$
together say that $(\widehat{a}^{(1)},\widehat{\otimes}^{(1)})$ is
a 2-cocycle.

Let us now suppose that
$((\widehat{a}^{(1)})',(\widehat{\otimes}^{(1)})')$ is another
2-cocycle defining a $\psi$-equivalent first order unitary
(tensorator,associator)-deformation of \Cg. We need to see that
both 2-cocycles are actually cohomologous. Indeed, from
Proposition~\ref{equivalents}, it immediately follows that they are
$\psi$-equivalent deformations if and only if there exists
2-morphisms $\widehat{\psi}^{(1)}(f,g):f\otimes g\Longrightarrow
f\otimes g$ such that
\begin{description}

\item[E$\widehat{\psi}$1]
The $\widehat{\psi}^{(1)}(f,g)$ are natural in $(f,g)$, i.e., they
define an element $\widehat{\psi}^{(1)}\in X^{0,1}(\text\Cg)$

\item[E$\widehat{\psi}$2]
$(\widehat{\otimes}^{(1)})'-\widehat{\otimes}^{(1)}=-\delta_h(\widehat{\psi}^{(1)})$.

\item[E$\widehat{\psi}$3]
$\widehat{\psi}^{(1)}$ is special.

\item[E$\omega$1]
$(\widehat{a}^{(1)})'-\widehat{a}^{(1)}=\delta_v(\widehat{\psi}^{(1)})$.
\end{description}
(in this case, the condition coming from equation $(E\omega 2)$ is
empty). The first and third conditions together say that
$\widehat{\psi}^{(1)}\in
X^{0,1}_s(\text\Cg)=X^1_{tens,ass}(\text\Cg)$, while the second and
fourth express nothing more than the fact that
$$
((\widehat{a}^{(1)})',(\widehat{\otimes}^{(1)})')-
(\widehat{a}^{(1)},\widehat{\otimes}^{(1)})=\delta_{tens,ass}(\widehat{\psi}^{(1)})
$$
as required.
\end{proof}

\subsection{}
With the above results, it is easy to obtain a cochain complex
whose cohomology describes the infinitesimal {\sl
associator-deformations} of $(\text\Cg,\otimes)$, i.e., those
deformations where only the associator is deformed. Indeed, such a
deformation is given by 2-isomorphisms of the form
\begin{align*}
\widehat{\otimes}_{\epsilon}((f',g'),(f,g))&=\widehat{\otimes}((f',g'),(f,g))
\\ (\otimes_0)_{\epsilon}(X,Y)&=1_{id_{X\otimes Y}}
\\ \widehat{a}_{\epsilon}(f,g,h)&=1_{f\otimes g\otimes h}+\widehat{a}^{(1)}(f,g,h)\epsilon
\\ (\pi_{\epsilon})_{X,Y,Z,T}&=1_{id_{X\otimes Y\otimes Z\otimes T}}
\end{align*}
According to the proof of the previous theorem, they define a
first-order associator-deformation of \Cg\ if and only if the
$\widehat{a}^{(1)}(f,g,h)$ define an element $\widehat{a}^{(1)}\in
X_s^{0,2}(\text\Cg)$ which moreover satisfies that
$$
\delta_v(\widehat{a}^{(1)})=0=\delta_h(\widehat{a}^{(1)})
$$
The first equality just says that $\widehat{a}^{(1)}$ is a
2-cocycle of the cochain complex
$(X^{0,\bullet}_s(\text\Cg),\delta_v)$, while the second one serves
to define a subcomplex of this complex. More explicitly, let's
define
$$
X^n_{ass}(\text\Cg):=\text{Ker}(\delta_h:X^{0,n}_s(\text\Cg)\To
X^{1,n}_s(\text\Cg))\ \ \ \ \ \ n\geq 0
$$
Since $X^{\bullet,\bullet}_s(\text\Cg)$ is a double complex, its
horizontal coboundary map $\delta_h$ is a morphism of complexes.
But the kernel of a morphism of complexes is a subcomplex. Hence,
the subspaces $X^n_{ass}(\text\Cg)$, $n\geq 0$, together with the
corresponding restriction of $\delta_v$, which will be denoted by
$\delta_{ass}$, define a cochain complex. Let us call it the {\sl
associator-deformation complex} of the Gray semigroup
$(\text\Cg,\otimes)$. Then, if $H^{\bullet}_{ass}(\text\Cg)$ denote
the corresponding cohomology groups, the following result is an
immediate consequence of the previous theorem:

\begin{thm}
If $(\text\Cg,\otimes)$ is a $K$-linear Gray semigroup, the
$\psi$-equivalence classes of its first order
associator-deformations are in bijection with the elements of
$H^2_{ass}(\text\Cg)$.
\end{thm}
\begin{proof}
Notice that the $\widehat{\psi}^{(1)}\in X^{0,2}_s(\text\Cg)$ must
be such that
$\delta_h(\widehat{\psi}^{(1)})=-(\widehat{\otimes}^{(1)})'+\widehat{\otimes}^{(1)}=0$,
i.e., $\widehat{\psi}^{(1)}\in X^1_{ass}(\text\Cg)$, as required.
\end{proof}

\subsection{}
In a similar way, we can easily get a cochain complex describing
the infinitesimal {\sl unitary tensorator-deformations} of
$\text\Cg$, i.e., those deformations where only the tensor product
is deformed, and in such a way that it remains unitary.  Indeed, such a
deformation is given by 2-isomorphisms of the form
\begin{align*}
\widehat{\otimes}_{\epsilon}((f',g'),(f,g))&=\widehat{\otimes}((f',g'),(f,g))+
\widehat{\otimes}^{(1)}((f',g'),(f,g))\epsilon
\\ (\otimes_0)_{\epsilon}(X,Y)&=1_{id_{X\otimes Y}}
\\ \widehat{a}_{\epsilon}(f,g,h)&=1_{f\otimes g\otimes h}%+\widehat{a}^{(1)}(f,g,h)\epsilon
\\ (\pi_{\epsilon})_{X,Y,Z,T}&=1_{id_{X\otimes Y\otimes Z\otimes T}}
\end{align*}
Again going back to the proof of
Theorem~\ref{tensorator,associator-deformations}, it immediately
follows that they define a first-order unitary
tensorator-deformation of
\Cg\ if and only if the $\widehat{\otimes}^{(1)}((f',g'),(f,g))$ define an
element $\widehat{\otimes}^{(1)}\in X_s^{1,1}(\text\Cg)$ which
moreover satisfies that
$$
\delta_h(\widehat{\otimes}^{(1)})=0=\delta_v(\widehat{\otimes}^{(1)})
$$
Now, recall that $X_s^{n,1}(\text\Cg)=X_s^{n+1}(\otimes)$, the
special $n$-cochains of the purely pseudofunctorial deformation
complex of the (unitary) pseudofunctor $\otimes$ (see Section 6).
Then, the first equality just says that $\widehat{\otimes}^{(1)}$
is a 2-cocycle of this cochain complex $X^{\bullet}_s(\otimes)$,
while the second one serves to define a subcomplex of this complex.
More explicitly, let's define
$$
X^n_{tens}(\text\Cg):=\text{Ker}(\delta_v:X^{n+1}_s(\otimes)\To
X^{n,1}_s(\text\Cg))\ \ \ \ \ \ n\geq 0
$$
The same argument as before shows that these subspaces define a
subcomplex. Let us call it the {\sl unitary tensorator-deformation
complex} of the Gray semigroup $(\text\Cg,\otimes)$. Then, if
$H^{\bullet}_{tens}(\text\Cg)$ denote the corresponding cohomology
groups, the following result is an immediate consequence of the
previous theorem:

\begin{thm}
If $(\text\Cg,\otimes)$ is a $K$-linear Gray semigroup, the
$\psi$-equivalence classes of its first order unitary
tensorator-deformations are in bijection with the elements of
$H^2_{tens}(\text\Cg)$.
\end{thm}

\section{Cohomology theory for the generic unitary deformations}

\subsection{}
In Sections 7 and 8 we have constructed complexes
$X^{\bullet}_{pent}(\text\Cg)$ and
$X^{\bullet}_{tens,ass}(\text\Cg)$ whose cohomologies separately
describe the infinitesimal unitary deformations of the
pentagonator, on the one hand, and the tensor product and the
associator, on the other, of a $K$-linear Gray semigroup
$(\text\Cg,\otimes)$ (actually, of an arbitrary $K$-linear
semigroupal 2-category in the case of the deformations of the
pentagonator). The goal of this section is to obtain a cohomology
describing the simultaneous unitary deformations of both sets of
structural 2-isomorphisms. To do that, we will need to go back to
the bigger cochain complex
$\widetilde{X}^{\bullet}_{pent}(\text\Cg)\supseteq
X^{\bullet}_{pent}(\text\Cg)$ introduced in Section 7. The
appropriate cohomology for the generic unitary deformations turns
out to be the total complex of a modified version of the double
complex of $(\text\Cg,\otimes)$ introduced in the previous section.
The modification consists of the substitution of the first column
$X^{1,\bullet}(\text\Cg)$ in this double complex by a suitable cone
of that column and the general pentagonator-deformation complex
$\widetilde{X}^{\bullet}_{pent}(\text\Cg)$.

\subsection{}
Recall that given two cochain complexes $(A^{\bullet},\delta_A)$
and $(B^{\bullet},\delta_B)$ and a morphism of complexes
$\varphi:(A^{\bullet},\delta_A)\To (B^{\bullet},\delta_B)$, the
cone of $A^{\bullet}$ and $B^{\bullet}$ over $\varphi$ is the cochain complex
$(C^{\bullet}_{\varphi}(A^{\bullet},B^{\bullet}),\delta_C)$ defined
by
$$
C^n_{\varphi}(A^{\bullet},B^{\bullet})=B^n\oplus A^{n+1}
$$
and
$$
\delta_C(b,a)=(\delta_B(b)+\varphi(a),-\delta_A(a)),\ \ \ \ (b,a)\in B^n\oplus A^{n+1}.
$$
The minus sign in $\delta_A$ is to ensure that
$\delta_C\circ\delta_C=0$. For more details see, for example,
Weibel\cite{cW94}.

\subsection{}
Let $(\text\Cg,\otimes)$ be a $K$-linear Gray semigroup, and let's
consider the general pentagonator-deformation complex
$\widetilde{X}^{\bullet}_{pent}(\text\Cg)$ defined in Section 7.
The associator $a$ of \Cg\ being trivial, this complex reduces to
$$
\widetilde{X}^{n-1}_{pent}(\text\Cg)=\text{PseudMod}({\mathbf 1}_n,{\mathbf 1}_n),\ \ \ \ \
n\geq 1
$$
where ${\mathbf 1}_n$ is the pseudonatural isomorphism ${\mathbf
1}_n:\otimes(n)\Longrightarrow\otimes(n)$ whose structural 1- and
2-isomorphisms are all identities, while the coboundary is given by
$$
(\delta_{pent}(\text\nn))_{X_0,\ldots,X_n}=1_{id_{X_0}}\otimes\text\nn_{X_1,\ldots,X_n}+
\sum_{i=1}^n(-1)^{i}\text\nn_{X_0,\ldots,X_{i-1}
\otimes X_i,\ldots,X_n}+(-1)^{n+1}\text\nn_{X_0,\ldots,X_{n-1}}\otimes 1_{id_{X_n}}
$$
(in this case, the padding operators are superfluous). Notice that
a generic element
$\text\nn\in\widetilde{X}^{n-1}_{pent}(\text\Cg)$, $n\geq 1$, is
just a collection of 2-morphisms
$\text\nn_{X_1,\ldots,X_n}:id_{X_1\otimes\cdots\otimes
X_n}\Longrightarrow id_{X_1\otimes\cdots\otimes X_n}$, for all
objects $(X_1,\ldots,X_n)$ of $\text\Cg^n$.

Let us consider, on the other hand, the complex
$(X^{0,\bullet}(\text\Cg),\delta_v)$ corresponding to the first
column  of the extended double complex of $(\text\Cg,\otimes)$. We
can define $K$-linear maps
$\varphi:\widetilde{X}^{n-1}_{pent}(\text\Cg)\To
X^{0,n-1}(\text\Cg)$, $n\geq 1$ by
$$
(\varphi(\text\nn))(f_1,\ldots,f_n)=\text\nn_{X'_1,\ldots,X'_n}\circ
1_{f_1\otimes\cdots\otimes f_n}-1_{f_1\otimes\ldots\otimes
f_n}\circ\text\nn_{X_1,\ldots,X_n}
$$
for all 1-morphisms $(f_1,\ldots,f_n):(X_1,\ldots,X_n)\To
(X'_1,\ldots,X'_n)$ of $\text\Cg^n$.

\begin{prop}
The above maps
$\varphi_{\bullet}:\widetilde{X}^{\bullet}_{pent}(\text\Cg)\To
X^{0,\bullet}(\text\Cg)$ define a morphism of cochain complexes.
\end{prop}
\begin{proof}
Let $\text\nn\in\widetilde{X}^{n-1}_{pent}(\text\Cg)$. Then, we
have
\begin{align*}
(\delta_v(\varphi(\text\nn)))(f_0,\ldots,f_n)&=1_{f_0}\otimes
(\varphi(\text\nn))(f_1,\ldots,f_n)+
\\ &\ \ \ +\sum_{i=1}^n(-1)^i(\varphi(\text\nn))
(f_0,\ldots,f_{i-1}\otimes f_i,\ldots,f_n)+ \\ &\ \ \ +
(-1)^{n+1}(\varphi(\text\nn))(f_0,\ldots,f_{n-1})\otimes 1_{f_n}
\\ &=-1_{f_0}\otimes(1_{f_1\otimes\cdots\otimes f_n}
\circ\text\nn_{X_1,\ldots,X_n})+ \\ &\ \ \ +1_{f_0}\otimes(\text\nn_{X'_1,\ldots,X'_n}\circ
1_{f_1\otimes\cdots\otimes f_n})- \\ &\ \ \
-\sum_{i=1}^n(-1)^i1_{f_0\otimes\cdots\otimes f_n}
\circ\text\nn_{X_0,\ldots,X_{i-1}\otimes X_i,\ldots,X_n}+ \\ &\ \ \ +
\sum_{i=1}^n(-1)^i\text\nn_{X'_0,\ldots,X'_{i-1}\otimes X'_i,\ldots,X'_n}\circ
1_{f_0\otimes\cdots\otimes f_n}- \\ &\ \ \
-(-1)^{n+1}(1_{f_0\otimes\cdots\otimes f_{n-1}}
\circ\text\nn_{X_0,\ldots,X_{n-1}})\otimes 1_{f_n}+ \\ &\ \ \
+(-1)^{n+1}(\text\nn_{X'_0,\ldots,X'_{n-1}}\circ
1_{f_0\otimes\cdots\otimes f_{n-1}})\otimes 1_{f_n}
\end{align*}
(the reader may easily check that the padding operators appearing
in the definition of $\delta_v$ are indeed trivial in this case).
Now, since \Cg\ is a Gray semigroup, it is
$\widehat{\otimes}((f_0,f_1\otimes\cdots\otimes
f_n),(id_{X_0},id_{X_1\otimes\cdots\otimes
X_n}))=\widehat{\otimes}((id_{X_0},id_{X_1\otimes\cdots\otimes
X_n}),(f_0,f_1\otimes\cdots\otimes f_n))=1_{f_0\otimes\cdots\otimes
f_n}$. Therefore, using Equation $(A\widehat{\otimes}1)$, we get
\begin{eqnarray*}
&1_{f_0}\otimes(1_{f_1\otimes\cdots\otimes f_n}
\circ\text\nn_{X_1,\ldots,X_n})=1_{f_0\otimes\cdots\otimes f_n}
\circ(1_{id_{X_0}}\otimes\text\nn_{X_1,\ldots,X_n})
\\
&1_{f_0}\otimes(\text\nn_{X'_1,\ldots,X'_n}\circ
1_{f_1\otimes\cdots\otimes
f_n})=(1_{id_{X_0}}\otimes\text\nn_{X'_1,\ldots,X'_n})\circ
1_{f_0\otimes\cdots\otimes f_n}
\end{eqnarray*}
The last two terms are treated similarly. It follows immediately
that
$\delta_v(\varphi(\text\nn))=\varphi(\delta_{pent}(\text\nn))$, as
required.
\end{proof}

\begin{rem}
Notice that the pentagonator-deformation subcomplex
$X^{\bullet}_{pent}(\text\Cg)\subseteq\widetilde{X}^{\bullet}_{pent}(\text\Cg)$
is nothing more that the kernel of this morphism $\varphi$.
\end{rem}

Associated to this cochain map, there is the corresponding cone
complex, which will be denoted by
$(X^{\bullet}_{pent,ass}(\text\Cg),\delta_{pent,ass})$. By
definition
$$
X^n_{pent,ass}(\text\Cg)=X^{0,n}(\text\Cg)\oplus
\widetilde{X}^{n+1}_{pent}(\text\Cg),\ \ \ \ \ \ \ n\geq 0
$$
with coboundary map $\delta_{pent,ass}:X^n_{pent,ass}(\text\Cg)\To
X^{n+1}_{pent,ass}(\text\Cg)$ given by
$$
\delta_{pent,ass}(\phi,\text\nn)=(\delta_v(\phi)+\varphi(\text\nn),-\delta_{pent}(\text\nn))
$$
for all $(\phi,\text\nn)\in X^n_{pent,ass}(\text\Cg)$. We have the
following modified version of the extended double complex of
$(\text\Cg,\otimes)$:
\begin{prop}
Let's substitute the first column
$(X^{0,\bullet}(\text\Cg),\delta_v)$ of the extended double complex
$(X^{\bullet,\bullet}(\text\Cg),\delta_h,\delta_v)$ of
$(\text\Cg,\otimes)$ for the previous cone complex
$(X^{\bullet}_{pent,ass}(\text\Cg),\delta_{pent,ass})$, and the
coboundary maps $\delta_h:X^{0,n}(\text\Cg)\To X^{1,n}(\text\Cg)$,
$n\geq 0$, for the maps $\delta'_h:X^n_{pent,ass}(\text\Cg)\To
X^{1,n}(\text\Cg)$ given by the projection to $X^{0,n}(\text\Cg)$
followed by $\delta_h$. Then, the resulting collection of
$K$-vector spaces and linear maps is a double complex (see Fig.
~\ref{figura_complex_doble_modificat}).
\end{prop}

\begin{figure}
\centering
\input{mdouble.pstex_t}
\caption{The modified double complex of a Gray semigroup}
\label{figura_complex_doble_modificat}
\end{figure}

\begin{proof}
By the way the $\delta'_h$ are defined, it is clear that the new
rows are cochain complexes. So, we only need to see that the
squares on the left of Fig???? still commute. Let
$(\phi,\text\nn)\in X^{n-1}_{pent,ass}(\text\Cg)$. Since the
coboundary maps $\delta_h,\delta_v$ commute, we have
\begin{align*}
\delta'_h(\delta_{pen,ass}(\phi,\text\nn))&=
\delta'_h(\delta_v(\phi)+\varphi(\text\nn),-\delta_{pent}(\text\nn))
\\ &=\delta_h(\delta_v(\phi))+\delta_h(\varphi(\text\nn))
\\ &=\delta_v(\delta_h(\phi))+\delta_h(\varphi(\text\nn))
\\ &=\delta_v(\delta'_h(\phi,\text\nn))+\delta_h(\varphi(\text\nn))
\end{align*}
Therefore, the proof reduces to show that the term
$\delta_h(\varphi(\text\nn))$ is zero for all
$\text\nn\in\widetilde{X}^n_{pent}(\text\Cg)$. Now, by definition
of $\delta_h$ and $\varphi$, we have
\begin{align*}
\delta_h(\varphi(\text\nn))((f'_0,\ldots,f'_n),(f_0,\ldots,f_n))&=
\lceil 1_{f'_0\otimes\cdots\otimes f'_n}\circ(\varphi(\text\nn))(f_0,\ldots,f_n)\rceil
\\ &\ \ \ \ \ \ \ \ \ -\lceil(\varphi(\text\nn))(f'_0\circ f_0,\ldots,f'_n\circ f_n)\rceil
\\ &\ \ \ \ \ \ \ \ \ +\lceil(\varphi(\text\nn))(f'_0,\ldots,f'_n)\circ
1_{f_0\otimes\cdots\otimes f_n}\rceil
\\ &=-\lceil 1_{f'_0\otimes\cdots\otimes f'_n}\circ
(1_{f_0\otimes\cdots\otimes
f_n}\circ\text\nn_{X_0,\ldots,X_n})\rceil
\\ &\ \ \ \ \ \ \ \ \ +\lceil 1_{f'_0\otimes\cdots\otimes f'_n}\circ(\text\nn_{X'_0,\ldots,X'_n}\circ
1_{f_0\otimes\cdots\otimes f_n})\rceil
\\ &\ \ \ \ \ \ \ \ \ +\lceil 1_{(f'_0\circ f_0)\otimes\cdots\otimes (f'_n\circ f_n)}
\circ\text\nn_{X_0,\ldots,X_n}\rceil
\\ &\ \ \ \ \ \ \ \ \ -\text\nn_{X''_0,\ldots,X''_n}\circ
1_{(f'_0\circ f_0)\otimes\cdots\otimes (f'_n\circ f_n)}\rceil
\\ &\ \ \ \ \ \ \ \ \ -\lceil (1_{f'_0\otimes\cdots\otimes f'_n}\circ\text\nn_{X'_0,\ldots,X'_n})\circ
1_{f_0\otimes\cdots\otimes f_n}\rceil
\\ &\ \ \ \ \ \ \ \ \ +\lceil(\text\nn_{X''_0,\ldots,X''_n}\circ
1_{f'_0\otimes\cdots\otimes f'_n})\circ 1_{f_0\otimes\cdots\otimes
f_n}\rceil
\end{align*}
Since \Cg\ is a 2-category, the second and fifth terms clearly
cancel out each other. On the other hand, using again that \Cg\ is
a 2-category and the interchange law and the naturality of the
2-morphisms
$\widehat{\otimes(n+1)}((f'_0,\ldots,f'_n),(f_0,\ldots,f_n))$ in
$(f'_0,\ldots,f'_n),(f_0,\ldots,f_n)$, we have
\begin{align*}
\lceil 1_{f'_0\otimes\cdots\otimes f'_n}\circ&(1_{f_0\otimes\cdots\otimes
f_n}\circ\text\nn_{X_0,\ldots,X_n})\rceil=
\\ &=\widehat{\otimes(n+1)}((f'_0,\ldots,f'_n),(f_0,\ldots,f_n))\cdot
(1_{(f'_0\otimes\cdots\otimes f'_n)\circ(f_0\otimes\cdots\otimes
f_n)}\circ\text\nn_{X_0,\ldots,X_n})
\\ &=(\widehat{\otimes(n+1)}((f'_0,\ldots,f'_n),(f_0,\ldots,f_n))\circ
1_{id_{X_0\otimes\cdots\otimes X_n}})\cdot
(1_{(f'_0\otimes\cdots\otimes f'_n)\circ(f_0\otimes\cdots\otimes
f_n)}\circ\text\nn_{X_0,\ldots,X_n})
\\ &=(\widehat{\otimes(n+1)}((f'_0,\ldots,f'_n),(f_0,\ldots,f_n))\cdot
1_{(f'_0\otimes\cdots\otimes f'_n)\circ(f_0\otimes\cdots\otimes
f_n)})\circ \text\nn_{X_0,\ldots,X_n}
\\ &=(1_{(f'_0\circ f_0)\otimes\cdots\otimes (f'_n\circ f_n)}\cdot
\widehat{\otimes(n+1)}((f'_0,\ldots,f'_n),(f_0,\ldots,f_n)))\circ\text\nn_{X_0,\ldots,X_n}
\\ &=(1_{(f'_0\circ f_0)\otimes\cdots\otimes (f'_n\circ f_n)}\circ\text\nn_{X_0,\ldots,X_n})
\cdot\widehat{\otimes(n+1)}((f'_0,\ldots,f'_n),(f_0,\ldots,f_n)))
\\ &=\lceil 1_{(f'_0\circ f_0)\otimes\cdots\otimes (f'_n\circ f_n)}
\circ\text\nn_{X_0,\ldots,X_n} \rceil
\end{align*}
Therefore, the first and third term above also cancel out each
other. The same thing can be shown similarly for the fourth and
sixth terms.
\end{proof}

This new double complex will be called the {\sl modified extended
double complex} of $(\text\Cg,\otimes)$, and denoted by
$X^{\bullet,\bullet}_{mod}(\text\Cg)$. So
$$
X^{m,n}_{mod}(\text\Cg)=\left\{ \begin{array}{ll}
X^{0,n}(\text\Cg)\oplus\widetilde{X}^{n+1}_{pent}(\text\Cg) &
\text{if}\ m=0,n\geq 0 \\ X^{m,n}(\text\Cg) &
\text{if}\ m\geq 1,n\geq 0
\end{array}\right.
$$
For short, the corresponding horizontal and vertical coboundary
maps will be denoted by $\delta'_h$ and $\delta'_v$, respectively.
But remind that $\delta'_v=\delta_v$ except for the first column,
where $\delta'_v=\delta_{pent,ass}$, and that $\delta'_h=\delta_h$
except for $m=0$, where it is the projection to the first component
followed by $\delta_h$.

\subsection{}
Let us proceed now as in Section 8 and consider the double complex
obtained from $X^{\bullet,\bullet}_{mod}(\text\Cg)$ after deleting
the first row. It will be the {\sl modified double complex} of
$(\text\Cg,\otimes)$. Furthermore, let's take the subcomplex
$X^{\bullet,\bullet}_{mod,s}(\text\Cg)$ of this modified double
complex corresponding to the special elements. These are defined in
the same way as in Section 8 for all $m\geq 1,n\geq 1$, while
$(\phi,\text\nn)\in X^{0,n}_{mod}(\text\Cg)$, $n\geq 1$, is called
special whenever $\phi$ is special. We leave to the reader to check
that special elements are preserved by the coboundary maps
$\delta'_h,\delta'_v$, so that
$X^{\bullet,\bullet}_{mod,s}(\text\Cg)$ is indeed a double complex
(the {\sl modified special double complex} of
$(\text\Cg,\otimes)$). Then, we can consider the associated total
complex, which will be denoted by $X^{\bullet}_{unit}(\text\Cg)$,
and called the {\sl unitary deformation complex} of the Gray
semigroup $(\text\Cg,\otimes)$. By definition
$$
X^q_{unit}(\text\Cg)=\bigoplus_{\begin{array}{c} m+n=q
\\ m\geq 0,n\geq 1 \end{array}}X^{m,n}_{mod,s}(\text\Cg),
\ \ \ \ \ q\geq 1
$$
and the coboundary operator $\delta_{unit}:X^q_{unit}(\text\Cg)\To
X^{q+1}_{unit}(\text\Cg)$ is
$$
\delta_{unit}=\bigoplus_{\begin{array}{c} m+n=q \\ m\geq 0,n\geq 1 \end{array}}
((-1)^n\delta'_h+\delta'_v),\ \ \ \ \ q\geq 1
$$
If $H^{\bullet}_{unit}(\text\Cg)$ denote its cohomology groups, we
have then the following final theorem, which says that this is the
right cochain complex describing the generic unitary deformations:

\begin{thm}
Given a $K$-linear Gray semigroup $(\text\Cg,\otimes)$, the
equivalence classes of its first order unitary deformations are in
bijection with the elements of $H^2_{unit}(\text\Cg)$.
\end{thm}
\begin{proof}
Let's consider 2-isomorphisms of the form
\begin{align*}
\widehat{\otimes}_{\epsilon}((f',g'),(f,g))&=\widehat{\otimes}((f',g'),(f,g))+
\widehat{\otimes}^{(1)}((f',g'),(f,g))
\\ (\otimes_0)_{\epsilon}(X,Y)&=1_{id_{X\otimes Y}} \\
\widehat{a}_{\epsilon}(f,g,h)&=1_{f\otimes g\otimes h}+\widehat{a}^{(1)}(f,g,h)\epsilon \\
(\pi_{\epsilon})_{X,Y,Z,T}&=1_{id_{X\otimes Y\otimes Z\otimes
T}}+(\pi^{(1)})_{X,Y,Z,T}\epsilon
\end{align*}
with $\widehat{\otimes}^{(1)}((f',g'),(f,g)):(f'\otimes g')\circ
(f\otimes g)\Longrightarrow (f'\circ f)\otimes (g'\circ g)$,
$\widehat{a}^{(1)}(f,g,h):f\otimes g\otimes h\Longrightarrow
f\otimes g\otimes h$ and $\pi^{(1)}_{X,Y,Z,T}:id_{X\otimes Y\otimes
Z\otimes T}\Longrightarrow id_{X\otimes Y\otimes Z\otimes T}$. In
particular, the $\pi^{(1)}_{X,Y,Z,T}$ clearly define an element
$$
\pi^{(1)}\in\widetilde{X}^3_{pent}(\text\Cg)\subseteq X^{0,2}_s(\text\Cg)\oplus
\widetilde{X}^3_{pent}(\text\Cg)=X^{0,2}_{mod,s}(\text\Cg)\subseteq X^2_{unit}(\text\Cg)
$$
Then, applying Proposition~\ref{estructura_semigrupal}, it follows
that the above 2-isomorphisms define a semigroupal structure on
$\text\Cg_{(1)}^0$ (hence, a first order unitary deformation of
$\text\Cg$) if and only if:
\begin{description}

\item[A$\widehat{\otimes}$1]
The $\widehat{\otimes}^{(1)}((f',g'),(f,g))$ define an element
$\widehat{\otimes}^{(1)}\in
X^{1,1}(\text\Cg)=X^{1,1}_{mod}(\text\Cg)$.

\item[A$\widehat{\otimes}$2]
$\widehat{\otimes}^{(1)}$ is such that
$\delta_h(\widehat{\otimes}^{(1)})=0$.

\item[A$\widehat{\otimes}$3]
$\widehat{\otimes}^{(1)}$ is special.

\item[A$\widehat{a}$1]
The $\widehat{a}^{(1)}(f,g,h)$ define an element
$\widehat{a}^{(1)}\in X^{0,2}(\text\Cg)\subset
X^{0,2}(\text\Cg)\oplus\widetilde{X}^3_{pent}(\text\Cg)=X^{0,2}_{mod}(\text\Cg)$.

\item[A$\widehat{a}$2]
$\widehat{a}^{(1)}$ and $\widehat{\otimes}^{(1)}$ are such that
$\delta_h(\widehat{a}^{(1)})+\delta_v(\widehat{\otimes}^{(1)})=0$.

\item[A$\widehat{a}$3]
$\widehat{a}^{(1)}$ is special.

\item[A$\pi$1]
$\widehat{a}^{(1)}$ and $\pi^{(1)}$ are such that
$\delta_v(\widehat{a}^{(1)})+\varphi(\pi^{(1)})=0$.

\item[A$\pi$2]
$\pi^{(1)}$ is such that $\delta_{pent}(\pi^{(1)})=0$.

\end{description}

Now, $\pi^{(1)}\in\widetilde{X}^3_{pent}(\text\Cg)$ together with
$(A\widehat{\otimes} 1)$, $(A\widehat{\otimes} 3)$ and
$(A\widehat{a} 1)$, $(A\widehat{a} 3)$ say that
$$
((\widehat{a}^{(1)},\pi^{(1)}),\widehat{\otimes}^{(1)})\in
(X^{0,2}_s(\text\Cg)\oplus\widetilde{X}^3_{pent}(\text\Cg))\oplus
X^{1,1}_s(\text\Cg)=X^2_{unit}(\text\Cg)
$$
On the other hand, we have
$$
\delta_{unit}((\widehat{a}^{(1)},\pi^{(1)}),\widehat{\otimes}^{(1)})=
((\delta_v(\widehat{a}^{(1)})+\varphi(\pi^{(1)}),-\delta_{pent}(\pi^{(1)})),
\delta_h(\widehat{a}^{(1)})+\delta_v(\widehat{\otimes}^{(1)}),
-\delta_h(\widehat{\otimes}^{(1)}))
$$
Hence, $(A\pi 1)$, $(A\pi 2)$, $(A\widehat{a}2)$ and
$(A\widehat{\otimes}2)$ together exactly say that
$((\widehat{a}^{(1)},\pi^{(1)}),\widehat{\otimes}^{(1)})$ is a
2-cocycle.

Let's consider
$(((\widehat{a}^{(1)})',(\pi^{(1)})'),(\widehat{\otimes}^{(1)})')$
another 2-cocycle defining an equivalent first order unitary
deformation of
\Cg. Then, by Proposition~\ref{equivalents} applied to our situation,
there exists 2-morphisms $\widehat{\psi}^{(1)}(f,g):f\otimes
g\Longrightarrow f\otimes g$ and
$(\omega^{(1)})_{X,Y,Z}:id_{X\otimes Y\otimes Z}\Longrightarrow
id_{X\otimes Y\otimes Z}$ (in particular,
$\omega^{(1)}\in\widetilde{X}^2_{pent}(\text\Cg)$) such that
\begin{description}

\item[E$\widehat{\psi}$1]
The  $\widehat{\psi}^{(1)}(f,g)$ define an element
$\widehat{\psi}^{(1)}\in X^{0,1}(\text\Cg)$.

\item[E$\widehat{\psi}$2]
$(\widehat{\otimes}^{(1)})'-\widehat{\otimes}^{(1)}=-\delta_h(\widehat{\psi}^{(1)})$.

\item[E$\widehat{\psi}$3]
$\widehat{\psi}^{(1)}$ is special.

\item[E$\omega$1]
$(\widehat{a}^{(1)})'-\widehat{a}^{(1)}=\delta_v(\widehat{\psi}^{(1)})+\varphi(\omega^{(1)})$

\item[E$\omega$2]
$(\pi^{(1)})'-\pi^{(1)}=-\delta_{pent}(\omega^{(1)})$

\end{description}

The first and third conditions together say that
$(\widehat{\psi}^{(1)},\omega^{(1)})\in
X^{0,1}_s(\text\Cg)\oplus\widetilde{X}^2_{pent}(\text\Cg)=
X^{0,1}_{mod,s}(\text\Cg)=X^1_{unit}(\text\Cg)$. On the other hand,
the reader may check that
$$
\delta_{unit}(\widehat{\psi}^{(1)},\omega^{(1)})=
((\delta_v(\widehat{\psi}^{(1)})+\varphi(\omega^{(1)}),-\delta_{pent}(\omega^{(1)})),
-\delta_h(\widehat{\psi}^{(1)}))
$$
so that the remaining conditions just say that
$$
(((\widehat{a}^{(1)})',(\pi^{(1)})'),(\widehat{\otimes}^{(1)})')-
((\widehat{a}^{(1)},\pi^{(1)}),\widehat{\otimes}^{(1)})=
\delta_{unit}(\widehat{\psi}^{(1)},\omega^{(1)})
$$
Hence, both 2-cocycles are cohomologous, as required.
\end{proof}

\section{Concluding remarks}

The present work is not intended to give the complete picture of
the theory of infinitesimal deformations of a monoidal 2-category.
Indeed, various points still remain for future work. Among them,
let us mention the following:

\begin{description}
\item[1]
There is the all-important question of the higher-order
obstructions. This turned out to be the most difficult point in the
cohomological deformation theory for monoidal categories first
initiated by Crane and Yetter \cite{CY981} and further developped by
Yetter \cite{dY01}. We guess that our cohomological description also
fits nicely into
the general picture established by Gerstenhaber for a good
deformation theory \cite{mG64}. But as already mentioned, this is left for a future work.

\item[2]
As stated at the beginning of this work, the ultimate goal
should be to get a cohomological description of the infinitesimal
deformations of a {\it monoidal} 2-category. Hence, it also
deserves further work the question of how to take into account the
additional unital structure in the whole theory.
In the case of monoidal categories, Yetter
\cite{dY01} has shown that the deformations of this additional
structure are already determined by those of the semigroupal
structure. It seems possible that the same situation reproduces in
our framework.

\item[3]
Finally, in the present work we have restricted our attention to
those infinitesimal deformations of the semigroupal 2-category
$(\text\Cg,\otimes,a,\pi)$ such that the bicategory structure of
$\text\Cg$ remains undeformed. But, as pointed out previously, this
structure can also be deformed. As regards this point, notice that
the elements $\phi\in X^{2,0}(\text\Cg)$ of our extended double
complex are of the form
$$
\phi(h,g,f):h\circ g\circ f\Longrightarrow h\circ g\circ f
$$
(we are thinking of a Gray semigroup, so that parenthesis are not
needed here). This suggests that the possible deformations of the
bicategory structure of \Cg\ may be related to such elements
$\phi\in X^{2,0}(\text\Cg)$. In this sense, the situation can once
more resemble that encountered in the deformation theory of a
bialgebra. Indeed, it can be shown (see
\cite{GS92},\cite{sS92},\cite{SS93}) that the right cochain complex
describing the deformations of a bialgebra as a {\it
quasibialgebra} (i.e., coassociative only up to conjugation) is
precisely that associated to the double complex obtained after
adding the bottom row of the full Gerstenhaber-Schack complex,
which had to be deleted to study the deformations in the bialgebra
setting. In our case, the weakening of the coassociativity
condition should correspond to the weakening of the 2-category
condition $\alpha_{h,g,f}=1_{h\circ g\circ f}$. It seems possible,
then, that taking into account the deleted bottom row in our double
complex is just the only step needed to consider these more general
deformations of
\Cg.
\end{description}

Another important point not addressed in this paper, and which we
are currently working on, is the question of examples. In
particular, a simple example of a $K$-linear Gray semigroup where
our theory can be applied is that introduced by Mackaay
\cite{mM00}, denoted by {\bf N}$(G,H,K^*)$, and associated to a
pair of finite groups $G,H$ (with $H$ an abelian group); it
includes as a special case the (semistrict version of) the monoidal
2-category of 2-vector spaces. More interesting examples, however,
are expected to come from the 2-categories of representations of
the Hopf categories associated to quantum groups, whose construction was sketched by Crane and Frenkel \cite{CF94}.

Finally, let us finish by mentioning the interest our work may have
for homotopy theory. Indeed, since Grothendieck \cite{aG83}, it was
suspected that homotopy $n$-types were somewhat equivalent to
certain algebraic structures called {\it weak n-groupoids}, which
should be a particular kind of {\it weak n-categories}
characterized by the fact that all morphisms are invertible up to
suitable equivalence. Recently, Tamsamani \cite{zT99} realized this idea,
giving a precise definition of a weak $n$-groupoid, for
any non-negative integer $n$, together with a suitable notion of
equivalence, and showing that the equivalence classes of weak
$n$-groupoids bijectively correspond to homotopy classes of
$n$-anticonnected CW-complexes. Since weak 3-groupoids with one object should correspond to
a special type of semigroupal 2-categories, it naturally
raises the question about the meaning our cohomology theory has in
this topological setting. Such a relation between certain monoidal
2-categories and homotopy 3-types has recently been discussed by
Mackaay \cite{mM00}, who conjectures that the classification of
{\it semi-weak} monoidal 2-category structures on the above
mentioned 2-category {\bf N}$(G,H,K^*)$ boils down to the
classification up to homotopy equivalence of connected
3-anticonnected ($>$1-simple) CW-complexes $X$ with $\pi_1(X)=G$,
$\pi_2(X)=H$ and $\pi_3(X)=K^*$. Via Postnikov's theory, this leds
him to conjecture a bijection between the equivalence classes of
semi-weak monoidal 2-category structures on {\bf N}$(G,H,K^*)$ and
pairs of cohomology classes $\alpha\in H^3(B_G,H)$ and $\beta\in
H^4(W(\alpha),K^*)$ ($W(\alpha)$ denotes a certain CW-complex
induced by $\alpha$; see \cite{mM00}).

\vspace{1 truecm}

\noindent{{\it Acknowledgements.}} First of all I
want to thank Louis Crane (my coadvisor together with Sebasti\`{a} Xamb\'o) for suggesting me
the topic of this work for my PhD, and Louis Crane and David Yetter for the
helpful discussions with them. I also want to thank the Department
of Mathematics of Kansas State University (KSU) and Louis Crane
(and his wife) in particular for their hospitality during my stay
in Manhattan, where I started this work. Finally, I want to thank
the Departament d'Universitats, Recerca i Societat de la
Informaci\'o (DURSI) of the Generalitat the Catalunya and the
Universitat Polit\`{e}cnica de Catalunya (UPC) for their financial
support which enabled me to visit KSU.

\bibliographystyle{amsplain}
\bibliography{cdtm2c1}

\end{document}